\documentstyle{amltd}
\begin{document}
\input amssym.def
\input amssym.tex

\annalsline{151}{2000}
\received{May 13, 1998}
\startingpage{193}
\def\bye{\end{document}}
 \font\tenrm=cmr10


\catcode`\@=11
\font\twelvemsb=msbm10 scaled 1100
\font\tenmsb=msbm10
\font\ninemsb=msbm10 scaled 800
\newfam\msbfam
\textfont\msbfam=\twelvemsb  \scriptfont\msbfam=\ninemsb
  \scriptscriptfont\msbfam=\ninemsb
\def\msb@{\hexnumber@\msbfam}
\def\Bbb{\relax\ifmmode\let\next\Bbb@\else
 \def\next{\errmessage{Use \string\Bbb\space only in math
mode}}\fi\next}
\def\Bbb@#1{{\Bbb@@{#1}}}
\def\Bbb@@#1{\fam\msbfam#1}
\catcode`\@=12

 \catcode`\@=11
\font\twelveeuf=eufm10 scaled 1100
\font\teneuf=eufm10
\font\nineeuf=eufm7 scaled 1100
\newfam\euffam
\textfont\euffam=\twelveeuf  \scriptfont\euffam=\teneuf
  \scriptscriptfont\euffam=\nineeuf
\def\euf@{\hexnumber@\euffam}
\def\frak{\relax\ifmmode\let\next\frak@\else
 \def\next{\errmessage{Use \string\frak\space only in math
mode}}\fi\next}
\def\frak@#1{{\frak@@{#1}}}
\def\frak@@#1{\fam\euffam#1}
\catcode`\@=12

\def\overset#1#2{%
  {\mathop{\kern 0pt #2}\limits^{#1}}}
\def\underset#1#2{%
  {\mathop{\kern 0pt #2}\limits_{#1}}}

\newcommand{\fcirc}{\overset{\circ}{F}}

\def\rank{\mathop{\rm rank}\nolimits} 
\def\ind{\mathop{\rm ind}\nolimits} 
\def\Sing{\mathop{\rm Sing}\nolimits}
\def\supp{\mathop{\rm supp}\nolimits} 
\def\Aut{\mathop{\rm Aut}\nolimits} 
\def\Map{\mathop{\rm Map}\nolimits} 
\def\Reg{\mathop{\rm Reg}\nolimits} 
\def\ev{\mathop{\rm ev}\nolimits} 
\def\Fom{\mathop{\rm Hom}\nolimits} 
\def\Res{\mathop{\rm Res}\nolimits} 
\def\Euler{\mathop{\rm Euler}\nolimits} 
\def\spin{\mathop{\rm spin}\nolimits} 
\def\SL{\mathop{\rm SL}\nolimits} 
\def\Lie{\mathop{\rm Lie}\nolimits} 
\def\End{\mathop{\rm End}\nolimits} 
\def\tr{\mathop{\rm tr}\nolimits} 
\def\Inj{\mathop{\rm Inj}\nolimits} 
\def\Rad{\mathop{\rm Rad}\nolimits} 
\def\RRe{\mathop{\rm Re}\nolimits} 
\def\divv{\mathop{\rm div}\nolimits} 
\def\sup{\mathop{\rm sup}\limits} 
\def\curl{\mathop{\rm curl}\limits} 

\newcommand{\grad}{\mathop{\rm grad\,}\nolimits}
\newcommand{\vol}{\mathop{\rm vol\,}\nolimits}

\def\mathcal#1{{\cal #1}}
\def\mathbb#1{{\Bbb #1}}
\def\A{{\mathcal A}}
\def\B{{\mathcal B}}
\def\C{{\mathcal C}}
\def\E{{\mathcal E}}
\def\F{{\mathcal F}}
\def\G{{\mathcal G}}
\def\K{{\mathcal K}}
\def\L{{\mathcal L}}
\def\M{{\mathcal M}}
\newcommand{\R}{\mathcal R}
\def\U{{\mathcal U}}
\def\BB{{\mathbb B}}
\def\CC{{\mathbb C}}
\def\QQ{{\mathbb Q}}
\def\RR{{\Bbb R}}
\def\ZZ{{\mathbb Z}}
\def\mathfrak{\frak}
\def\sl{\mathfrak{sl}}

\makeatletter
\@addtoreset{equation}{section}
\makeatother
\def\spn{\speqnu}

 \title{Gauge theory and calibrated geometry, I}
\author{Gang Tian}
\institutions{Massachusetts Institute of Technology,
 Cambridge, MA\\
{\eightpoint {\it E-mail}\/: tian@math.mit.edu}}

\smallbreak \centerline{\bf Contents}
\smallbreak
\def\bri{\vskip1pt\noindent\hglue15pt}
\def\smb#1{\vskip4pt\noindent{\bf #1}}
\bri 0.1.\ Introduction
\smb{1.  Preliminaries}
\bri 1.1.\ The Yang-Mills functional
\bri 1.2.\ Anti-self-dual instantons
\bri 1.3.\ Complex anti-self-dual instantons
\bri 1.4.\ Instantons on $G_2$-manifolds
\smb{2. Consequences of a monotonicity formula}
\bri 2.1.\ A monotonicity formula
\bri 2.2.\ Curvature estimates
\bri 2.3.\ Admissible Yang-Mills connections
\smb{3. Rectifiability of blow-up loci}
\bri 3.1.\ Convergence of Yang-Mills connections
\bri 3.2.\ Tangent cones of blow-up loci
\bri 3.3.\ Rectifiability
\smb{4. Structure of blow-up loci}
\bri 4.1.\ Bubbling Yang-Mills connections
\bri 4.2.\ Blow-up loci of anti-self-dual instantons
\bri 4.3.\ Calibrated geometry and blow-up loci
\bri 4.4.\ Cayley cycles and complex anti-self-dual instantons
\bri 4.5.\ General blow-up loci

\smb{5. Removable singularities of Yang-Mills equations}
\bri 5.1.\ Stationary properties of Yang-Mills connections
\bri 5.2.\ A removable singularity theorem
\bri 5.3.\ Cone-like Yang-Mills connections

\smb{6. Compactification of moduli spaces}
\bri 6.1.\ Compactifying moduli spaces
\bri 6.2.\ Final remarks

\centerline{\bf \rm 0.1. \ Introduction}
\bigbreak
 
The geometry of submanifolds is intimately related to the theory of
functions and vector bundles. It has been of fundamental
importance to find out how those two objects interact in many
geometric and physical problems. A typical example of this
relation is that the Picard group of line bundles on an algebraic
manifold is isomorphic to the group of divisors, which is
generated by holomorphic hypersurfaces modulo linear equivalence.
A similar correspondence can be made between the K-group of
sheaves and the Chow ring of holomorphic cycles. There are two more 
very recent examples of such a relation. The mirror symmetry in string
theory has revealed a deeper phenomenon  involving
special Lagrangian cycles (cf.\ [SYZ]). On the other hand, C.
Taubes has shown that the Seiberg-Witten invariant coincides with
the Gromov-Witten invariant on any symplectic $4$-manifolds.

In this paper, we will show another natural interaction between
Yang-Mills connections, which are critical points of
a Yang-Mills action associated to a vector bundle, 
and minimal submanifolds, which have been studied extensively
for years in classical differential geometry and the
calculus of variations. 

Let $M$ be a manifold with a Riemannian metric $g$.
Let $E$ be a vector bundle over $M$ with a compact Lie
group as its structure group. For instance, 
$E$ may be a complex bundle and $G$ is then
a unitary group. A connection $A$ of $E$ can be given by specifying
a covariant derivative 

$$
D _A: C^\infty (E) \mapsto C^\infty(E\otimes \Omega^1M).
$$
\vglue4pt\noindent 
In local trivializations of $E$, $D_A$ is of the form 
$d + a$ for some ${\rm Lie}(G)$-valued $1$-form $a$. The curvature
of $A$ is a ${\rm Lie}(G)$-valued $2$-form $F_A$, which is equal to $D_A^2$.
As usual, it measures deviation from the symmetry of second derivatives.  
Such a connection $A$ is Yang-Mills if $D^*_A F_A=0$,
where $D^*_A$ is the adjoint of $D_A$ with respect to
the metric $g$. By the second
Bianchi identity, we also have $D_AF_A=0$. 
The system $D^*_A F_A=0, D_AF_A=0$ is called the
Yang-Mills equation and is invariant under so-called gauge transformations,
which are locally made of $G$-valued functions.

The moduli space of Yang-Mills connections is the quotient
of the set of solutions of the Yang-Mills equation by the gauge group, 
which consists of all gauge transformations.
It is well-known that this moduli space may not be compact. 
Given any sequence of Yang-Mills connections $\{A_i\}$ with 
a uniformly bounded $L^2$-norm of curvature, 
Uhlenbeck (also see [Na]) proved that by taking a
subsequence if necessary, $A_i$ converges to, modulo gauge transformations, 
a Yang-Mills connection $A$ in smooth topology outside a closed subset $S_b(\{A_i\})$
of Hausdorff codimension at least $4$. In fact, for any compact $K\subset M$,
$S_b(\{A_i\})\cap K$ has finite $(n-4)$-dimensional Hausdorff measure.
Furthermore, by taking subsequences if necessary,
we may assume that as measures, 
$|F_{A_i}|^2 dV_g$ converges weakly to $|F_A|^2dV_g
+ \Theta H^{n-4}\lfloor S_b(\{A_i\})$, where $\Theta \ge 0$ is a function
and  is called the  multiplicity of $S_b(\{A_i\})$, and $H^{n-4}\lfloor S_b(\{A_i\})$ is the 
$(n-4)$-dimensional Hausdorff measure restricted to $S_b(\{A_i\})$.
The set $S_b(\{A_i\})$ is the union of two closed subsets $S_b$ and
$S([A])$, where $S([A])$ consists of all points in $M$ where the 
$(n-4)$-dimensional density of $|F_A|^2dV_g$ is positive, and $S_b$ is the closure
of $S_b(\{A_i\})\backslash S([A])$. One can show that $\Theta H^{n-4}\lfloor S_b(\{A_i\})$ 
coincides with $\Theta H^{n-4}\lfloor S_b$ and $S([A])$ has vanishing $(n-4)$-dimensional
Hausdorff measure. Presumably, $S([A])$ is the singular set of $A$ modulo
gauge transformations. We will call $S_b$ with multiplicity $\Theta$ the blow-up locus 
of $\{A_i\}$ converging to $A$. If $M$ is
a $4$-dimensional compact manifold, the blow-up locus $S_b$
consists of finitely many points, $S([A])=\emptyset$ and the limiting connection
$A$ can be extended to be a Yang-Mills connection on the whole manifold 
with smaller $L^2$-norm of curvature [Uh1]. In particular, it follows
that the moduli space of anti-self-dual instantons on a $4$-manifold
(see the following
for the definition) can be compactified by adding all smaller 
anti-self-dual instantons together with finitely many points on $M$. 
This compactified moduli space plays
a fundamental role in the theory of Donaldson invariants.

With $M$   of higher dimension,  
little has been  known   about the blow-up locus $S_b$ itself.
Without further knowledge on the structure of $S_b$, one can not achieve
a reasonable compactification of the moduli space of Yang-Mills connections 
as we had in the case of $4$-manifolds. 
The main theme of this paper is to show that blow-up loci
of Yang-Mills connections have natural geometric structures and
introduce a natural compactification for moduli space of 
anti-self-dual instantons on higher dimensional manifolds by
adding cycles with appropriate geometric structure. We believe that
such a compactification will play an important role in our searching for
new invariants of Donaldson type for higher dimensional manifolds.

In this paper, we will first show that {\it any blow-up locus
$S_b$ is rectifiable}\/; {\it i{\rm .}e{\rm .,} except for a subset of $(n-4)$\/{\rm -}\/dimensional
Hausdorff measure zero{\rm ,} it is contained in a countable union
of $C^1$\/{\rm -}\/smooth submanifolds of dimension $n-4$}
(cf.\ Proposition 3.3.3). It is equivalent to saying that
{\it $S_b$ has a unique tangent space $T_xS_b$ for} 
$H^{n-4}$-a.e.\ $x$ {\it in} $S_b$. It can be thought of  as a 
rough regularity for $S_b$. We will show that $S_b$
inherits a nice geometric structure (Chapter 4). 
We will also prove a removable singularity theorem for 
the limiting Yang-Mills connection $A$ (Chapter 5).
It follows that $A$ can be extended smoothly to the complement
of $S([A])$ modulo gauge transformations.

Let $\Omega$ be a closed differential 
form of degree $n-4$ on $M$. Then one can define a linear operator 
$T=-\ast \Omega \wedge $ acting
on $2$-forms, where $\ast$ denotes the Hodge operator of the metric $g$. 
A connection $A$ is $\Omega$-anti-self-dual if its curvature
form $F_A$ is annihilated by $T - {\rm Id}$. One can also define 
the $\Omega$-anti-self-duality for more general connections (cf.\ Section 1.2).
We observe that the $\Omega$-anti-self-duality implies the Yang-Mills equation
and is invariant under gauge transformations. Furthermore, if
$M$ is a compact manifold without boundary and $A$ is $\Omega$-anti-self-dual, 
there is an {\it a priori}  $L^2$ bound on
$F_A$, which depends only on $E$, $M$ and $\Omega$. 

We will prove that 
{\it if $\{A_i\}$ is a sequence of $\Omega$-anti-self-dual
instantons converging to $A$ with blow-up locus $S_b$ with multiplicity
$ \Theta$, then $(S_b, \Theta)$ defines a closed integral current calibrated by 
$\Omega$} (Theorem 4.2.3). In particular, $\Theta$ is integer-valued
and $\Omega$ restricts to the induced volume form on each tangent 
space $T_xS_b$. 
If $\Omega$ has co-mass one, then  this implies 
that the blow-up locus $(S_b, \Theta)$ is area-minimizing (cf.\ [HL]).
Known regularity theorems in geometry measure theory further imply that $S$ is the closure of
a smooth submanifold calibrated by $\Omega$.
We will also prove a removable singularity theorem for any stationary
Yang-Mills connections (Theorem 5.2.1). Particularly, this implies that 
{\it the limiting connection $A$ 
extends to become a smooth connection on $M\backslash S$ for a
closed subset $S$ with vanishing $(n-4)$-dimensional
Hausdorff measure $H^{n-4}(S)=0$} (Theorem 5.2.2).

Now we can introduce a natural compactification of the moduli space 
${\cal M}_{\Omega, E}$ of $\Omega$-anti-self-dual instantons of $E$ on $M$.

A generalized $\Omega$-anti-self-dual instanton is a pair 
$(A, C)$ satisfying: (1) $A$ is $\Omega$-anti-self-dual on
$M\backslash S(A)$ with $(n-4)$-dimensional
Hausdorff measure $H^{n-4}(S(A)) =0$; (2) $C=(S,\Theta)$ is a closed,
integral current calibrated by $\Omega$; (3) The second Chern class
$C_2(E)$ of $E$ is the same as $[C_2(A)] + [C_2(S, \Theta)]$, where
$C_2(A)$ denotes the second Chern form of $A$ and $[C_2(S, \Theta)]$
denotes the Poincar\'e dual of the homology class represented by the current
$(S, \Theta)$. If the co-norm $|\Omega| \le 1$, it follows from a result of
F. Almgren that $C$ is of the form 
$\sum_{a=1}^{l(C)} m_a C_a$ ($l(C)$ may be zero),
such that each $m_a$ is a positive integer and
$C_a$ is the closure of a submanifold calibrated by $\Omega$. 

Two generalized $\Omega$-anti-self-dual instantons
$(A, C)$, $(A', C')$ are equivalent if and only if
$C=C'$ and there is a gauge transformation
$\sigma$ such that $\sigma(A)= A'$ on $M\backslash S(A) \cup S(A')$. 
We denote by $[A, C]$ the equivalence class represented by
$(A,C)$. Clearly, $[A, C]\in {\cal M}_{\Omega, E}$
if and only if $C=0$ and $A$ extends smoothly to $M$ modulo
a gauge transformation.

We define $\overline {\cal M}_{\Omega, E}$ to be the set of all
equivalence classes of generalized $\Omega$-anti-self-dual instantons of $E$.

The topology of $\overline {\cal M}_{\Omega, E}$ can be defined as follows:
a sequence $[A_i, C_i]$ converges to $[A,C]$ in 
$\overline {\cal M}_{\Omega, E}$ if and only if (1) 
$C_i$ converges to a closed integral current $C_\infty \subset C$ with respect to 
the weak topology for currents; (2) 
There are gauge transformations $\sigma_i$ such that $\sigma_i(A_i)$
converges to $A$ outside $S(A)\cup (C\backslash C_\infty)$.
One can show that this topology makes $\overline {\cal M}_{\Omega, E}$
a Hausdorff space. 

\vglue-2pt
 It follows from results in Chapters 4 and 5 that
{\it $\overline {\cal M}_{\Omega, E}$ is compact with respect  to
this topology on any compact manifold $M$} (Theorem 6.1.1).

Clearly,  $\overline {\cal M}_{\Omega, E}$ coincides with Uhlenbeck's
compactification of the moduli space of anti-self-dual instantons on a 
$4$-manifold $M$.

There are two important cases of such $\Omega$-anti-self-dual instantons,
which are worth being mentioned. In the first case,
let $(M,\omega)$ be a complex $m$-dimensional
K\"ahler manifold with the K\"ahler
form $\omega$. For any connection $A$, 
its curvature $F_A$ decomposes into (2,0), (1,1) and (0,2)-parts
$F_A^{2,0}$, $F_A^{1,1}$ and $F_A^{0,2}$.
Put $\Omega = \frac{\omega^{m-2}}{(m-2)!}$.
Then $A$ is $\Omega$-anti-self-dual if and only if
$F_A^{0,2}=0$ and $F^{1,1}_A\cdot \omega=0$; i.e., 
$A$ is a Hermitian-Yang-Mills connection. 
Combining Theorem 4.2.3 with a result
of King or Harvey and Shiffman, we obtain that
{\it blow-up loci of Hermitian-Yang-Mills
connections are effective holomorphic integral cycles consisting of
complex subvarieties of codimension two} (Theorem 4.3.3). Consequently,
the compactification $\overline {\cal M}_{\frac{\omega^{m-2}}{(m-2)!}, E}$
is the collection of equivalence classes $[A,C]$,
where $A$ is a Hermitian-Yang-Mills connection and $C$
is a holomorphic integral cycle of complex 
dimension $m-2$. A holomorphic integral cycle
is a formal sum of irreducible subvarities with positive coefficients.
In view of the Donaldson-Uhlenbeck-Yau theorem that
each (irreducible) Hermitian-Yang-Mills connection corresponds to
a stable bundle, our generalized Hermitian-Yang-Mills connection
$[A,C]$ should correspond to a stable sheaf. 
We would like to point out that our method can
be applied to more general situations where the connections
are not necessarily Hermitian-Yang-Mills. In order to
conclude the holomorphic property of the blow-up
locus, we only need that the $(0,2)$-part of curvature
is much smaller compared to the full curvature tensor
during the limiting process.

One of our motivations in this work is to 
carry out part of the program proposed in [DT] in a rigorous way. 
The program is to build up a gauge theory in higher dimensions.
If one is less ambitious, one may just want to
construct new holomorphic invariants for
Calabi-Yau $4$-folds in terms of complex anti-self-dual instantons.
Complex anti-self-dual instantons are anti-self-dual with respect to
appropriate $4$-form $\Omega$ on $M$. Since they have been discussed before by Donaldson and 
Thomas, we refer the readers to [DT] and its references. In contrast to the previous case, we can prove
that {\it blow-up loci of complex anti-self-dual instantons
are Cayley cycles} (cf.\ Theorem 4.4.3). 
A Cayley cycle is a rectifiable set such
that its tangent spaces are Cayley with respect to the given
K\"ahler form and the  holomorphic (4,0)-form on the underlying
Calabi-Yau $4$-fold (cf.\ [HL]). Notice that
special Lagrangian submanifolds used in
[SYZ] are special cases of Cayley cycles.
This allows us to compactify
the moduli space of complex anti-self-dual instantons
in terms of Cayley cycles as we did in the above.
Our methods may also  be  used to produce Cayley cycles,
which seem to be elusive with our existing knowledge.

One implication of our results here is that minimal submanifolds
can be considered as limiting solutions of the Yang-Mills equation.
Bearing this in mind, we may expect to construct Yang-Mills connections 
from minimal submanifolds in general position. Indeed, near a  minimal submanifold,
one can construct approximated solutions of the Yang-Mills equation,
whose curvature concentrates near the submanifold,
in a suitable sense.

An outline of this paper is as follows:
In Chapter 1, we give general discussions on Yang-Mills
connections, particularly, $\Omega$-anti-self-dual instantons.
We   analyze the $\Omega$-anti-self-duality in a few important cases.
In Chapter 2, we will derive a slight generalization of the mononicity formula 
of P. Price, a basic curvature estimate of K. Uhlenbeck. Then we apply Uhlenbeck's 
estimate to defining Chern-Weil forms for admissible Yang-Mills connections,
which are kinds of singular connections.
In Chapter 3, we prove rectifiability of blow-up loci. In Chapter 4,
we prove that blow-up loci of anti-self-dual instantons are calibrated, closed
integral currents. We will also analyze a few important special cases.
Chapter 5 contains a new removable singularity theorem.
In  the last chapter, we discuss compactification of moduli space
of anti-self-dual instantons and some related problems. 

All the results of this paper can be generalized to the case
of the Yang-Mills-Higgs equation. The details will appear elsewhere.

The author would like to thank I. Singer for bringing 
Cayley submanifolds to his attention. I also thank J. Cheeger and T. Colding
for some useful conversations. Part of this work was done when the author
was visiting ETH in Z\"urich, Switzerland and the  Institute for Advanced Study, Princeton.
The author is grateful to both places for providing excellent research environments.
 
\section{Preliminaries}
 
\demo{{\rm 1.1.}\ The Yang-Mills functional}
 Let $\pi~\hbox{:~}E\to M$ be a vector bundle of rank $r$ over a
differentiable manifold $M$ with a Lie group $G$ as its structure
group.  Then there is an open covering ${U_{\alpha}}$ of $M$, such
that for each $\alpha$, there is a local trivialization
\begin{eqnarray}
\label{eq:1.1.1}
\pi^{-1}(U_{\alpha})&\stackrel{\varphi_{\alpha}}{\longrightarrow}&U_{\alpha}
       \times\RR^r \spn{1.1.1}\\[1ex]
\pi\downarrow&\quad &\downarrow p_1\nonumber\\[1ex]
U_{\alpha}&\stackrel{=}{\longrightarrow} &U_{\alpha}\nonumber
\end{eqnarray}
where $p_1$ is the projection onto the first factor.  Note that
each $\varphi_\alpha$ is a diffeomorphism.  Furthermore, if
$U_{\alpha}\cap U_{\beta}\not=\emptyset $, then one can write 
\begin{eqnarray}
\label{eq:1.1.2}
\varphi_{\alpha}\cdot\varphi_{\beta}^{-1}: (U_{\alpha}\cap
U_{\beta})\times\RR^r &\longrightarrow& (U_{\alpha}\cap
U_{\beta})\times\RR^r,\spn{1.1.2}\\[2ex]
(x,\upsilon) &\longrightarrow& (x,g_{\alpha\beta}(x)\upsilon) \nonumber
\end{eqnarray} 
for some function $g_{\alpha\beta}:U_{\alpha}\cap U_{\beta}\to
G\subset {\rm GL}(r,\RR)$.  Such a function $g_{\alpha\beta}$ is  
called a transition function of $\pi:E\to M$.

Examples we often use in this paper include complex vector
bundles with a hermitian structure. For those bundles, the
structure group $G$ is $U(r/2)$.

A connection $A$ on $E$ is defined by specifying a covariant
derivative
$$
D=D_A : C^{\infty}(E)\to C^{\infty}(E\otimes\Omega^1 M) .
$$
Here $C^\infty(E)$ denotes the space of $C^\infty$ sections of
the bundle $E$. In a local trivialization $(U_\alpha,
\varphi_\alpha)$ of $E$, the covariant derivative takes the form
\begin{equation}
\label{eq:1.1.3}
D=d+A_\alpha, ~~A_\alpha: U_\alpha\to T^\ast U_\alpha\otimes  {\rm Lie}(G), \spn{1.1.3}
\end{equation} 
where ${\rm Lie}(G)$ denotes the Lie algebra of the structure group $G$. 
If $G$ is a  unitary group, Condition \ref{eq:1.1.3} is
equivalent to saying that $D$ preserves the corresponding
hermitian structure of $E$.

Note that $A_\alpha$ usually has no global description on $M$.
If $(U_\beta,\varphi_\beta)$ is another local trivialization and
$g_{\alpha\beta}$ is the corresponding transition function, then
\begin{equation}
\label{eq:1.1.4}
A_\beta=g_{\alpha\beta}^{-1}dg_{\alpha\beta}+g_{\alpha\beta}^{-1}A_\alpha
g_{\alpha\beta}. \spn{1.1.4}
\end{equation}

The curvature of the connection $A$ is determined by
$D^2:\Omega^0(E)\to\Omega^2(E)$.  It is a tensor, usually denoted
by $F_A$ or simply $F$ if no confusion occurs.  Formally, the
curvature tensor $F_A$ can be written as
$$
F_A=dA+A\wedge A,
$$
which actually means that in each local trivialization
$(U_\alpha,\varphi_\alpha)$,
\begin{equation}
\label{eq:1.1.5}
F_\alpha=dA_\alpha+A_\alpha\wedge A_\alpha. \spn{1.1.5}
\end{equation}
If $\{x_1,\cdots, x_n\}$ is a local coordinate system for $U_\alpha$, then we
have
\begin{equation}
\label{eq:1.1.6}
A_\alpha=A_{\alpha,i}dx_i,\quad A_{\alpha,i}\in {\rm Lie}(G), \spn{1.1.6}
\end{equation}
and
\begin{eqnarray}
\label{eq:1.1.7}
F_\alpha&=&\frac{1}{2}\sum_{i,j}F_{\alpha , ij}dx_i\wedge dx_j, \spn{1.1.7}\\[2ex]
F_{\alpha , ij}&=&\frac{\partial A_{\alpha,j}}{\partial
  x_i}-\frac{\partial A_{\alpha,i}}{\partial x_j}+[A_{\alpha i},
A_{\alpha j}].\nonumber
\end{eqnarray}
It follows that
\begin{equation}
\label{eq:1.1.8}
F_\beta  =g_{\alpha\beta}^{-1}F_\alpha g_{\alpha\beta} . \spn{1.1.8}
\end{equation}
Hence, $F_A\in\Omega^2(\End(E))$.

From now on, we assume that $G$ is a compact Lie group. 
We denote by
$\langle \cdot,\cdot \rangle$ the Killing form of its Lie algebra ${\rm Lie}(G)$.  
If $G=U(r/2)$, we have $\phantom{\displaystyle\sum_\int}$
\begin{equation}
\label{eq:1.1.9}
 \langle a, b\rangle=-\tr(ab) \, , \qquad a,b\in u\left(r/2\right)
=\Lie\left(U(r/2)\right).  \spn{1.1.9}\\
\end{equation}
\vglue6pt\noindent 
We can easily extend $\langle \cdot,\cdot \rangle$ to a product on differential
forms with values in ${\rm Lie}(G)$ as follows: if $\phi$ and $\psi$ are
differential forms of degree $p$ and $q$, respectively, we define
$$
\langle \phi, \psi \rangle = \sum _{i_1,\cdots, i_p, j_1, \cdots ,j_q} 
\langle \phi_{i_1\cdots i_p}, \psi_{j_1\cdots j_q} \rangle dx_{i_1}\wedge \cdots \wedge 
dx_{i_p}\wedge dx_{j_1}\wedge \cdots \wedge dx_{j_q},
$$
where 
\begin{eqnarray}
\phi & =& \sum _{i_1, \cdots,
i_p} \phi_{i_1\cdots i_p} dx_{i_1}\wedge \cdots \wedge 
dx_{i_p}, ~\phi_{i_1\cdots i_p} \in {\rm Lie}(G) , \nonumber \\[2ex]
\psi &=& \sum _{j_1, \cdots ,j_q} 
 \psi_{j_1\cdots j_q} dx_{j_1}\wedge \cdots \wedge dx_{j_q},~
\psi_{j_1\cdots j_q} \in {\rm Lie}(G).\nonumber
\end{eqnarray}

Let us also fix a Riemannian metric $g$ on $M$ and denote by $dV$
its volume form.  Then we can define
$$
|F_A|^2=\sum_{i,j,k,l}\langle F_{\alpha ij},F_{\alpha kl} \rangle g^{ik}g^{jl}
$$
in terms of local trivializations, where $(g_{ij})$ is the metric
tensor of $g$ in $x_1,\ldots,x_n$ and $(g^{ij})$ is its inverse
matrix.

The Yang-Mills functional of $E$ is defined by
\begin{equation}
\label{eq:1.1.10}
YM(A)=\frac{1}{4\pi^2}\int_M |F_A|^2 dV_g. \spn{1.1.10}
\end{equation}

Let $\G$ be the gauge group of $E$, which consists of all smooth sections
of the bundle $P(E)\times_{\rm Ad} G$ associated to the adjoint representation 
${\rm Ad}$ of $G$, where $P(E)$ denotes the principal bundle of $E$.
In terms of those trivializations $\{U_\alpha, \varphi_\alpha \}$, 
any $\sigma$ in $\G$ is given by a family of $G$-valued functions $\sigma_\alpha$ satisfying: 
$$\sigma_\alpha = g_{\alpha \beta} \cdot \sigma _\beta \cdot
g_{\alpha \beta}^{-1} \enspace{\rm on}\enspace U_\alpha\cap U_\beta.$$
Let $\sigma (A)$ be the connection with
$D_{\sigma (A)}=\sigma\cdot D_A \cdot\sigma^{-1}$; i.e.,
in each $U_\alpha$, 
$$D_{\sigma (A)}= d - d\sigma _\alpha \cdot \sigma^{-1}_\alpha 
+ \sigma _\alpha \cdot
A_\alpha \cdot \sigma _\alpha^{-1}.$$
Two smooth connections $A_1$ and $A_2$ of $E$ are equivalent
if there is a gauge transformation $\sigma $ such that $A_2 = \sigma (A_1)$.
A simple observation is: if there is a gauge transformation $\tau$ of $E$
over an open-dense subset $U$ such that $A_2 = \tau (A_1)$ in $U$, then
$\tau$ extends to $M$ and $A_1$, $A_2$ are equivalent.

One can easily show
\begin{equation}
\label{eq:1.1.11}
YM(\sigma(A))=YM(A), \spn{1.1.11}
\end{equation}
where $\sigma (A)$ is the connection with
$\sigma(D_A)=\sigma\cdot D_A \cdot\sigma^{-1}$.

The Euler-Lagrange equation of $YM$ is
\begin{equation}
\label{eq:1.1.12}
D_A^\ast F_A=0, \spn{1.1.12}
\end{equation}
where $D_A^\ast$ denotes the adjoint operator of $D_A$ with respect to
the Killing form of $G$ and the Riemannian metric $g$ on $M$.  On
the other hand, by the second Bianchi identity, we have
\begin{equation}
\label{eq:1.1.13}
D_AF_A=0. \spn{1.1.13}
\end{equation}
This, together with (\ref{eq:1.1.12}), implies that if $A$ is a
critical point of $YM$, then $F_A$ is harmonic.  In this case, we
say the $A$ is a Yang-Mills connection.  It follows from (1.11)
that if $A$ is a
Yang-Mills connection, so is $\sigma(A)$ for any gauge
transformation $\sigma$. In other words, both equations
(\ref{eq:1.1.12}) and (\ref{eq:1.1.13}) are invariant under the
action of the gauge group.
\enddemo

\vglue6pt
\demo{{\rm 1.2.}\ Anti-self-dual instantons}
 In this section, we discuss a special class of solutions 
to the Yang-Mills equation (1.1.12).
This class includes Hermitian-Yang-Mills connections on
a K\"{a}hler manifold.

Let $\pi : E\mapsto M$ be a unitary bundle of complex rank $r$,
and $\Omega $ be a closed
form of degree $n-4$, where $n=\dim M$. 
As before, we fix a Riemannian metric $ g$ on $M$.
We denote by $\ast $ the Hodge operator acting on forms with values
in ${\rm Lie}(G)$; i.e., for any $\phi$, $\psi$ in $\Omega ^p(\Lie (G))$,
$\ast \psi \in \Omega ^{n-p}({\rm Lie}(G))$ and 
\begin{equation}
\label{eq:1.2.1}
\langle \phi\wedge \ast \psi \rangle = ( \phi, \psi )  dV_g, \spn{1.2.1}
\end{equation}
where $( \cdot, \cdot )$ denotes the 
inner product on $\Omega ^p({\rm Lie} (G))$ induced by
$g$ and the Killing form $\langle \cdot, \cdot\rangle $.

Let $\tr$ be the standard trace on unitary matrices. 
For any unitary connection $A$ of the bundle $E$ over $M$,
we have a well-defined $\tr(F_A)$ in $\Omega ^2(M)$. 
It follows from the second
Bianchi identity that $\tr (F_A)$ is in fact a closed 
$2$-form.
In fact, $\frac{\sqrt{-1}}{2\pi} \tr (F_A)$
represents the first Chern class $C_1(E)$ in $H^2(M, \RR)$.
\enddemo

\vglue6pt

\specialnumber{1.2.1} \proclaim{Lemma}
\label{lemma:1.2.1}
Let $A$ be a unitary connection of $E$ over $M$ such that
$\tr(F_A)$ is a harmonic $2$\/{\rm -}\/form and 
\begin{equation}
\label{eq:1.2.2}
\Omega \wedge (F_A - \frac{1}{r}\tr(F_A) {\rm Id} )
= - \ast (F_A- \frac{1}{r} \tr(F_A) {\rm Id}), \spn{1.2.2}
\end{equation}
then $A$ is a Yang\/{\rm -}\/Mills connection{\rm .} Moreover{\rm , }
if $M$ is a compact manifold without boundary{\rm ,} 

\phantom{top of page}
\vglue-24pt
\begin{eqnarray}
\label{eq:1.2.3}
&& \hskip-7pt\frac{1}{4\pi^2}
\int _M |F_A|^2 dV_g - \frac{1}{4 r \pi^2}\int _M |\tr(F_A)|^2 dV_g\spn{1.2.3} \\
=&&\hskip-7pt \left ( 2 C_2(E) - \frac{r-1}{r}C_1(E)^2 \right )\cdot[\Omega],  
 \nonumber\end{eqnarray}
where $[\Omega]$ denotes the cohomology class of $\Omega${\rm .}
\endproclaim 

\demo{Proof} Recall that $D_A^\ast = -\ast D_A \ast $, so that 
\begin{eqnarray}
 D_A^\ast F_A 
&=&\frac{1}{r} D_A^\ast (\tr(F_A) {\rm Id})
+ \ast D_A (\Omega \wedge (F_A - \frac{1}{r}\tr (F_A){\rm  Id}))\nonumber\\[1ex]
&=& \frac{1}{r} d^\ast (\tr (F_A)) Id + \ast (\Omega \wedge ( D_A F_A
- \frac{1}{r} d (\tr (F_A)){\rm  Id})\nonumber\\  [1ex]
&=& 0.\nonumber
\end{eqnarray}
\vglue-4pt\noindent 
Hence, $A$ is a Yang-Mills connection. 

Next, multiplying (1.2.2) by $F_A$ and integrating 
the resulting identity over $M$, we get
\begin{eqnarray}
 && \hskip-7pt \left ( 2 C_2(E) - \frac{r-1}{r} C_1(E)^2\right )\cdot[\Omega]\nonumber\\ [2ex]
  =&&\hskip-7pt \left ( -Ch_2(E) + \frac{1}{r} C_1(E)^2 \right ) \cdot [\Omega]
\nonumber\\[2ex]
 =&&\hskip-7pt  \frac{1}{4\pi^2}\int _M \tr \left ((F_A - \frac{1}{r} \tr ( F_A){\rm Id})\wedge
(F_A - \frac{1}{r} \tr ( F_A){\rm Id})\right ) \wedge \Omega\nonumber\\[2ex]
 =&&\hskip-7pt   -\frac{1}{4\pi^2}\int _M \tr \left ((F_A - \frac{1}{r} \tr ( F_A){\rm  Id})\wedge
\ast (F_A - \frac{1}{r} \tr ( F_A){\rm Id} )\right ) \nonumber\\[2ex]
 =&&\hskip-7pt  \frac{1}{4\pi^2}\int _M \left ( |F_A |^2 - \frac{1}{r} |\tr ( F_A) |^2
\right ) dV_g,\nonumber
\end{eqnarray}
where $C_i(E)$ denotes the $i^{\rm th}$ Chern class of $E$
and $Ch_i(E)$ denotes the $i^{\rm th}$ Chern character of $E$.
Then (1.2.3) follows.
\enddemo  

In general, (1.2.2) is an over-determined system and has no solutions.
However, if $A$ is a solution of (1.2.2) and the co-norm of $\Omega $ is less than one, 
then $A$ is an absolute minimizer of $YM$ (cf.\ [HL]).

We will call any solution $A$ of (1.2.2) an $\Omega$-anti-self-dual
instanton. If there is no possible confusion,
we will simply say that $A$ is an anti-self-dual instanton.

\demo{{R}emark {\rm 1}}
For a general compact Lie group, we can also define the\break $\Omega$-anti-self-duality
instantons simply as the solutions of
$-\ast (F_A \wedge \Omega) = F_A$.
\enddemo  

In the following and next two sections,
we will give some solutions of~(1.2.2).

Now we let $M$ be a complex $m$-dimensional K\"ahler manifold with
a K\"{a}hler metric $g$. As usual, we denote by
$\omega=\omega_g$ the associated K\"{a}hler form. Then
\begin{equation}
\label{eq:1.2.1} 
  dV_g=\frac{\omega^m}{m!} \, . \spn{1.2.4}
\end{equation}

For any $U(r)$-connection $A$ of a complex bundle $E$ over $M$, 
we can decompose
\begin{equation}
\label{eq:1.2.2}
F_A=F^{2,0}_{A}+F^{1,1}_{A}+F^{0,2}_{A} \spn{1.2.5}
\end{equation}
where $F^{0,2}_{A}$ denotes the $(0,2)$-part of $F_A$,
$F^{2,0}_{A}= - (F^{0,2}_{A})^\ast$ and
$F^{1,1}_{A}$ denotes the $(1,1)$-part of $F_A$.
  
By the Newlander-Nirenberg
theorem, the vanishing of $F^{0,2}_{A}$ is equivalent to the integrability of
$\bar{\partial}_A=D^{0,1}_A$, which is the $(0,1)$-part of $D_A$; that is,
$\pi:E\to M$ has a holomorphic
structure induced by $D^{0,1}_A$.

Since $A$ is unitary,  
\begin{equation}
\label{eq:1.2.3}
F_A^{1,1}= - (F_A^{1,1})^{\ast}~~{\rm and}~~|F_A|^2=|F^{1,1}_{A}|^2+2|F^{0,2}_{A}|^2. \spn{1.2.6}
\end{equation}

We introduce notation:
\begin{equation}
\label{eq:1.2.4}
H_A=(F^{1,1}_{A}\cdot\omega),\quad
\overset \circ F \!\!~^{1,1}_{A}=F^{1,1}_{A}
-\frac{1}{m}H_A\omega \spn{1.2.7}
\end{equation}
where $F^{1,1}_{A}\cdot\omega$ denotes the orthogonal projection
of $F^{1,1}_{A}$ in the $\omega $-direction.

Now we set 
$$
\Omega ~=~ \frac{\omega ^{m-2}}{(m-2)!} 
$$
and we have:
\specialnumber{1.2.2} \proclaim{Proposition}
\label{prop:1.2.2}
The unitary connection $A$ satisfies {\rm (1.2.2)} if and only if
$\tr (F_A)$ is harmonic and
$$
F^{0,2}_A = \frac{1}{r} \tr(F^{0,2}) {\rm Id},~~H_A - \frac{1}{r} \tr (F_A^{1,1}\cdot \omega) {\rm Id} =0.
$$ 
If $C_1(E)$ is of the type
$(1,1)${\rm ,} then $A$ satisfies {\rm (1.2.2)} if and only if 
$$
F^{0,2}_A = 0,~~~H_A = \lambda {\rm Id}, 
$$
where $\lambda = \frac{m (C_1(E)\cdot [\omega]^{m-1})}{r[\omega]^m}${\rm .}

Furthermore{\rm ,} $A$ is the absolute minimum of the
Yang\/{\rm -}\/Mills functional if 
$$
F^{0,2}_A = 0,~~~H_A = \lambda {\rm Id}. 
$$
In this case{\rm ,  }
\begin{equation}
\label{eq:1.2.5}
YM(A)=(2C_2(E)- C_1(E)^2)\cdot 
\frac{[\omega]^{m-2}}{(m-2)!} + 
\frac{m(C_1(E)\cdot [\omega]^{m-1})^2}{r (m-1)![\omega]^m}, 
\qquad\spn{1.2.8}
\end{equation}
where $[\omega]$ denotes the cohomology class represented by
$\omega${\rm .}
\endproclaim

The proof follows from (1.2.3) and direct computations,
\begin{eqnarray}
\label{eq:1.2.6}
&& 
\hskip-7pt 4\pi^2 (2C_2(E)-C_1(E)^2)\cdot [\Omega ] \spn{1.2.9}  \\ [2ex]
 = &&\hskip-7pt \int_M 
            \left(\mid\fcirc\!\!~_{A}^{1,1}\mid^2-2|
            F^{0,2}_{A}|^2-\frac{m-1}{m}|H_A|^2 \right)\frac{\omega^m}{m!}\nonumber \\[2ex]
 = &&\hskip-7pt \int_M \left(\mid F_{A}\mid^2-4|
            F^{0,2}_{A}|^2-|H_A|^2 \right)\frac{\omega^m}{m!} .\nonumber
\end{eqnarray} 

\demo{Definition {\rm 1.2.3}}
 We call $A$ a {\it Hermitian\/{\rm -}\/Yang\/{\rm -}\/Mills connection} of $E$
if $A$ is unitary and 
$$
F^{1,1}_{A}\cdot \omega= \lambda {\rm Id},~~~F^{0,2}_{A}=0,
$$ 
where $\lambda = \frac{m(C_1(E)\cdot [\omega]^{m-1})}{r[\omega]^m}$. 
\enddemo

It follows from Proposition 1.2.2 that the action $YM(A)$ of 
any Hermitian-Yang-Mills connection $A$ is uniquely
determined by $E$ and the K\"ahler class $[\omega]$.

As we said, each Hermitian-Yang-Mills connection gives
rise to a natural holomorphic structure on $E$. In fact, by the
Donaldson-Uhlenbeck-Yau theorem, irreducible
Hermitian-Yang-Mills connections are in one-to-one
correspondence with stable holomorphic bundles over $M$.

\demo{{\rm 1.3.}\ Complex anti-self-dual instantons}
 In this section, we will discuss complex anti-self-dual
instantons on $4$-dimensional Calabi-Yau manifolds, 
as well as instantons on manifolds with
special holonomy. Complex anti-self-dual
instantons were previously studied by both mathematicians and physicists,
notably Donaldson and Thomas. We recommend the readers to the excellent reference [DT]
for  a more complete history.

First we assume that $M$ is a Calabi-Yau $4$-fold with a K\"{a}hler
metric $\omega$ and a holomorphic $(4,0)$-form $\theta$.
Furthermore, we normalize
\begin{equation}
\label{eq:1.3.1}
\theta\wedge\bar{\theta}=\frac{\omega^4}{4!} . \spn{1.3.1}
\end{equation}
Note that such a $\theta$ is only unique modulo multiplication by
units in $\CC$.

We now choose $\Omega$ to be the parallel form
$$
 4 {\rm Re}(\theta) + \frac{1}{2} \omega^2.
$$
Then solutions of (1.2.2) can be described as follows.

Let $h$ be a fixed hermitian metric of $\pi:E\to M$.
Then one can define a complex Hodge operator
\begin{equation}
\label{eq:1.3.2}
\ast_\theta:\Omega^{0,2}(\End(E))\to \Omega^{0,2}(\End(E)) \spn{1.3.2}
\end{equation}
by the equation
\begin{equation}
\label{eq:1.3.3}
-\tr(\varphi\wedge
\ast_\theta\psi)= ( \varphi , \psi ) \,\bar{\theta},\quad\forall
\varphi , \psi\in\Omega^{0,2}(\End(E)), \spn{1.3.3}
\end{equation}
where $( \cdot , \cdot )$ denotes the inner product on
$\Omega^{0,2}(\End(E))$ induced by the K\"{a}hler metric $\omega$
and the hermitian metric $h$ on $E$. More explicitly,
let $\{\varphi_1, \varphi_2,
\varphi_3,\varphi_4\}$ be any unitary coframe of $\omega$ such that
\begin{eqnarray}
\omega &=&\frac{\sqrt{-1}}{2}
\sum _i \varphi_i\wedge\bar{\varphi}_i,\nonumber\\[2ex]
\theta &=&-\frac{1}{4}\varphi_1
\wedge \varphi_2\wedge\varphi_3\wedge\varphi_4.\nonumber 
\end{eqnarray}
Then
\begin{eqnarray}
\ast_\theta(\sigma\bar{\varphi}_1\wedge\bar{\varphi}_2)&=&\sigma^\ast
\bar{\varphi}_3\wedge\bar{\varphi}_4,\nonumber\\[2ex]
\ast_\theta(\sigma\bar{\varphi}_1\wedge\bar{\varphi}_3)&=&\sigma^\ast
\bar{\varphi}_4\wedge\bar{\varphi}_2,\nonumber\\[2ex]
\ast_\theta(\sigma\bar{\varphi}_1\wedge\bar{\varphi}_4)&=&\sigma^\ast
\bar{\varphi}_2\wedge\bar{\varphi}_3,\nonumber
\end{eqnarray}
where $\sigma\in \End(E)$ and $\sigma^\ast$ denotes its adjoint
with respect to the Hermitian metric $h$ on $E$.

Let $A$ be an $\Omega$-anti-self-dual connection, i.e.,
$\tr(F_A)$ is harmonic and
$$
\Omega \wedge (F_A - \frac{1}{r}\tr(F_A){\rm  Id})= - \ast (F_A - \frac{1}{r}\tr(F_A){\rm  Id}).
$$
As in last section, we decompose 
$$
F_A = F_A^{2,0} + F_A^{1,1} +F^{0,2}_A.
$$
Then by direct computations, one can show that
the above is equivalent to the system
\begin{eqnarray}
\label{eq:1.3.4}
F^{1,1}_{A}\cdot\omega&=&\lambda {\rm Id},\spn{1.3.4}\\[2ex]
(d + d^* )\tr(F_A^{0,2}) &=& 0,\nonumber\\[2ex]
(1+\ast_\theta)(F^{0,2}_{A}- \frac{1}{r}
\tr(F_A^{0,2}) {\rm Id} )&=&0,\nonumber
\end{eqnarray}
where 
\begin{equation}
\label{eq:1.3.5}
\lambda=\frac{4C_1(E)\cdot[\omega]^3}{r [\omega]^4} \spn{1.3.5}
\end{equation}
Note that $\ast_\theta$ induces a decomposition
of $H^{0,2}(M,\CC )$ into the self-dual part and anti-self-dual part.
For any solution $A$ of (1.3.4),
$(1+\ast_\theta )\tr(F_A^{0,2})$ is harmonic and represents
the self-dual part of $C_1(E)^{0,2}$. In particular, if
$C_1(E)^{0,2}$ is anti-self-dual, then (1.3.4) reduces to
\begin{equation}
\label{eq:1.3.6}
(1+\ast_\theta)F^{0,2}_{A}=0, ~~~F^{1,1}_{A}\cdot\omega=\lambda {\rm Id}. \spn{1.3.6}
\end{equation}
Following [DT], we say that $A$ is a complex anti-self-dual instanton 
associated to $(E,h)$, if $D_Ah=0$ and $F_A$ satisfies (1.3.4).

For such a connection $A$, we observe
\begin{equation}
\label{eq:1.3.7}
[\theta]\wedge(2C_2(E)- \frac{r-1}{r} C_1(E)^2)=\frac{1}{4\pi^2}
\int_M \left |F^{0,2}_{A}- \frac{1}{r} \tr(F^{0,2}_A) {\rm Id} \right |^2 
\frac{\omega^4}{4!}, \quad\spn{1.3.7}
\end{equation}
where $[\theta]$ denotes the cohomology class of $\theta$ in 
$ H^4(M, \CC )$.  
Hence, (1.3.4) has no solutions if $[\theta]\wedge(2C_2(E)- 
\frac{r-1}{r}C_1(E)^2)$ is not a nonnegative real number. 
Since $\theta$ is only unique modulo multiplication by
units in $\CC$, for any given complex bundle $\pi:E\to M$, 
we should normalize $\theta$ such that
\begin{equation}
\label{eq:1.3.8}
[\theta]\wedge(2C_2(E)-\frac{r-1}{r}C_1(E)^2)\ge 0. \spn{1.3.8}
\end{equation}
Clearly, if this is not zero, then such a $\theta$ is unique once
$\omega$ is fixed. Moreover, if $C_1(E)^2\cdot [\theta]=0$
and $C_2(E)\cdot[\theta] =0$, then any complex
anti-self-dual instanton of $E$ is automatically a
Hermitian-Yang-Mills connection, which can be thought of  as holomorphically
flat. The readers may compare it to the Chern number conditions on the
flatness of Hermitian-Yangs-Mills connections.

The following proposition can be proved by straightforward
computations. \enddemo

\specialnumber{1.3.1} \proclaim{Proposition}
\label{prop:1.3.1}
Assume that $\theta$ is chosen so that {\rm (1.3.8)} holds{\rm . }
Let $A$ be any complex anti\/{\rm -}\/self\/{\rm -}\/dual instanton{\rm ,} then
\begin{eqnarray}
\label{eq:1.3.9}
&& \spn{1.3.9}\\
YM(A)&=& (2C_2(E)-C_1(E)^2)\cdot
\frac{[\omega]^2}{2}+\frac{4(C_1(E)\cdot
  [\omega]^3)^2}{6 r[\omega]^4}\nonumber \\[2ex]
&& ~+\  4 (2C_2(E)-\frac{r-1}{r}C_1(E)^2)\cdot [\theta]+
\frac{1}{r \pi^2}\int_M|\tr(F^{0,2}_A)|^2 dV_g.\nonumber
\end{eqnarray}
\endproclaim

It follows that each complex anti-self-dual instanton attains
the absolute minimum of the Yang-Mills functional. Moreover,
its action depends only on $E$, $[\omega]$ and $[\theta]$.

Calabi-Yau $4$-folds have holonomy group ${\rm SU}(4)$, which is contained in
${\rm Spin}(7)$. It turns out that complex anti-self-dual instantons can also be
 defined on ${\rm Spin}(7)$-manifolds, which have ${\rm Spin}(7)$ as 
their holonomy group (cf.\ [DT]).

Now let $(M,g)$ be a ${\rm Spin}(7)$-manifold. Then ${\rm Spin}(7)$,
acting on $\wedge^4(M)$, the space of $4$-forms, leaves invariant
a parallel $4$-form $\Omega \not= 0$. More explicitly, in terms of
an orthonormal basis $\{e_i\}$, the form
\begin{eqnarray*}
\Omega&=&e_1\wedge e_2\wedge e_5\wedge e_6 +e_1\wedge e_2\wedge e_7\wedge e_8
+e_3\wedge e_4\wedge e_5\wedge e_6\\[2ex]
&&+\ e_3\wedge e_4\wedge e_7\wedge e_8 +e_1\wedge e_3\wedge e_5\wedge e_7
-e_1\wedge e_3\wedge e_6\wedge e_8\\[2ex]
&&-\  e_1\wedge e_4\wedge e_5\wedge e_7+e_2\wedge e_4\wedge e_6\wedge e_8
-e_1\wedge e_4\wedge e_5\wedge e_8\\[2ex]
&&-\  e_1\wedge e_4\wedge e_6\wedge e_7-e_2\wedge e_3\wedge e_5\wedge e_8
-e_2\wedge e_3\wedge e_6\wedge e_7\\[2ex]
&&+\ e_1\wedge e_2\wedge e_3\wedge e_4+e_5\wedge e_6\wedge e_7\wedge e_8.
\end{eqnarray*}
If $M$ happens to be a Calabi-Yau $4$-fold, then it is the same as the one given above.  

One observes that the operator $\phi \mapsto - \ast (\Omega \wedge \phi)$
is self-adjoint on\break $2$-forms and has eigenvalues $1$ and $-3$. Following
[BKS], we let $\Omega_{21}^2(M, {\rm End}(E))$ and $\Omega _+^2(M, {\rm End}(E))$
be its eigenspaces corresponding to eigenvalues $1$ and $-3$.
Given any connection $A$, we write $F_A = F_{A,-}+ F_{A,+}$
according to this decomposition. Then $A$ solves (1.2.2) if and
only if $F_{A,+} = \frac{1}{r} \tr(F_{A,+}) Id$ and $\tr (F_A)$ is harmonic. 
Moreover, we have the identity
\begin{eqnarray*}
&&\hskip-7pt (2C_2(E)- \frac{r-1}{r} C_1(E)^2)\cdot [\Omega] \\[2ex]
  = &&\hskip-7pt \frac{1}{4\pi^2}
\int_M \left ( |F_{A,-}-\frac{1}{r} \tr(F_{A,-}) Id |^2 - 
3 |F_{A,+}-\frac{1}{r} \tr(F_{A,+}) Id |^2 \right ) dV_g.
\end{eqnarray*}

Therefore:

\specialnumber{1.3.2} \proclaim{Proposition}
\label{prop:1.3.2}
Let $(M,g)$ be a ${\rm Spin}(7)$\/{\rm -}\/manifold{\rm ,} and $A$ be an\break 
$\Omega$\/{\rm -}\/anti\/{\rm -}\/self\/{\rm -}\/dual
instanton{\rm .} Then $F_{A,+} =\frac{1}{r} \tr(F_{A,+}) {\rm Id} $
and $YM(A)$ depends only on $M$ and $E${\rm .} In fact{\rm ,} 
\begin{equation}
\label{eq:1.3.10}
YM(A)=  (2C_2(E)-\frac{r-1}{r} C_1(E)^2)\cdot [\Omega] +\frac{1}{4 r \pi^2}
\int_M | \tr(F_{A})|^2 dV_g . \qquad \spn{1.3.10}
\end{equation}
\endproclaim 

\demo{{\rm 1.4.}\ Instantons on $G_2$-manifolds}
Let $(M,g)$ be a Riemannian manifold with holonomy group being
the exceptional group $G_2$. Then there is a parallel, hence closed,
3-form $\Omega$ which is invariant under the action of $G_2$.
In terms of an orthonormal basis $\{e_i\}$, this form
\begin{eqnarray*}
\Omega&=&e_1\wedge e_2\wedge e_3 +e_1\wedge e_4\wedge e_5
-e_1\wedge e_6\wedge e_7\\[2ex]
&&+\ e_2\wedge e_4\wedge e_6 + e_2\wedge e_5\wedge e_7+ e_3\wedge e_4\wedge e_7
- e_3\wedge e_5\wedge e_6.
\end{eqnarray*}
The operator $\phi \mapsto - \ast (\Omega \wedge \phi)$
is self-adjoint on $2$-forms and has eigenvalues $1$ and $-2$.
We denote by $\Omega_{12}^2(M, {\rm End}(E))$ and $\Omega _+^2(M, {\rm End}(E) )$
its eigenspaces corresponding to eigenvalues $1$ and $-2$.
Given any connection $A$, we write $F_A = F_{A,-}+ F_{A,+}$
according to this eigenspace decomposition. 
Then $A$ is an $\Omega$-anti-self-dual instanton 
if and only if $F_{A,+} = \frac{1}{r} \tr(F_{A,+}) Id$ and
$\tr(F_{A})$ is harmonic. 
Moreover, we have the identity
\begin{eqnarray*}
&&\hskip-7pt (2C_2(E)- \frac{r-1}{r}C_1(E)^2)\cdot [\Omega] \\[2ex]
 = &&\hskip-7pt \frac{1}{4\pi^2}\int_M \left ( |F_{A,-}-\frac{1}{r} \tr(F_{A,-}) 
Id|^2 - 2 |F_{A,+}-\frac{1}{r} \tr(F_{A,+}) {\rm Id} |^2 \right ) dV_g.
\end{eqnarray*}
Therefore: \enddemo

\specialnumber{1.4.1} \proclaim{Proposition}
\label{prop:1.4.1}
Let $(M,g)$ be a $G_2$-manifold{\rm ,} and $A$ be an $\Omega$\/{\rm -}\/anti\/{\rm -}\/self\/{\rm -}\/dual
instanton{\rm ,} where $\Omega$ is the above $3$\/{\rm -}\/form defining the $G_2$\/{\rm -}\/structure{\rm . }
Then $F_{A,+} =\frac{1}{r} \tr(F_{A,+}) {\rm Id}$ 
and $YM(A)$ depends only on $M$ and $E${\rm .} In fact{\rm ,} 
\begin{equation}
\label{eq:1.4.1}
YM(A)= (2C_2(E)- \frac{r-1}{r}C_1(E)^2)\cdot [\Omega] + \frac{1}{4 r\pi^2}
\int _M |\tr(F_{A})|^2 dV_g. \qquad \spn{1.4.1}
\end{equation}
\endproclaim 

\section{Consequences of a monotonicity formula}
\label{chap:2}

In this chapter, we discuss Price's monotonicity formula, Uhlenbeck's
curvature estimate and singular Yang-Mills connections of a certain type. 

\demo{{\rm 2.1.}\ A monotonicity formula}
In this section, we will derive a monotonicity formula for Yang-Mills
connections, which is essentially due to Price [Pr]. This formula 
will be used in establishing cone properties of blow-up loci.
Its proof follows Price's arguments with some
modifications. 

As before, $M$ denotes a Riemannian manifold with a metric $g$
and $E$ is a vector bundle over $M$ with compact structure group $G$.

For any connection $A$ of $E$, its
curvature form $F_A$ takes values in ${\rm Lie}(G)$. The
norm of $F_A$ at any $p\in M$ is given by
\begin{equation}
\label{eq:2.1.1}
|F_A|^2=\sum^n_{i,j=1} \left< F_A(e_i,e_j),F_A(e_i,e_j) \right>, \spn{2.1.1}
\end{equation}
where $\{e_i\}$ is any orthonormal basis of $T_pM$,
and $\langle \cdot,\cdot \rangle$ is the Killing form of ${\rm Lie}(G)$.  

Let $\{\phi_t\}_{|t|<\infty}$ be a one-parameter family of
diffeomorphisms of $M$, and $A_0$ be a fixed smooth connection of $E$
and $D$ be its associated covariant derivative. Then for any connection $A$,
we can define a family of connections $\phi^\ast_{t}(A)$ as follows:
Denote by $\tau^0_{t}$ the parallel transport of $E$ associated to $A_0$ along the
path ${\phi_s(x)}_{0\leq s\leq t}$, where $x\in M$.  More
precisely, for any $u\in E_x$ over $x\in M$, let $\tau^0_s(u)$ be
the section of $E$ over the path ${\phi_s(x)}_{0\leq s\leq t}$
such that
\begin{equation}
\label{eq:2.1.2}
D_{\frac{\partial}{\partial s}}\tau^0_{s}(u)=0,\quad \tau^0_{0}(u)=u. \spn{2.1.2}
\end{equation}
We define $A^t=\phi^\ast_{t}(A)$ by defining its associated
covariant derivative 
\begin{equation}
\label{eq:2.1.3}
D^t_{X}v=(\tau^0_t)^{-1}\left ( D_{d\phi_t(X)}(\tau^0_t(v))\right )  \spn{2.1.3}
\end{equation}
for any $X\in TM, v\in\Gamma(M,E)$, where $\Gamma(M,E)$ is the
space of sections of $E$ over $M$.

To see that $A^t$ is indeed a connection, it is sufficient to check
\begin{eqnarray}
&&\hskip-36pt D^t_X(fv)(x)  \nonumber \\[1ex]
&=&(\tau^0_{t})^{-1}\left (D_{d\phi_t(X)}((\phi^{-1}_{t})^{\ast}f 
         \cdot \phi_t )\tau^0_{t}(v))\right )(x)\nonumber\\ [1ex]
&=&(\tau^0_{t})^{-1}\left(f(x)D_{d\phi_t(X)}\tau^0_{t}(v)(\phi_t(x))
          +d\phi_t(x)\left((\phi^{-1}_{t})^{\ast}f\right)\tau^0_{t}(v)
(\phi_t(x))\right) 
            \nonumber\\[1ex]
&=&f(x)D^{t}_{X}v (x) + X(f)(x)v(x)\nonumber.
\end{eqnarray} 
The curvature form of $A^t$ is then given by 
\begin{equation}
\label{eq:2.1.4}
F_{A^{t}}(X,Y)=(\tau^0_{t})^{-1}\cdot
F_A(d\phi_t(X), d\phi_{t}(Y))\cdot \tau^0_{t}.  \spn{2.1.4}
\end{equation}
It follows that
\begin{eqnarray}
\label{eq:2.1.5}
 \qquad YM(A^t) 
&=&\frac{1}{4\pi^2}\int_M|F_{A^{t}}|^{2}dV_g\spn{2.1.5}\\[2ex]
&=&\frac{1}{4\pi^2}\int_M\sum^{n}_{i,j=1}|F_A ( 
d\phi_t(e_i), d\phi_t(e_j))|^{2}
         (\phi_t(x))dV_g(x) \, , \nonumber
\end{eqnarray}
where $dV_g$ denotes the volume form of $g$, and $\{e_i\}$ is any
local orthonormal basis of $TM$.

By changing variables, we obtain
$$
YM(A^t)=\frac{1}{4\pi^2}\int_M\sum^{n}_{i,j=1}|F_A (d\phi_t(e_i(\phi_t^{-1}(x))),
d\phi_t(e_j(\phi^{-1}_{t}(x))))|^{2}{\rm Jac}(\phi^{-1}_{t})dV_g.
$$
Let $X$ be the vector field 
$\frac{\partial\phi_{t}}{\partial t}|_{t=0}$ on $M$. Then we deduce 
from the above that
\begin{eqnarray}
\label{eq:2.1.6} &&\spn{2.1.6}\\
&&   \frac{d}{dt}YM(A^t)|_{t=0}\nonumber\\ 
&&\qquad  \enspace =\ -\frac{1}{4\pi^2}\int_M\left(|F_A|^{2}{\rm div} X+
      4\sum^{n}_{i,j=1}\langle 
F_A([X,e_i],e_j), F_A(e_i,e_j)\rangle\right)dV_g.\nonumber
\end{eqnarray}
Here we have used the formula  
$$\frac{d}{dt}\left(d\phi_t(e_i(\phi^{-1}_{t}(x)))\right)|_{x=0}=-[X,e_i].$$
Since $[X,e_i]=\nabla_Xe_i-\nabla_{e_i}X,$ where $\nabla$ denotes
the Levi-Civita connection of $g$, we obtain
\begin{eqnarray}
\label{eq:2.1.7}
&&\qquad\sum^{n}_{i,j=1}
      \left<F_A([X,e_i],e_j), F_A(e_i,e_j)\right> \spn{2.1.7}\\[2ex]
&&\qquad\quad =\ -\sum^{n}_{i,j=1}
    \Bigl(\langle F_A(\nabla_{e_i}X,e_j),
                  F_A(e_i,e_j) \rangle
         -\langle F_A(\nabla_{X}e_i,e_j),
                  F_A(e_i,e_j) \rangle \Bigr)\nonumber\\ [2ex]
&&\qquad\quad  =\ -\sum^{n}_{i,j=1} 
    \Bigl( \langle F_A(\nabla_{e_i}X,e_j),F_A(e_i,e_j) \rangle \nonumber\\
&&\qquad\quad\hskip1in 
- g(\nabla_Xe_i,e_k)
    \langle F_A(e_k,e_j),F_A(e_i,e_j) \rangle \Bigr).\nonumber 
\end{eqnarray}  
Since 
$$
g(\nabla_Xe_i,e_k)=-g(e_i,\nabla_Xe_k)=-g(\nabla_Xe_k,e_i),
$$ 
the second term in (2.1.7) vanishes.

Now suppose that $A$ is a Yang-Mills connection; then
\begin{equation}
\label{eq:2.1.8}
0=\int_M\left(|F_A|^{2}\divv X-4\sum^{n}_{i,j=1}\langle F_A(\nabla_{e_i}X,e_j),
  F_A(e_i,e_j)\rangle\right)dV_g. \hskip.25in \spn{2.1.8}
\end{equation}
The required monotonicity will follow from this variational formula.

Fix any $p\in M$, let $r_p$ be a positive number with properties:
there are normal coordinates $x_1, \cdots, x_n$ in the geodesic ball
$B_{r_p}(p)$ of $(M, g)$, such that $p=(0,\cdots,0)$ and for some constant $c(p)$,
\begin{eqnarray}
\label{eq:2.1.9}
|g_{ij}-\delta_{ij}| &\le& c(p)(|x_1|^{2}+\cdots +|x_n|^{2}), \spn{2.1.9}\\[2ex]
|dg_{ij}| &\le & c(p) \sqrt{|x_1|^{2}+\cdots +|x_n|^{2}},  \spn{2.1.10}
\end{eqnarray}
where
\begin{equation}
\label{eq.2.1.11}
g_{ij}=g \left( \frac{\partial}{\partial x_{i}},
                \frac{\partial}{\partial x_j} \right). \spn{2.1.11}
\end{equation}

\demo{{R}emark {\rm 2}}
\label{rem:2.1.1}
The constants $r_p$ and $c(p)$ can be chosen depending only on the injective
radius at $p$ and the curvature of $g$. If $M=\RR^n$ and $g$ is flat, we can take
$r_p = \infty$ and $c(p) =0$.
\enddemo

Let $r(x)$ be the distance function from $p$; i.e., 
$$
r(x) = \sqrt{x^2_1 + \cdots + x_n^2}.
$$
Let $\phi$ be a positive function on the unit sphere $S^{n-1}$.
Define
\begin{equation}
\label{eq:2.1.12}
X(x)=\xi(r)\phi(\frac{x}{r})r\frac{\partial}{\partial
  r}=\xi(r) \phi  ({x\over r} ) \left(\sum_{i}x_i\frac{\partial}{\partial x_{i}}\right), \spn{2.1.12}
\end{equation}
where $\xi$ is some smooth function with compact support in $B_{r_p}(p)$. 

Let $\{e_1,\,\cdots \,,e_n\}$ be any orthonormal basis near $p$ such that
$e_1=\frac{\partial}{\partial r}$. Since $x_1,\cdots ,x_n$ are normal
coordinates, we have
$$
\nabla _{\frac{\partial}{\partial r}}\frac{\partial}{\partial r} =0.
$$
It follows that 
\begin{equation}
\label{eq:2.1.13}
\nabla_{\frac{\partial}{\partial r}}X=(\xi
r)^{\prime}\phi (\theta) \frac{\partial}{\partial r}
=(\xi^{\prime}r+\xi)\phi (\theta)\frac{\partial}{\partial r}, \spn{2.1.13}
\end{equation}
where $\theta = \frac{x}{r}$. 
Moreover, for $i\geq 2$,  
\begin{equation}
\label{eq:2.1.14}
\nabla_{e_{i}}X=\xi r\nabla_{e_{i}}(\phi
\frac{\partial}{\partial r})
=\xi r e_{i}(\phi) \frac{\partial}{\partial r}+\ \xi \phi \sum _{j=1}^n b_{ij} e_j, \hskip.75in\spn{2.1.14}
\end{equation}
where $|b_{ij} -\delta _{ij}| = O(1) c(p) r^2$. 
We will always denote by $O(1)$
a quantity bounded by a constant depending only on $n$.

Applying (2.1.13) and (2.1.14) to the first variational formula
(2.1.8), we obtain
\begin{eqnarray}
\label{eq:2.1.15}
&&\hskip-4pt \int_M |F_A|^{2}(\xi^{\prime}r+(n-4)\xi+O(1)c(p) r^2 \xi)\phi dV_g\spn{2.1.15} \\[2ex]
&&\qquad\qquad =\ 4\int _M \left (\xi^{\prime}r\phi
|\frac{\partial}{\partial r}\rfloor F_{A}|^{2}+\xi r
\langle \frac{\partial}{\partial r}\rfloor F_{A}, \nabla \phi \rfloor
F_A \rangle\right)dV_g, \nonumber
\end{eqnarray}
where $\frac{\partial}{\partial r}\rfloor F_{A}=F_A(\frac{\partial}{\partial r},\cdot)$. 

We choose, for any $\tau$ small enough,
$\xi(r)=\xi_{\tau}(r)=\eta(\frac{r}{\tau})$, where $\eta$ is smooth
and satisfies: $\eta(r)=1$ for $r\in[0,1],\eta(r)=0$ for
$r\in[1+\varepsilon,\infty), \varepsilon>0$ and $\eta^{\prime}(r)\leq 0$. 
Then
\begin{equation}
\label{eq:2.1.16}
\tau\frac{\partial}{\partial\tau}\left(\xi_{\tau}(r)\right)=
    -r \xi^{\prime}_{\tau}(r). \spn{2.1.16}
\end{equation}
Plugging this into (\ref{eq:2.1.15}), we obtain
\begin{eqnarray}
\label{eq:2.1.17}
&&\spn{2.1.17}\\
&&\hskip-9pt\tau\frac{\partial}{\partial\tau}\left(\int_M\xi_{\tau}\phi | F_A|^2dV_g\right)
     +\left((4-n)+O(1)c(p) \tau ^2
\right)\int_M\xi_{\tau}\phi |F_A|^{2}dV_g \nonumber  \\[2ex]
&&\qquad = 4\tau\frac{\partial}{\partial\tau}
\left(\int_M\xi_{\tau}\phi |\frac{\partial}{\partial r}\rfloor F_{A}|^{2}dV_g\right)
- 4\int_{M}\xi_\tau r \langle \frac{\partial}{\partial
          r}\rfloor F_A , \nabla\phi \rfloor F_A\rangle dV_g. \nonumber
\end{eqnarray} 
Choose a nonnegative number $a \ge O(1) c(p)$. Then we deduce from the above
\begin{eqnarray}
\label{eq:2.1.18}  
&&\spn{2.1.18}\\
 &&\hskip-12pt \frac{\partial}{\partial\tau}\left(\tau^{4-n}e^{\pm a\tau^{2}}\int_M
   \xi_{\tau}\phi |F_A|^{2}dV_{g}\right)\nonumber \\[2ex]
 = && \hskip-12pt  4 \tau^{4-n}e^{\pm a\tau^{2}}
  \left
(\frac{\partial}{\partial\tau}\left(\int_M\xi_{\tau}\phi |\frac{\partial}
{\partial r}\rfloor F_A|^{2}dV_{g} \right) \right.%
  \nonumber\\ [2ex]
&& \hskip-12pt
 +\  (-O(1)c(p)\pm 2a)\tau \int_M\xi_{\tau}\phi |F_A|^{2}dV_g
 -\ \tau^{-1}
\left .\int_{M}\xi_\tau \langle \frac{\partial}{\partial r}
\rfloor F_A , \nabla\phi \rfloor F_A\rangle dV_g\right ). \nonumber  
\end{eqnarray} 
Then, by integrating on $\tau$ and letting $\varepsilon$ 
tend to zero, we prove:
\specialnumber{2.1.1} \proclaim{Theorem}
\label{th:2.1.1}
Let $r_p$, $c(p)$ and $a$ be as above{\rm . } Then for any
$0 < \sigma\break < \rho < r_p${\rm ,}  we have
\begin{eqnarray}
\label{eq:2.1.19}
&&\pm \rho^{4-n}e^{\pm a\rho^2}\int_{B_{\rho}(p)}\phi |F_{A}|^2 dV_g\mp \sigma^{4-n}e^{\pm
a\sigma^2}\int_{B_{\sigma}(p)}\phi |F_{A}|^2dV_g
\spn{2.1.19}\\
&&\mp\ 4\int_{B_\rho(p)\backslash B_\sigma(p)}r^{4-n}e^{\pm\ ar^2}\phi 
|\frac{\partial}{\partial r}\rfloor F_{A}|^2dV_g \nonumber\\
&&\ge    - 4\int^{\rho}_{\sigma}
 \tau^{3-n}e^{\pm\ a\tau^2} 
d\tau\int_{B_{\tau}(p)}|\frac{\partial}{\partial r}\rfloor F_{A}|
|\nabla\phi \rfloor F_{A}| dV_g. \nonumber
\end{eqnarray}
\endproclaim 

This inequality is needed for establishing the existence of tangent cones of blow-up loci.
Taking $\phi =1$, we obtain:
\proclaimtitle{Price}
\specialnumber{2.1.2} \proclaim{Theorem} 
\label{th:2.1.2}
Let $r_p$, $c(p)$ and $a$ be as above{\rm .} Then for any
$0 < \sigma < \rho < r_p${\rm ,} we have
\begin{eqnarray}
\label{eq:2.1.20}
&&\hskip-24pt \rho^{4-n}e^{a\rho^{2}}\int_{B_{\rho}(p)}|F_A|^{2}dV_g-
    \sigma^{4-n}e^{a\sigma^{2}}\int_{B_{\sigma}(p)}|F_A|^{2}dV_g \spn{2.1.20}\\[2ex]
&&\qquad\qquad\ge\ 4\int_{B_\rho(p)\backslash B_{\sigma}(p)}r^{4-n}e ^{a r^2}
|\frac{\partial}{\partial r}\rfloor
F_A|^{2}dV_g. \nonumber
\end{eqnarray} 
Moreover{\rm ,} if $M=\RR^n$ and $g$ is flat{\rm ,} then 
the equality holds in {\rm (2.1.20)} for $\rho\in(0,\infty)$ and $a=0${\rm .}
\endproclaim  

\demo{{R}emark {\rm 3}}
Both (2.1.19) and (2.1.20) still hold when $A$ is only a smooth Yang-Mills connection on
$M\backslash \{p\}$ with
$$
\int _M |F_A|^2 dV_g < \infty.
$$
To see this, we replace $\eta$ in defining $\xi$ in (2.1.16) by 
$\eta_\varepsilon$ for sufficiently small $\varepsilon$, where 
$\eta_\varepsilon (t) =0$ for either $t\le \varepsilon$ or $t \ge 1+\varepsilon$,
and $\eta_\varepsilon(t)=1$ for $t\in (\varepsilon, 1-\varepsilon)$. Then 
we can follow the same arguments from (2.1.16) on and obtain both
(2.1.19) and (2.1.20) for such a Yang-Mills connection $A$ with isolated singularity at
$p$.
\enddemo

It follows from this theorem that $\rho^{4-n}e^{a\rho^{2}}\int_{B_{\rho}(p)}|F_A|^{2}dV_g$
is a nondecreasing function in $(0, r_p)$.
Another simple corollary of (\ref{eq:2.1.20}) is the following:
\specialnumber{2.1.3} \proclaim{{C}orollary}
\label{cor:2.1.3} 
Let $A$ be a Yang\/{\rm -}\/Mills $G$\/{\rm -}\/connection
of the trivial bundle $(\RR^{n}\backslash\{0\})\times\RR^{r}${\rm ,}
such that $\rho^{4-n}\int_{B_{\rho}(0)}|F_A|^{2}dV_{g_{0}}$ is independent of\break
$\rho\in(0,\infty)${\rm ,} where $g_0$ is a flat metric on $\RR ^n${\rm .  }
Then $A$ is gauge equivalent to $d+ A_s${\rm ,} where $A_s:S^{n-1}\longrightarrow
T^{\ast}S^{n-1}\otimes {\rm Lie}(G)$ is a ${\rm Lie}(G)$-valued $1$\/{\rm -}\/form{\rm .}
\endproclaim
\demo{Proof} By (\ref{eq:2.1.20}) for the flat metric $g_0$, we obtain 
\begin{equation}
\label{eq:2.1.21}
\frac{\partial}{\partial r} \rfloor F_A\equiv 0. \spn{2.1.21}
\end{equation}
Let $D$ be the associated covariant derivative of $A$.
Write $D=d+\tilde{a}dr+\tilde{A}$, where
$\tilde{a}:\RR^{n}\backslash\{0\}\to {\rm Lie}(G)$ and
$\tilde{A}:\RR^{n}\backslash\{0\}\to T^{\ast}S^{n-1}\otimes {\rm Lie}(G)$.
For any gauge transformation $\sigma:\RR^{n}\backslash\{0\}\to G$,
we have a similar representation $d+\tilde{a}_{\sigma}dr+\tilde{A}_{\sigma}$ 
for $D_{\sigma(A)}$; moreover,  

\centerline{${\displaystyle
\tilde{a}_{\sigma}=(\sigma \cdot \tilde{a}
-\frac{\partial\sigma}{\partial r})  \cdot\sigma^{-1}. 
}$}
\smallbreak\noindent
Choose $\sigma$ by solving the ordinary differential equation on $G$: 
$$
\sigma \cdot \tilde{a}
-\frac{\partial\sigma}{\partial r}=0.
$$
Then $\tilde{a}_{\sigma}=0$. 
Together with (2.1.21), we deduce
$\frac{\partial\tilde{A}_{\sigma}}{\partial r}=0$;
i.e., $\tilde{A}_{\sigma}(r, \theta )=A_s(\theta )$ for some
$A_s:S^{n-1}\to T^\ast S^{n-1}\otimes {\rm Lie}(G)$. 
The corollary is proved.
\enddemo

\demo{{\rm 2.2.}\ Curvature estimates}
In this section, we give a basic curvature estimate for
Yang-Mills connections. This estimate was first derived by
K. Uhlenbeck (also see [Na]). Since it is crucial
to us here, we will outline its proof for the reader's convenience.

We will adopt the notation of the last section.
\proclaimtitle{K. Uhlenbeck}
\specialnumber{2.2.1} \proclaim{Theorem}
\label{th:2.2.1}
Let $A$ be any Yang\/{\rm -}\/Mills connection of a $G$\/{\rm -}\/bundle $E$ over $M${\rm .}
Then there are $\varepsilon = \varepsilon(n) > 0$ and $C=C(n) > 0${\rm , }
which depend only on $n$ and $M${\rm ,} such 
that for any $p\in M$ and $\rho < r_p${\rm ,} whenever
$$ 
\rho^{4-n}\int_{B_{\rho}(p)}|F_A|^2dV_g\break\leq
\varepsilon,
$$
\vglue-4pt\noindent 
then
\begin{equation}
\label{eq:2.2.1}
|F_A|(p)\leq\frac{C}{\rho^{2}}\left(\rho^{4-n}\int_{B_{\rho}(p)}
|F_A|^{2}dV_g\right)^{\frac{1}{2}}. \spn{2.2.1}
\end{equation}
\endproclaim
\phantom{woops}
\vglue-6pt

Our proof here uses R. Schoen's arguments in [Sc] for harmonic maps.  

By scaling, we may assume that $\rho=1$.
Define a function
\begin{equation}
\label{eq:2.2.2}
f(r)=(1-2r)^{2}\sup_{x\in B_{r}(p)}|F_A|(x),\qquad r\in [0,\frac{1}{2}]. \spn{2.2.2}
\end{equation}
Then $f(r)$ is continuous in $[0,\frac{1}{2}]$ with
$f(\frac{1}{2})=0$, so that  $f$ attains
its maximum at a certain $r_0$ in $[0, \frac{1}{2}]$.  

First we claim that $f(r_0)\leq 64$ if $\varepsilon$ is sufficiently small.
Assume that $f(r_0)>64$. Put $b =\sup_{x\in B_{r_{0}}(p)}|F_A|(x)=|F_A|(x_0)$; 
then taking $\sigma= \frac{1}{4} (1-2r_0)$, we get
\begin{eqnarray}
\label{eq:2.2.3}
 \sup_{x\in B_{\sigma}(x_{0})}|F_A|
&\leq&\sup_{x\in
  B_{r_{0}+\sigma}(p)}|F_A|(x)\spn{2.2.3}\\[2ex]
&\leq& \frac{(1-2r_0)^{2}}{(1-2r_0-2\sigma)^{2}}\sup_{x\in
  B_{r_0}(p)}|F_A|(x)=4b. \nonumber
\end{eqnarray} 
Clearly, $16\sigma^{2}b \geq 64$; i.e.,
$\sigma\sqrt{b}\geq 2$.  Define a scaled metric
$\tilde{g}=b g$; then with respect to $\tilde{g}$, the norm
$|F_A|_{\tilde{g}}$ of $F_A$ is equal to $b^{-1}|F_A|$. Hence,
\begin{equation}
\label{eq:2.2.4}
\sup_{x\in B_{2}(x_{0}, \tilde{g})}|F_A|_{\tilde{g}}\leq 4, \spn{2.2.4}
\end{equation}
where $B_2(x_{0},\tilde{g})$ denotes the geodesic ball of
$\tilde{g}$ with radius $2$ and center at $x_0$.
  
Since $A$ is a Yang-Mills connection,
by the second Bianchi identity and straightforward computations,
we can derive the following equation:
\begin{equation}
\label{eq:2.2.5}
\frac{1}{2}\Delta_{\tilde{g}}|F_A|^{2}_{\tilde{g}}=|\tilde{\nabla}
F_A|^{2}_{\tilde{g}}-2F_A{\#}F_A{\#}R(\tilde{g})-2F_A\ast F_A\ast
F_A, \hskip.5in \spn{2.2.5}
\end{equation}
where $F_A{\#}F_A{\#}R(\tilde{g})$ and $F_A\ast F_A\ast F_A$ are defined
as follows:  in any orthonormal basis
$e_1,\ldots ,e_n$ of $\tilde{g}$,
\begin{eqnarray}
\label{eq:2.2.6}
&& \spn{2.2.6}\\
 F_A{\#}F_A{\#}R(\tilde{g})
&=&\sum_{l,k,i,j}\left (\langle F_A(e_l,e_k),F_A(e_i,e_j)\rangle \phantom{\sum_{l}}
\right.\nonumber\\[2ex]
&&\left. -\sum _{m}\langle F_A(e_l,e_m),F_A(e_i,e_m)\rangle \delta_{jk}
\right )R(\tilde{g})(e_l,e_j,e_k,e_i),\nonumber
\end{eqnarray} 
and
\begin{equation}
\label{eq:2.2.7}
 F_A\ast F_A\ast F_A   \sum_{i,j,k}\langle \left[F_A(e_i,e_j),F_A(e_j,e_k)
               \right],F_A(e_k,e_i)\rangle .\spn{2.2.7}
\end{equation} 
It follows from (\ref{eq:2.2.5})--(\ref{eq:2.2.7}) that 
there are uniform constants $c_1$, $c_2$, depending only on $n$, such that
\begin{equation}
\label{eq:2.2.8}
-\Delta_{\tilde{g}}|F_A|_{\tilde{g}}\leq c_1|F_A|_{\tilde{g}}+
   c_2|F_A|^{2}_{\tilde{g}}. \spn{2.2.8}
\end{equation}
Using (2.2.4), we deduce from
(2.2.8) that in $B_{2}(x_{0},\tilde{g})$,
\begin{equation}
\label{eq:2.2.9}
-\Delta_{\tilde{g}}|F_A|_{\tilde{g}}\leq(c_1+4c_2)|F_A|_{\tilde{g}}. \spn{2.2.9}
\end{equation}
Then, by using either the mean-value theorem or the standard Moser
iteration, we obtain
\begin{equation}
\label{eq:2.2.10}
1=|F_A|_{\tilde{g}}(x_0)\leq \tilde{c}\left(\int_{B_{1}(x_{0},
\tilde{g})}|F_A|^{2}_{\tilde{g}}dV_{\tilde{g}}\right)^{\frac{1}{2}}, \spn{2.2.10}
\end{equation}
where $\tilde{c}$ is some uniform constant.

However, by the monotonicity (Theorem 2.1.1),  
\begin{eqnarray}
 \int_{B_{1}(x_0,\tilde{g})}|F_A|^{2}_{\tilde{g}}dV_{\tilde{g}} 
&=& (\sqrt{b})^{n-4}\int_{B_{\frac{1}{\sqrt{b}}}(x_{0})}|F_A|^{2}dV_g\nonumber\\[2ex]
&\leq& \left(\frac{1}{2} \right)^{4-n}e^{\frac{a}{4}}
      \int_{B_{\frac{1}{2}}(x_0)}|F_A|^{2}dV_g\ \leq\ \varepsilon 2^{n-4}e^{\frac{a}{4}}.\nonumber
\end{eqnarray}  
Combining this with (\ref{eq:2.2.10}), we obtain
$$
1\leq \tilde{c}\varepsilon 2^{n-4}e^{\frac{a}{4}}.
$$
It is impossible when $\varepsilon=\varepsilon(n)$ is
sufficiently small. The claim is proved.

Thus, we have
\begin{equation}
\label{eq:2.2.11}
\sup_{x\in B_{\frac{1}{4}}(p)}|F_A|(x)\leq 4f(r_{0})\leq 256. \spn{2.2.11}
\end{equation}
It follows from this and (\ref{eq:2.2.5}) with $\tilde{g}$
replaced by $g$ that for some uniform constant $c^\prime$,
\begin{equation}
\label{eq:2.2.12}
-\Delta _{g}|F_A|\leq c^{\prime}|F_A| . \spn{2.2.12}
\end{equation}
Then (\ref{eq:2.2.1}) follows
from (\ref{eq:2.2.12}) and a standard Moser iteration.
\enddemo

\demo{{\rm 2.3.}\ Admissible Yang-Mills connections}
\label{sec:2.3}
In order to compactify the moduli space of Yang-Mills connections, we
need to use singular Yang-Mills connections of a certain type.  
Those singular connections
behave like the  usual smooth connections in many ways; for
instance, one can define the first two terms of the Chern character by
using their curvature forms.

An admissible Yang-Mills connection is a smooth connection $A$
defined outside a closed subset $S(A)$ in $M$, such that (1)
$H^{n-4}(S(A)\cap K)<\infty$ for any compact subset $K\subset M$, 
where $H^{n-4}(\cdot)$ stands for the $(n-4)$-dimensional Hausdorff measure (cf.\ [Si2]); 
(2) $A$ is Yang-Mills on $M\backslash S(A)$; (3) $A$ satisfies
\begin{equation}
\label{eq:2.3.1}
\int_{M\backslash S(A)}|F_A|^2\,dV_g<\infty. \spn{2.3.1}
\end{equation}

Together with (2.3.1), this implies that for any smooth 
${\rm Lie}(G)$-valued\break $1$-form $u$ over $M$ with compact support,
\begin{equation}
\label{eq:2.3.2}
\int_M\left<F_A, du\right>dV_g=0. \spn{2.3.2}
\end{equation}

Clearly, $A$ is smooth on $M$ if $S(A)=\emptyset$.
We will call $S(A)$ the singular set of $A$. This is not invariant
under gauge transformations. Even if $S(A)\not= \emptyset$, there may be
a gauge transformation $\sigma$ on $M\backslash S(A)$ such that $\sigma(A)$
extends to become a smooth connection on $M$. 

Two admissible connections $A_1$ and $A_2$ are
gauge equivalent if there is a gauge transformation $\sigma$ of
$E$ over $M\backslash S(A_1)\cup S(A_2)$ such that $\sigma(A_1)=A_2$ 
outside $S(A_1)\cup S(A_2)$. This new gauge equivalence extends 
the previous one for smooth connections.

Similarly, by requiring that $A$ be $\Omega$-anti-self-dual outside
$S(A)$, we can also define admissible $\Omega$-anti-self-dual instantons.

Now let us assume that $G$ is a unitary group.
By the standard Chern-Weil theory, associated to each smooth
connection $A$, we have closed forms $\frac{\sqrt{-1}}{2\pi}\tr(F_A)$ and
$\left(\frac{\sqrt{-1}}{2\pi}\right)^2\tr(F_A\wedge F_A)$ of
degree $2$ or $4$. If $M$ is compact, they represent the first two Chern characters
${\rm Ch}_1(E)$ or ${\rm Ch}_2(E)$ respectively.  
We now extend these to  admissible Yang-Mills connections.

Let $A$ be an admissible Yang-Mills connection with the singular set $S=S(A)$. 
Then $\tr(F_A)$ and $\tr(F_A\wedge F_A)$ are closed forms on $M\backslash S$.
Because of    (3) above, we can extend them
to forms on $M$ in the sense of distribution. Clearly, these
forms are invariant under gauge transformations.
\specialnumber{2.3.1} \proclaim{Proposition}
\label{prop:2.3.1}
  The extended forms $\frac{\sqrt{-1}}{2\pi}\tr(F_A)$ and 
$(\frac{\sqrt{-1}}{2\pi})^2\tr(F_A\wedge F_A)$ are closed
on $M${\rm .} They are  denoted by ${\rm Ch}_1(A)$ and ${\rm Ch}_2(A)${\rm .}
\endproclaim
\demo{Proof}
We only show the closedness of $\left ( \frac{\sqrt{-1}}{2\pi} \right )^2
\tr(F_A\wedge F_A)$ here. The
other case is easier. We will always denote by $C$ a uniform 
constant in this proof.

It is sufficient to show that for any smooth form $\varphi$ of degree
$n-5$ with compact support in $M$,
\begin{equation}
\label{eq:2.3.3}
\int_M d\varphi\wedge\tr\left(F_A\wedge F_A\right)=0. \spn{2.3.3}
\end{equation}
Note that this is well-defined since $F_A$ is $L^2$-integrable.

Without loss of generality, we may assume that $M$ 
is a ball in $\RR^n$ and $E$ is a trivial
bundle over $M$. Let $K$ be a compact subset in $M$ containing 
$\supp(\varphi)$ in its interior.

Fixing any $\varepsilon \le \varepsilon(n)$, as given in Theorem 2.2.1,
we define
\begin{equation}
\label{eq:2.3.4}
E_r=\{x\in K\,|\,r^{4-n}e^{a r^2}
\int_{B_r(x)}|F_A|^2dV_g\geqslant\varepsilon\}, \spn{2.3.4}
\end{equation}
where $a$ is as in
Theorem~\ref{th:2.1.2}.
By Theorem 2.1.2, $E_r\supset E_{r^{\prime}}$ whenever $r\geqslant r'$.  
We can find a finite covering 
$\{B_{2r}(x_k)\}_{1\leqslant k \leqslant L_r}$ of
$E_r$ such that (1) $x_k\in E_r$;\break (2) $B_r(x_k)\cap
B_r(x_l)=\emptyset$ for $k\neq l$. Next
we expand $\{B_{2r}(x_k)\}_{1\le k\le L_r}$
to a covering $\{B_{2r}(x_k)\}_{1\leqslant k\leqslant L^{\prime}_{r}}$
($L^{\prime}_{r}\geqslant L_r$) 
of $(S\cap K)\cup E_r$, such that
$x_k \in (S\cap K)\cup E_r$, $B_r(x_k)\cap
B_r(x_l)=\emptyset$ for $k\neq l$. Note that
for any $k$, the number of $x_l$ with 
$B_{8r}(x_k)\cap B_{8r}(x_l)\not=\emptyset$ is bounded
by a constant depending only on $n$ and $M$.

For any $x\in\hspace{-2.0ex}\backslash\bigcup^{L^{\prime}_r}_
{k=1}B_{2r}(x_{k})$, 
\begin{equation}
\label{eq:2.3.5}
r^{4-n}\int_{B_{r}(x)}\left|F_A\right|^2\, dV_g<
\varepsilon. \spn{2.3.5}
\end{equation}
It follows from Uhlenbeck's estimate (Theorem~\ref{th:2.2.1})
that
\begin{equation}
\label{eq:2.3.6}
\left|F_A\right|(x)\leqslant\frac{C}{r^{2}}\left ( r^{4-n} \int _{B_r(x)} |F_A|^2 dV_g
\right )^{\frac{1}{2}}
\leqslant\frac{C\sqrt{\varepsilon}}{r^{2}}. \spn{2.3.6}
\end{equation}
Then, by using Theorem~ 1.2.7 in [Uh1, p.~18], we can construct
a gauge transformation $\sigma_x$ over $B_r(x)$ for any $x \in 
M\backslash N_{3r}\left((S\cap K)\cup E_r\right)$, such that
\begin{equation}
\label{eq:2.3.7}
\left|\sigma_x(A)\right|(y)\leqslant\frac{C}{r}
\left ( r^{4-n} \int _{B_r(x)} |F_A|^2 dV_g
\right )^{\frac{1}{2}},
~~\forall y \in B_r(x).\hskip.5in  \spn{2.3.7}
\end{equation}
Note that for any $\delta > 0$ and any subset $S' \subset M$,
$$
N_{\delta}( S')=\{x\in M~|~ d(x, S')\le \delta \},
$$
where $d(\cdot, \cdot)$ denotes the distance function of the metric $g$.

Gluing these $\sigma_x$ appropriately, we can construct
a gauge transformation\break $\sigma_k$ over each $B_{8r}(x_k)\backslash
N_{3r}\left((S\cap K)\cup E_r\right)$, such that
\begin{equation}
\label{eq:2.3.8}
\left|\sigma_k(A)\right|(x)\leqslant\frac{C}{r}
\left ( r^{4-n} \int _{B_r(x)} |F_A|^2 dV_g
\right )^{\frac{1}{2}},  \spn{2.3.8}
\end{equation}
whenever $x\in B_{8r}(x_k)\backslash N_{3r}((S\cap K)\cup E_r)$.
One gets from (2.3.8) that on the overlap 
$\left (B_{8r}(x_l)\cap B_{8r}(x_k)\right )\backslash N_{3r}\left((S\cap K)\cup E_r\right)$,
$$
|d \sigma_k\cdot \sigma_l^{-1}| \le \frac{2 C \sqrt{\varepsilon}}{r} .
$$
Hence, by modifying $\sigma_k$ slightly on overlaps, we may
assume that $\sigma_k\cdot\sigma^{-1}_{l}$ is constant on each
connected component of $B_{8r}(x_k)\cap B_{8r}(x_l)\backslash
N_{3r}\left((S\cap K)\cup E_r\right)$ for any $k\not= l$.

Let $\eta\!:\!\RR^1\!\to\!\RR^1$ be a cut-off $C^{\infty}$-function
satisfying: $\eta(t)=0$ for $t~\leqslant~1$, $\eta(t)=1$ for
$t\geqslant 2$, and $0\leqslant\eta^\prime(t)\leqslant 1$. Then
\begin{equation}
\label{eq:2.3.9}
 \int_{M}d\varphi\wedge\tr(F_A\wedge F_A)=\lim_{r\to 0}\int_{M}\eta\left(\frac{d(x,(S\cap K)\cup E_{r})}{3r}\right)
d\varphi\wedge\tr(F_A\wedge F_A). \spn{2.3.9}
\end{equation}  
For each $k\leqslant L^{\prime}_{r}$, 
\begin{eqnarray}
\label{eq:2.3.10}
\spn{2.3.10}\\
 \tr\left(F_A\wedge F_A\right)(x) 
&=&\tr\left(F_{\sigma_{k}(A)}\wedge
  F_{\sigma_{k}(A)}\right)(x)\nonumber \\ [2ex]
&=&d\tr\left(\sigma_{k}(A)\wedge
  F_{\sigma_{k}(A)}+\frac{1}{3}\sigma_{k}(A)
  \wedge\sigma_{k}(A)\wedge\sigma_{k}(A)\right)(x),\nonumber 
\end{eqnarray}
where $x\in B_{8r}(x_{k})\backslash N_{3r}\left((S\cap K)\cup E_{r}\right)$.

Since $\sigma_{k}\cdot\sigma_{l}^{-1}$ is piecewise constant,  we have
\begin{eqnarray}
&& \hskip-7pt\tr \left(\sigma_k(A)\wedge
F_{\sigma_{k}(A)}+\frac{1}{3}\sigma_k(A)\wedge\sigma_k(A)\wedge
\sigma_k(A)\right)\nonumber\\
= && \hskip-7pt \tr\left(\sigma_l(A)\wedge F_{\sigma_{l}(A)}+ 
\frac{1}{3}\sigma_l(A)\wedge\sigma_l(A)\wedge\sigma_l(A)\right)\nonumber
\end{eqnarray}  
on the overlap $B_{8r}(x_k)\cap B_{8r}(x_l)\backslash
N_{3r}\left((S\cap K)\cup E_r\right)$. Therefore, there is  a
globally defined Chern-Simon transgression form $\Psi$ outside
$N_{3r}(S\cup E_r)$,
such that
$$
d\Psi = \tr(F_A\wedge F_A)
$$
and
$$
\Psi(x)=\tr\left(\sigma_k(A)\wedge F_{\sigma_k(A)}+\frac{1}{3}\sigma_k(A)
\wedge\sigma_k(A)\wedge\sigma_k(A)\right)
$$
whenever $x\in B_{8r}(x_k)$.
For each $k$ and any $x\in B_{6r}(x_k)\backslash B_{3r}(x_k)$,
$$
|\psi(x)| \le C r^{-3} \left ( r^{4-n} \int_{B_r(x)} |F_A|^2 dV_g\right )^{\frac{3}{2}}
\le C r^{1-n} \int_{B_{8r}(x_k)} |F_A|^2 dV_g.
$$
It follows that 
\begin{eqnarray}
 \left|\int_{M}d\varphi\wedge\tr(F_A\wedge F_A)\right| 
&=&\lim_{r\to 0}\left|\int_{M}\eta\left(\frac{d(x, (S\cap K)\cup
      E_r)}{3r}\right)d\varphi\wedge d\Psi\right|\nonumber\\
&\leq&\lim_{r\to 0}\int_{3r\leqslant d(x, (S\cap K)\cup
  E_r)\leqslant 6r}\frac{1}{3r}|\Psi||d\varphi|\,dV_g\nonumber\\
&\leq&C \lim_{r\to
  0}\left\{\sup_{M}|d\varphi|\sum_{k=1}^{L^{\prime}_{r}}
\int _{B_{8r}(x_k)} |F_A|^2 dV_g
\right\}\nonumber\\
&\le&C \sup_{M}|d\varphi| 
\lim_{r\to 0} \int_{N_{8r}(S\cup E_r)}|F_A|^2 dV_g. \nonumber
\end{eqnarray} 
Since $\bigcap _{r > 0} N_{8r}(S\cup E_r) \subset S$ and 
$N_{8r}(S\cup E_r) \subset N_{8r'}(S\cup E_{r'})$ for $r\le r'$,
the last integral converges to zero as $r$ tends to $0$.
Therefore, we have
$$
\int_{M}d\varphi\wedge\tr(F_A\wedge F_A)=0,
$$
so that $\tr(F_A\wedge F_A)$ is closed in the sense of distribution.
\enddemo 

Let $C_1$ and $C_2$ denote the Chern-Weil polynomials defining
the first two Chern classes. Then $C_1(A)={\rm Ch}_1(A)$ is well-defined.

On $M\backslash S(A)$,  
\begin{equation}
\label{eq:2.3.11}
C_2(A)=\frac{1}{8\pi^2} \left (\tr (F_A\wedge F_A) - \tr(F_A)\wedge \tr(F_A)\right ).
\spn{2.3.11}
\end{equation}
Then $C_2(A)$ extends to a form, still denoted by $C_2(A)$,
on $M$ in the sense of distribution.

\specialnumber{2.3.2} \proclaim{{C}orollary}
\label{cor:2.3.2}
The extended form $C_2(A)$ is closed{\rm .}
\endproclaim
\demo{Proof}
Since $\tr(F_A)$ is harmonic outside $S(A)$ and $L^2$-bounded,
by the standard elliptic theory, it extends to be a smooth form
on $M$. Then this corollary follows from the  last proposition.
\enddemo \enddemo

\section{Rectifiability of blow-up loci}
\label{chap:3}

We study the blow-up set of Yang-Mills connections which converge to an 
admissible Yang-Mills connection.

\demo{{\rm 3.1.}\ Convergence of Yang-Mills connections}
Given any sequence of admissible Yang-Mills connections
${A_i}$,
we say that the  $A_i$ converge weakly to an admissible Yang-Mills connection $A$ 
(modulo gauge transformations), if $\int_M |F_{A_i}|^2 dV_g \le c$ for some
uniform constant $c$ and there are a closed subset $S$
and gauge transformations $\sigma_i$ of the $G$-bundle $E$ over
$M\backslash S$, such that for any compact 
$K \subset M\backslash S$,  
$\sigma_i(A_i) $ extend smoothly across $K$ for $i$ sufficiently large
and converge to $A$ in the $C^\infty$-topology in $K$ as $i$ tends to infinity. 
Obviously, $S$ contains $S(A)$.
In particular, this implies that for any smooth form $\varphi$ with compact support in $M$,
$$
\lim_{i\to \infty} \int _M (F_{\sigma_i(A_i)} , d\varphi ) dV_g 
= \int_M (F_A, d\varphi ) dV_g.
$$
This is exactly what the weak convergence is.
Clearly, we have the next result:
\enddemo
\specialnumber{3.1.1} \proclaim{Lemma}
\label{lem:3.1.1}
Weak limits of admissible connections $\{A_i\}$ are unique
modulo gauge transformations{\rm .}
\endproclaim

From now on, we always assume that $\{A_i\}$
is a sequence of smooth Yang-Mills connections
with $YM(A_i) \le \Lambda$. All the discussions in this section
also work for admissible Yang-Mills connections with slight
modification.

\specialnumber{3.1.2} \proclaim{Proposition}
\label{prop:3.1.2}
There is a subsequence
$\{A_{i_{j}}\}$ which converges weakly to some admissible
Yang\/{\rm -}\/Mills connection $A$ on $M${\rm .}
\endproclaim

\demo{Proof} Let $\varepsilon$ be as in Theorem 2.2.1
and $a$ be as in Theorem 2.1.2. We define a
closed subset for each $i$ and $r>0$ sufficiently
small:
\begin{equation}
\label{eq:3.1.3}
E_{i,r}=\{x\in M \,| \, e^{a r^2} r^{4-n}\int_{B_r(x)}
|F_{A_i}|^2dV_g\geq\varepsilon\}.\spn{3.1.1}
\end{equation}

It follows from the monotonicity formula
(Theorem 2.1.2) that $E_{i,r} \subset E_{i, r'}$ for any $r \le r'$.

By the standard diagonal process, we can choose a subsequence 
$\{i_j\}$ of $\{i\}$ such that for each $k$, the  $E_{i_j, 2^{-k}}$ converge
to a closed subset $E_{2^{-k}}$. Then $E_{2^{-k}} \subset E_{2^{-l}}$ for
$k\ge l$. Put $S= \bigcap _k E_{2^{-k}}$.

We first claim that $S$ is of Hausdorff codimension at least $4$.
Given $\delta >0$ sufficiently small and any compact subset $K$
of $M$, 
let $\{B_{4\delta}(x_\alpha)\}$ be any finite covering of
$S\cap K$ such that (1) $x_\alpha \in 
S\cap K$; (2) $B_{2\delta}(x_\alpha)
\cap B_{2\delta}(x_\beta) =\emptyset$
for $\alpha \not= \beta$. Take $k$ big enough such that $2^{-k} < \delta$. 
Then for $j$ sufficiently large,
there are $y_\alpha \in E_{i_j, 2^{-k}}$ such that $d(x_\alpha, y_\alpha) < \delta$.
Then $\{B_{5\delta }(y_\alpha)\}$ is a finite covering of $S\cap K$
and $B_{\delta}(y_\alpha)
\cap B_{\delta}(y_\beta) =\emptyset$ for $\alpha \not= \beta$.
By Theorem 2.1.2, 
$$
e^{a\delta ^2}\delta^{4-n}\int_{B_{\delta}(y_\alpha)}
|F_{A_{i_j}}|^{2}dV_g \ge 
e^{a 2^{-2k} }2^{(n-4)k}\int_{B_{2^{-k}}(y_{\alpha})}
|F_{A_{i_j}}|^{2}dV_g \geq\varepsilon.
$$
Hence,
$$
\sum_\alpha \delta^{n-4}\leq\frac{e^a}{\varepsilon}\sum_\alpha
\int_{B_{\delta}(y_\alpha)}|F_{A_{i_j}}|^2
dV_g\leq\frac{e^a}{\varepsilon}\int_M|F_{A_{i_j}}|^2 dV_g\leq
\frac{ce^a}{\varepsilon}.
$$
It follows that $H^{n-4}(S\cap K)$, and consequently,
$H^{n-4}(S)$, is no more than $\frac{5^{n-4} e^a c}{\varepsilon}$.
This proves the claim.

Now we prove that $A_{i_j}$ converges to some $A$ outside $S$
modulo gauge transformations. To save the notation, we assume
$\{i_j\}=\{i\}$. 

We notice that for any $r> 0$, there is an $i(r)> 0$, such that
for any $i \ge i(r)$ and $x\in M$ with $d(x, E_{2^{-k}}) \ge  r$, 
where $2^{-k-1} \le r \le 2^{-k}$,
\begin{equation}
\label{eq:3.1.4}
e^{ar^2}r^{4-n}\int_{B_r(x)}
|F_{A_{i}}|^{2}dV_g < \varepsilon. \spn{3.1.2}
\end{equation}
This is equivalent to saying that $x\in M\backslash E_{i,r}$.
By Theorem 2.2.1, we deduce from (3.1.2) that for any $x\in M\backslash
B_r(E_r)$,
\begin{equation}
\label{eq:3.1.5}
|F_{A_{i}}| (x) < {\frac{C\sqrt {\varepsilon}}{r^2}}. \spn{3.1.3}
\end{equation}
It follows from Theorem 3.6 in [Uh2] that there exists a subsequence
$\{i'\} \subset \{i\}$ and gauge transformations $\sigma(i')$, such that 
$\sigma(i')(A_{i'})$ converge to a smooth connection $A$ in
$C^1$-topology on any compact subset outside $S$. Since    $A_i$ are
Yang-Mills connections, by the standard elliptic theory, $A$
is a Yang-Mills connection and   $\sigma(i')(A_{i'})$ converge to $A$ smoothly outside $S$.
\enddemo 

In the following, we always assume 
that the sequence $A_i$ converges to an admissible
Yang-Mills connection $A$ with $\int_M |F_{A_i}|^2 dV_g\le \Lambda$.

\specialnumber{3.1.3} \proclaim{Lemma}
\label{lem:3.1.3}
Define
\begin{equation}
\label{eq:3.1.6}
 {\it S}_b(\{A_i\})=\bigcap_{r>0}\{x\in M|\lim_{i\to\infty}\inf  e^{a r^2}r^{4-n}\int_{B_r(x)}
|F_{A_i}|^2dV_g\geq\varepsilon\}, \hskip.25in \spn{3.1.4}
\end{equation}
where $\varepsilon$ is as  given in Theorem~{\rm \ref{th:2.2.1}.}
Then {\rm (i)} ${\it S}_b(\{A_i\})$ is closed and contained in the above $S${\rm ;}
{\rm (ii)} Its Hausdorff measure $H^{n-4}({\it S}_b(\{A_i\}))\leq C$ for some
constant $C$ depending only on $M$ and $\Lambda${\rm ;
(iii)} $A$ extends to   a smooth connection on $M\backslash {\it S}_b(\{A_i\})${\rm .}
\endproclaim

\demo{Proof} Suppose $x_0\in M\backslash {\it S}_b(\{A_i\})$;
then there is an $r_0>0$ such that
$$
r_0^{4-n}\int_{B_{r_{0}}(x_0)}|F_{A_{n_{i}}}|^2 dV_g<\varepsilon
$$
for some subsequence $n_i\to\infty$.  By Theorem 2.2.1, 
$$
\sup_{n_{i}}\sup_{x\in B_{\frac{r_0}{2}}(x_0)}|F_{A_{n_{i}}}|
\leq\frac{c_0\sqrt{\varepsilon}}{r_0^2}
$$
for some constant $c_0=c_0(n,M)$.  In particular,
$$
\sup_{n_{i}}\sup_{x\in
  B_{\frac{r_0}{4}}(x_0)}r^{4-n}\int_{B_{r}(x)}|F_{A_{n_{i}}}|^2
dV_g
\leq\frac{\varepsilon}{2}
$$
whenever $r\leq r_0/\sqrt[4]{c+1}$ for some 
constant $c$ depending only on $g$. Hence,
$B_{\frac{r_0}{4}}(x_0)\subset M\backslash {\it S}_b(\{A_i\})$,
and consequently, ${\it S}_b(\{A_i\})$ is closed. This also implies
that $A$ is a limit of some subsequence of $\{A_{n_i}\}$ 
(modulo gauge transformations) in $B_{\frac{r_0}{4}}(x_0)$
in the  $C^\infty$-topology. Then (iii) follows. 

For any $x_0 \in M\backslash S$, if $r$ is sufficiently small,

$$
r^{4-n} \int_{B_r(x_0)} |F_A|^2 dV_g < \varepsilon_0.
$$
\vglue4pt\noindent This implies that for $i$ sufficiently large,

$$
r^{4-n} \int_{B_r(x_0)} |F_{A_i}|^2 dV_g < \varepsilon_0.
$$
\vglue4pt\noindent Hence, $x_0 \in M\backslash {\it S}_b(\{A_i\})$. This shows that
${\it S}_b(\{A_i\}) \subset S$.

The estimate on $H^{n-4}({\it S}_b(\{A_i\}))$ follows from
the proof of Proposition 3.1.2.
\enddemo

Since $A$ can be extended smoothly to $M\backslash {\it S}_b(\{A_i\})$,
we may assume that $S(A) \subset {\it S}_b(\{A_i\})$.
If ${\it S}_b(\{A_i\}) = \emptyset$, then there is a subsequence of $\{A_i\}$
which converges to $A$ smoothly on $M$.

Consider the Radon measures 
$\mu_i=|F_{A_{i}}|^2 dV_g$ $(i=1,2,\cdots )$. By taking a subsequence
if necessary, we may assume that $\mu_i\to\mu$ weakly on $M$ as Radon measures; i.e.,
for any continuous function $\varphi$ with compact support in $M$,
\begin{equation}
\label{eq:3.1.7}
\lim_{i\to \infty} \int_M \varphi |F_{A_i}|^2  dV_g = \int _M \varphi d\mu. \spn{3.1.5}
\end{equation}

\vglue4pt
Let us write (by Fatou's lemma)
$$
\mu=|F_A|^2 dV_g+\nu
$$
for some nonnegative Radon measure $\nu$ on $M$.

\specialnumber{3.1.4} \proclaim{Lemma}
\label{lem:3.1.4}
When $\nu(x)=\Theta(x)H^{n-4}\lfloor {\it S}_b(\{A_i\}), x\in M${\rm ,} 
for $H^{n-4}${\rm -a.e.} $x\in {\it S}_b(\{A_i\})${\rm ,} then
$$
\varepsilon\leq \Theta(x)\leq\ 4^{n-4}
r_x^{4-n} e^{a r_x^2}\Lambda ,
$$
where $r_x${\rm ,} $a$ are as given in Theorem {\rm 2.1.2.}
\endproclaim
\demo{Proof}
First we observe:
\begin{itemize}
\item[(a)] For any $x\in M$, $e^{a r^2} r^{4-n}\mu(B_r(x))$ is a 
  nondecreasing function of $r$ sufficiently small; thus the density 
\end{itemize}
\begin{equation}
\label{eq:3.1.8}
\Theta(\mu,x)=\lim_{r\to 0+}r^{4-n}\mu(B_r(x)) \spn{3.1.6}
\end{equation}
\begin{itemize} \item[]
exists for every $x\in M$;
\item[(b)] $x\in {\it S}_b(\{A_i\})\mbox{~if and only if~}\Theta(\mu,x)\geq\varepsilon$;
\item[(c)] For $H^{n-4}$-a.e. $ x\in {\it S}_b(\{A_i\})$,
$$
\lim_{r\to 0+}r^{4-n}\int_{B_r(x)}|F_A|^2dV_g=0.
$$
\end{itemize}
Indeed, (a) follows from the monotonicity formula in
Section 2.1. The statement (b) follows from the
definition of ${\it S}_b(\{A_i\})$ and (a).  To prove (c), we define
\begin{equation}
\label{eq:3.1.9}
E_j=\{x~|~\overline{\lim}_{r\to 0+}r^{4-n}\int_{B_{r}(x)}|F_A
|^2dV_g>\frac{1}{j}\}. \spn{3.1.7}
\end{equation}
It suffices to show that $H^{n-4}(E_j)=0$ for each $j\geq1$.  For
any $\delta>0$, there is a covering of $E_j$ by balls
$B_{2r_{\alpha}}(x_\alpha)$ with $x_\alpha\in E_j$ and
$2r_{\alpha}\leq\delta$,
such that
$$
r_{\alpha}^{4-n}\int_{B_{r_{\alpha}}(x_\alpha)}
|F_A|^2 dF_g>\frac{1}{j},
$$
and $B_{r_{\alpha}}(x_{\alpha})\cap
B_{r_{\beta}}(x_\beta)=\emptyset$. Then
$$H^{n-4}(E_j)-\psi(\delta)\leq\sum_{\alpha}(2r_{\alpha})^{n-4}\leq
j2^{n-4}\int_{N_{\delta}({\it S}_b(\{A_i\}))}|F_A|^2 dV_g,$$
where $N_\delta({\it S}_b(\{A_i\}))$ denotes the $\delta$-tubular neighborhood of
${\it S}_b(\{A_i\})$, and\break $\psi(\delta)\to 0$ as $\delta\to0$.  It follows that
$H^{n-4}(E_j)=0$, since $\delta$ is arbitrarily small.

From the monotonicity formula we obtain
\begin{equation}
\label{eq:3.1.10}
r^{4-n}\mu(B_{r}(x))\leq C \spn{3.1.8}
\end{equation}
for some constant $C$ depending only on $\Lambda$ and $M$.
Therefore, $\mu |_{{\it S}_b(\{A_i\})}$ is absolutely continuous with respect to
$H^{n-4}\lfloor {\it S}_b(\{A_i\})$; consequently, by the Radon-Nikodym theorem, we have
\begin{equation}
\label{eq:3.1.11}
\mu |_{{\it S}_b(\{A_i\})}(x)=\Theta(x) H^{n-4}\lfloor {\it S}_b(\{A_i\}) \spn{3.1.9}
\end{equation}
for $H^{n-4}$-a.e. $x\in {\it S}_b(\{A_i\})$. Then by (c),  
$$
\nu(x)=\Theta(x)H^{n-4}\lfloor {\it S}_b(\{A_i\})
$$
for $H^{n-4}$-a.e. $x\in {\it S}_b(\{A_i\})$. 
We notice that $\mu$ is a Borel regular measure.
The estimates of $\Theta(x)$ follow the
above density estimate (b) and the fact that for
$H^{n-4}$-a.e. $x\in {\it S}_b(\{A_i\})$,
\begin{equation}
\label{eq:3.1.12}
2^{4-n} \leq
\overline {\lim}_{r\to 0+}\frac {H^{n-4}({\it S}_b(\{A_i\})\cap B_{r}(x))}{r^{n-4}}
\leq 1, \spn{3.1.10}
\end{equation}
which can be easily proved (cf.\ [Si, Th.~3.6)].
\enddemo

Define
\begin{equation}
\label{eq:3.1.13}
S_b = \overline {\{x\in {\it S}_b(\{A_i\})\,|\, \Theta(\mu, x) > 0,~
\lim_{r\to 0+}r^{4-n}\int_{B_r(x)}|F_A|^2dV_g=0\}}. \quad\spn{3.1.11}
\end{equation}
Then ${\it S}_b(\{A_i\})= S_b \cup S(A)$. We call $(S_b, \Theta)$ the blow-up
locus of the weakly convergent sequence $\{A_i\}$. Here, $S_b$ is the support of
the blow-up locus and $\Theta$ is its multiplicity. If no confusion can occur,
we may simply say that $S_b$ is the blow-up locus.
 
\demo{{\rm 3.2.}\ Tangent cones of blow-up loci}
We adopt the notation of the  last section unless specified otherwise. For simplicity, we write
$S=S_b$ for the blow-up locus.  
In this section, we study the properties of tangent cones of $S$.

Recall that $\mu$ is the limit Radon measure of $\mu_{i}=|F_{A_{i}}|^2dV_g$.  
For any $y\in M$ and sufficiently small $\lambda$, we define the scaled
measure $\mu_{y,\lambda}$ as follows: for any $E$ in $T_yM$,
\begin{equation}
\label{eq:3.2.1}
\mu_{y,\lambda}(E)=\lambda^{4-n}\mu(\exp_y(\lambda E)), \spn{3.2.1}
\end{equation}
where $\exp_y:T_yM\to M$ is the exponential map of the metric $g$
and
\begin{equation}
\label{eq:3.2.2}
\lambda E =\{x \in T_yM\,|\,\lambda^{-1}x\in E\}. \spn{3.2.2}
\end{equation} \enddemo

\specialnumber{3.2.1} \proclaim{Lemma}
\label{lem:3.2.1}
Let $\{\lambda_k\}$ be any sequence with
$\displaystyle{\lim_{k\to\infty}\lambda_k=0}${\rm . } 
Then there exist  a subsequence
$\{\lambda^{\prime}_{k}\}$ and a Radon measure $\eta$ on $T_yM$ such
that $\mu_{y,\lambda^{\prime}_{k}}$ converges to $\eta$ weakly{\rm .}
Moreover{\rm ,} $\eta_{0,\lambda}=\eta$ for each $\lambda>0${\rm ;} i{\rm .}e{\rm .,}
$\eta$ is a cone measure{\rm .}
\endproclaim

\demo{Proof}
Define a connection on $T_yM$ for each $y$ and $\lambda$ by
\begin{equation}
\label{eq:3.2.3}
A_{i,y,\lambda}= \tau _\lambda ^\ast \exp_y^\ast A_i, \spn{3.2.3}
\end{equation}
where $\tau _\lambda : T_yM \mapsto T_yM$ maps $v$ to $\lambda v$.
Then $A_{i,y,\lambda}$ is a Yang-Mills connection with respect
to the metric $\lambda^{-2}\exp^{\ast}_{y}g$, which will be denoted by
$g_{y,\lambda}$.  Clearly, $g_{y,\lambda}$
converges to the flat metric $g_{y,0}= g|_{T_yM}$ on $T_yM$ as $\lambda\to 0+$.
Moreover, by the monotonicity (Theorem 2.1.2), for any small $r > 0$,  
\begin{eqnarray}
\label{eq:3.2.4}
&&\hskip-.25in e^{a \lambda^2 r^2} r^{4-n}
\int_{B_r(0,g_{y,\lambda})}|F_{A_{i,y,\lambda}}|^2
dV_{y,\lambda} \spn{3.2.4}\\
&&\qquad\quad=
e^{a \lambda^2 r^2}(\lambda r)^{4-n}\int_{B_{\lambda r}(y)}|F_{A_{i}}|^{2}dV_g\leq
C(M,\Lambda), \nonumber
\end{eqnarray}
where $C(M,\Lambda)$ denotes a constant depending only on
$M,\Lambda$, and
$B_r(0,g_{y,\lambda})$ denotes the geodesic
ball of $g_{y,\lambda}$ in $T_yM$ with radius $r$ and center
at $0$.  Clearly, $|F_{A_{i,y,\lambda}}|^2 dV_{y,\lambda}$
converges to $\mu_{y,\lambda}$ weakly on
$B_{\lambda^{-1}r_0}(0,g_{y,\lambda})$, where $r_0$ depends
only on $y$. Letting $i$ go to   infinity, we obtain
\begin{equation}
\label{eq:3.2.5}
\mu_{y,\lambda}(B_{r}(0,g_{y,\lambda}))\leq C(M,\Lambda)r^{n-4} \spn{3.2.5}
\end{equation}
for any $r\leq\lambda^{-1}r_0$.  Hence, we may
find a subsequence 
$\{\lambda^{'}_{k}\}\subset\{\lambda_k\}$ such
that  the $\mu_{y,\lambda^{'}_{k}}$ converge to $\eta$ weakly as Radon
measures on $T_yM$. Then there are (by the standard diagonal
process) $i_k\to\infty$ such that
\begin{equation}
\label{eq:3.2.6}
|F_{A_{i_{k},y,\lambda_{k}}}|^2 dV_{y,\lambda}\to\eta. \spn{3.2.6}
\end{equation}
Since $\mu$ is the weak limit of $\mu_i=|F_{A_{i}}|^2 dV_g$, 
 for $0<\sigma<\rho$ sufficiently small,
\begin{equation}
\label{eq:3.2.7}
e^{a\sigma^{2}}\sigma^{4-n}\mu(B_{\sigma}(y))\leq e^{a\rho^{2}}\rho^{4-n}
\mu(B_{\rho}(y)). \spn{3.2.7}
\end{equation}
This implies that $\lim_{r\to 0+}r^{4-n}\mu(B_{r}(y))=\Theta(\mu,y)$, 
and consequently, for any $r>0$,
\begin{eqnarray}
\label{eq:3.2.8}
 r^{4-n}\eta(B_{r}(0,g_{y,0})) 
&=&\lim_{\lambda^{\prime}_{k}\to
  0}r^{4-n}\mu_{y,\lambda^{\prime}_{k}}(B_{r}
(0,g_{y,\lambda_k',}))\spn{3.2.8}\\ 
&=& \lim_{\lambda^{'}_{k}\to
  0}(\lambda^{\prime}_{k}r)^{4-n}
\mu(B_{\lambda^{\prime}_{k}r}(y)).\nonumber
\end{eqnarray}
That is,
$$
\eta(B_{r}(0,g_{y,0}))\ =\ \Theta(\mu,y)r^{n-4}.  \hskip.7in
$$ 
This indicates that $\eta$ is a cone measure. To prove it rigorously, 
we first observe that for any
$0<\sigma<\rho<\infty$,
\begin{equation}
\label{eq:3.2.9}
\int_{B_{\rho}(0,g_{y,\lambda^{\prime}_{k}})\backslash B_{\sigma}(0,g_{y,\lambda^{\prime}_{k}})}{r^{4-n}|
\frac{\partial}{\partial r}}\rfloor
F_{A_{i_{k},y,\lambda^{'}_{k}}}|^2dV_{y,\lambda^{\prime}_{k}}\to 0
\mbox{ as } k\to\infty. \qquad \spn{3.2.9}
\end{equation}
Here we have used Theorem 2.1.2 and (3.2.8).

Let $\phi(\theta)$ be any positive function on the unit sphere
$S^{n-1}\subset T_yM$.
It follows from Theorem 2.1.1 that for any
$0<\sigma<\rho$ and $\lambda^{\prime}_k$ sufficiently small,
\begin{eqnarray}
 &&\hskip-1in \sigma^{4-n}e^{a(\lambda^{\prime}_{k}\sigma)^2}
       \int_{B_{\sigma}(0,g_{y,\lambda^{\prime}_{k}})}
           |F_{A_{i_{k},y,\lambda^{\prime}_{k}}}|^2\phi (\theta)dV_{y,
            \lambda_{k^{\prime}}} \nonumber\\[2ex]
&\le&\rho^{4-n}e^{a(\lambda^{\prime}_{k}\rho)^2}\int_{B_{\rho}(0,g_{y,
         \lambda^{\prime}_{k}})}|F_{A_{i_{k},y,\lambda^{\prime}_{k}}}|^2
    \phi (\theta)
         dV_{y,\lambda^{\prime}_{k}}\nonumber\\[2ex]
&&- 4\int_{B_{\rho}(0,g_{y, \lambda^{\prime}_{k}})\backslash
B_\sigma(0,g_{y,\lambda^{\prime}_{k}})}r^{4-n}|\frac{\partial}{\partial
           r}\rfloor F_{A_{i_{k},y,\lambda^{\prime}_{k}}}|^2\phi dV_{y,
            \lambda^{\prime}_{k}}\nonumber\\[2ex]
&& - 4\int^{\rho}_{\sigma}\tau^{3-n}d\tau
          \int_{B_{\tau}(0,g_{y, \lambda^{\prime}_{k}})}|\nabla\phi||
            \frac{\partial}{\partial
           r}\rfloor F_{A_{i_{k},y,\lambda^{\prime}_{k}}}|^2 
dV_{y, \lambda^\prime_k}.\nonumber
\end{eqnarray} 
Letting $k$ go to infinity, by (\ref{eq:3.2.9}), we obtain
\begin{equation}
\label{eq:3.2.15}
\sigma^{4-n}\int_{B_{\sigma}(0,g_{y,0})}\phi
d\eta=\rho^{4-n}\int_{B_{\rho}(0,g_{y,0})}\phi d\eta. \spn{3.2.10}
\end{equation}
Differentiating (\ref{eq:3.2.15}) on $\rho$, we have
\begin{equation}
\label{eq:3.2.16}
\rho^{n-4}\int_{\partial B_{\rho}(0,g_{y,0})}\phi
d\xi=(n-4)\int_{B_{\rho}(0,g_{y,0})}\phi d\eta, \spn{3.2.11}
\end{equation}
where $d\eta(r,\theta)=r^{n-5}drd\xi(r,\theta)$. Combining
(\ref{eq:3.2.15}) and (\ref{eq:3.2.16}), we get
\begin{equation}
\label{eq:3.2.17}
\int_{\partial
  B_{\rho}(0,g_{y,0})}\phi(\theta)d\xi(\rho,\theta)=
        \int_{\partial B_{\sigma}(0,g_y,0)}\phi(\theta)d\xi(\sigma,\theta), \hskip.75in\spn{3.2.12}
\end{equation}
for any $0<\sigma<\rho<\infty$.  This implies
$$
drd\xi(r+r_1,\theta)=drd\xi(r,\theta)
$$
for any $r_1>0$.  That is, $r^{5-n}d\eta(r,\theta)$ is
translation invariant in $r$, or
$d\eta(r,\theta)=r^{n-5}dr d\xi(\theta)$ for some
Radon measure $d\xi(\theta)$ on $S^{n-1}$.
\enddemo

Next we study the tangent cones $\eta $ with support in $T_yM$ at
$H^{n-4}$-a.e. $y\in S$.
First we recall two elementary lemmas about the Radon measure
$\mu$ given above.

\specialnumber{3.2.2} \proclaim{Lemma}
\label{lem:3.2.2}
The density function $\Theta(\mu,x)$ is $H^{n-4}$\/{\rm -}\/approximately
continuous at $H^{n-4}${\rm -a.e.} $x$ in $S${\rm .} Here $\Theta(\mu, \cdot )$ is
$H^{n-4}$\/{\rm -}\/approximately continuous at $ x\in S$ if for any $\varepsilon >0${\rm ,}
\begin{equation}
\label{eq:3.2.18}
\lim_{r\to 0}\frac{H^{n-4}(\{y\in B_r(x)\cap S\,|\,|\Theta(\mu,y)-\Theta(\mu,x)|>\varepsilon\})}{r^{n-4}}=0.\qquad
\spn{3.2.13}
\end{equation}
\endproclaim
 
\demo{Proof}
The density function
$\Theta(\mu,x)(x\in S)$ is upper-semi-continuous, so that 
$E_c = \{ x \, |\, \Theta(\mu,x) < c\}$ is open, and consequently,
for any $c_1 < c_2$, $E_{c_2}\backslash E_{c_1}$ is a Borel
set and thus measurable. Now we define
$$
E_i = \{ x\in S ~|~\frac{(i-1)\varepsilon }{2}
\le \Theta (\mu, x) < \frac{i \varepsilon}{2}\}.
$$
Clearly, each $E_i$ is contained in $S$ and $H^{n-4}(S\backslash \bigcup_i E_i)=0$.
Then for any $x\in E_i$, we have
\begin{eqnarray*}
&&\hskip-1in \lim_{r\to 0}\frac{H^{n-4}(\{y\in B_r(x)\cap 
S \,|\,|\Theta(\mu,y)-\Theta(\mu,x)|>\varepsilon\})}{r^{n-4}}\\ \noalign{\vskip6pt}
&&\qquad =\ \overline {\lim}_{r\to 0}\frac{H^{n-4}(B_r(x)\cap 
(S\backslash E_i ))}{r^{n-4}} 
\ = \ 0.
\end{eqnarray*}
Here we have used Theorem 3.5 in [Si2].
Thus the lemma follows.
\enddemo 

\specialnumber{3.2.3} \proclaim{Lemma}
\label{lem:3.2.3}
Let $x\in S$ be such that $\Theta(\mu,x)\geq\varepsilon_0>0$ and
$\Theta(\mu,\cdot)$ is $H^{n-4}$\/{\rm -}\/approximately continuous at
$x${\rm .}  Then there is a $r_x> 0${\rm ,} such
that for each $r\in (0,r_x)${\rm ,} we may find $n-4$ points $x_1,\ldots
,x_{n-4}$ in $B_r(x)\cap S$ satisfying\/{\rm :}\/
\begin{itemize}
\item[{\rm (i)}]
$\Theta(\mu,x_j)\geq\Theta(\mu,x)-\varepsilon(r)$
for $j=1,2,\cdots,n-4${\rm ,} where $\varepsilon(r)\to 0$
as $r$ tends to zero\/{\rm ;}\/
\item[{\rm (ii)}]
Let $\exp_x$ be the exponential map of $(M,g)$ at $x${\rm .}  Then for
some\break $s\in(0,\frac{1}{2})$ depending only on $n${\rm ,} $d(x_1,x)\geq
sr$ and $d(x_k,\exp_x(V_{k-1}))\geq sr$\break for $k\ge 2${\rm ,} where
$V_{k-1}$ denotes the subspace in $T_xM$ spanned by\break
$(\exp_x|_{B_{r(0)}})^{-1}(x_1),\ldots, (\exp_x|_{B_{r}(0)})^{-1}(x_{k-1})${\rm .}
\end{itemize}

\endproclaim

\demo{Proof} The arguments here are essentially due to F.H. Lin in [Li].
By the assumption, we may find a positive function $\varepsilon(r)$
for $0<r<r_x$ such that $\lim_{r\to 0}\varepsilon(r)=0$ and
\begin{equation}
 \label{eq:3.2.19}
\frac{H^{n-4}(\{y\in B_r(x)\cap S~|~|\Theta(\mu,y)-\Theta(\mu,x)|\geq\varepsilon(r)\})}{r^{n-4}}
\leq\frac{s(n)}{2}<\frac{1}{2} \enspace\spn{3.2.14}
\end{equation}
($s(n)>0$ will be determined later).

Suppose that the lemma is false.  Then there would be a
sufficiently small $r>0$ such that one cannot find $n-4$ points
$x_1,\cdots ,x_{n-4}$ inside the set
\begin{equation}
\label{eq:3.2.20}
\{y\in S\cap B_r(x)\,|\, |\Theta(\mu,y)-\Theta(\mu,x)|<\varepsilon(r)\} \spn{3.2.15}
\end{equation}
satisfying   condition (ii) of Lemma 3.2.3.  Therefore, the set
in (3.2.15) is contained in an $sr$-neighborhood of $\exp_x(L)$ for
some $(n-5)$-dimensional subspace $L$ in $T_xM$. In particular,
this  implies
\begin{eqnarray}
&&   \label{eq:3.2.21}
\hskip-36pt\mu(\{y\in B_r(x)\cap S \,|\,
  |\Theta(\mu,y)-\Theta(\mu,x)|<\varepsilon(r)\}) \spn{3.2.16}\\ 
&&\qquad\leq \  C(n)s(n)r^{n-4}\Theta(\mu,x), \nonumber
\end{eqnarray}
where $C(n)$ is some uniform constant independent of $s(n)$.  

On the other hand, by the upper-semi-continuity of
$\Theta(\mu,\cdot)$, we may assume that for any
$y\in B_r(x)$,
$$
\Theta(\mu,y)\leq 2\Theta(\mu,x).
$$

Thus
\begin{eqnarray}
&&\quad \mu(\{y\in B_r(x)\cap
  S\,|\,|\Theta(\mu,y)-\Theta(\mu,x)|\geq\varepsilon(r)\}) \spn{3.2.17}
\label{eq:3.2.22}
\\ 
&&\qquad\quad \leq  2\Theta(\mu,x)H^{n-4}(\{y\in B_r(x)\cap
        S||\Theta(\mu,y)-\Theta(\mu,y)|\geq\varepsilon(r)\})
        \nonumber\\ 
&&\qquad\quad \leq \Theta(\mu,x)s(n)r^{n-4}. \nonumber
\end{eqnarray}
Putting (\ref{eq:3.2.21}) and (\ref{eq:3.2.22}) together, we
obtain
\begin{eqnarray}
\label{eq:3.2.23}
\mu(B_r(x)\cap S)&\leq&
s(n)\left (1+C(n)\right)\Theta(\mu,y)r^{n-4}\spn{3.2.18}\\ 
                 &<&\frac{1}{2}\Theta(\mu,x)r^{n-4}, \nonumber
\end{eqnarray}
if we choose $s(n)<\frac{1}{2(1+C(n))}$.  However,
(\ref{eq:3.2.23}) is impossible for $r$
sufficiently small, since $\lim_{r\to 0}\frac{\mu(B_r(x)\cap
  S)}{r^{n-4}}= \Theta(\mu,x) > 0$.
\enddemo

Now we can state the main result of this section, which will be
used in proving the rectifiability of blow-up loci.
\specialnumber{3.2.4} \proclaim{Proposition}
\label{prop:3.2.4}
Let $\mu$ be the Radon measure given at the beginning of this
section{\rm .} Then for $H^{n-4}${\rm -}a{\rm .}e{\rm .}\ $x\in S\subset M${\rm ,} any
tangent cone measure $\eta$ on $T_xM$ of $\mu$ is of the form
$\Theta(\mu,x)H^{n-4}\lfloor F$ for some $(n-4)$\/{\rm -}\/dimensional subspace
$F$ in $T_xM${\rm .}
\endproclaim

Note that the existence of $\eta$ is assured by
Lemma~\ref{lem:3.2.1}. 

The rest of this section is devoted to proving this proposition.
First we recall (cf.\ Lemma~\ref{lem:3.1.4}) that
$\mu=|F_A|^2dV_g+ \Theta (\mu, \cdot) H^{n-4}\lfloor S$, 
where $A$ is the weak limit of a sequence
$\{A_i\}$. By Lemma~\ref{lem:3.1.4}
and the observation (c) in its proof, for
$H^{n-4}$-a.e. $ x\in S$, 
\begin{equation}
\label{eq:3.2.24}
\Theta(\mu,x)\geq\varepsilon_0>0,\quad\lim_{r\to
  0}r^{4-n}\int_{B_r(x)}|F_A|^2dV_g=0 .\hskip.5in\quad \spn{3.2.19}
\end{equation}
Furthermore, it follows from Lemma~\ref{lem:3.2.2} that
$\Theta(\mu,\cdot)$ is $H^{n-4}$-approximately continuous 
at $H^{n-4}$-a.e. $ x$ in $S$.

From now on, we fix a point $x\in S$ such that (\ref{eq:3.2.24})
holds and $\Theta(\mu,\cdot)$ is $H^{n-4}$-approximately continuous 
at $ x$.

Assume that $\eta$ is the weak limit of $\mu_{x,r_k}$, where
$\lim_{k\to\infty}r_k = 0$. For $k$ sufficiently large, by
Lemma~\ref{lem:3.2.3}, we may find $n-4$ points
$x^k_1,\cdots,x^k_{n-4}$ in $B_{r_{k}}(x)\cap S$, such that
for $j=1,2,\cdots, n-4$,
\begin{equation}
\label{eq:3.2.25}
\Theta(\mu,x^{k}_{j})\geq\Theta(\mu,x)-\varepsilon(r_k), \spn{3.2.20}
\end{equation}
\begin{equation}
\label{eq:3.2.26}
d(x^{k}_{j},\exp_x(V^{k}_{j-1}))\geq sr_k, \spn{3.2.21}
\end{equation}
where $V^{k}_{j-1}$ denotes the
$0$-dimensional space $\{0\}$ if $j=1$,
and the subspace in $T_xM$ spanned by
$\xi^{k}_{1}=
\exp_{x}^{-1}(x^{k}_{1})$,...,
$\xi^{k}_{j-1}=\exp^{-1}_{x}(x^{k}_{j-1})$ for $j\geq 2$.

As before, we denote by $g_{x,r_{k}}$ the scaled metric
$r^{-2}_{k}\exp^{\ast}_{x}g$ on $T_xM$, which converges to the flat
one $g_{x,0}$ as $k$ tends to $\infty$.  Clearly, 
$r_k^{-1}\xi^{k}_{j}\in
B_1(0,g_{x,r_{k}})$ for each $j$, so that by taking a subsequence of
$\{r_{k}\}$ if necessary, we may assume that as $k$ tends
to $\infty$, $r^{-1}_{k}\xi^{k}_{j}\in
B_1(0,g_{x,0})$ converges to $\xi_j$ with respect to a fixed
metric $g_{x,0}$.  By (\ref{eq:3.2.26}), $\xi_{1},\cdots
,\xi_{n-4}$ span an $(n-4)$-dimensional subspace $F$ in $T_xM$,
which is in fact the limit of $V^{k}_{n-4}$.  Moreover,
$d_{g_{x,0}}(\xi_{i},0)\geq s$ and $d_{g_{x,0}}(\xi_{i},\xi_{j})\geq s$ for any $i\not= j$.

From (\ref{eq:3.2.25}), we can deduce that for any $r > 0$,
\begin{eqnarray}
  r^{4-n}\mu_{x,r_{k}}(B_r(\xi^{k}_{j},g_{x,r_{k}}))
   &=&(rr_k)^{4-n}\mu(B_{rr_{k}}(x^{k}_{j}))\nonumber\\[2ex]
&\geq&\Theta(\mu,x^{k}_{j}) \geq  \Theta(\mu,x)-\varepsilon(r_{k}).\nonumber
\end{eqnarray}
Thus for all $r<0$,
\begin{equation}
\label{eq:3.2.27}
r^{4-n}\eta(B_r(\xi_{j},g_{x,0}))\geq\Theta(\mu,x)=\Theta(\eta,0). \spn{3.2.22}
\end{equation}
In particular,
\begin{equation}
\label{eq:3.2.28}
\Theta(\eta,\xi_j)\geq\Theta(\eta,0). \spn{3.2.23}
\end{equation}
On the other hand, for any $r, \tilde r >0$, it follows from the monotonicity, 
\begin{eqnarray}
 r^{4-n}\eta(B_r(\xi_j,g_{x,0})) 
&=&\lim_{k\to\infty}r^{4-n}\mu_{x,r_{k}}(B_r(\xi^{k}_{j},g_{x,r_{k}}))
\nonumber\\[2ex]
&=&\lim_{k\to\infty}rr_{k}^{4-n}\mu(B_{rr_{k}}(x^{k}_{j}))\nonumber\\[2ex]
&\leq&\lim_{k\to\infty}\left(e^{a\tilde{r}^{2}}\tilde{r}^{4-n}
\mu(B_{\tilde{r}}(x^{k}_{j}))\right)\nonumber\\[2ex]
&=&e^{a\tilde{r}^{2}}\tilde{r}^{4-n}\mu(B_{\tilde{r}}(x)).\nonumber
\end{eqnarray}
Since $\tilde{r}$ can be arbitrarily small,  
$$
r^{4-n}\eta(B_r(\xi_{j},g_{x,0}))=\Theta(\eta,0)
$$
for any $r> 0$. Then, using Theorem 2.1.1 as  in the proof of
Lemma~\ref{lem:3.2.1}, we can show that $\eta$ is a cone measure with
center at $\xi_{j}$ for each $j= 1, \cdots , n-4$; i.e.,
$$
d\eta(r_j,\theta)=r^{n-5}_{j}dr_{j}d\xi(\theta)
$$
for some Radon measure $d\xi_{j}(\theta)$ on the unit sphere
$\{\xi\in T_{x}M|r_{j}(\xi )=1\}$, where $r_{j}(\xi )=|\xi-\xi_{j}|$.
Clearly, it follows that
$$
\eta(y_1,\cdots,y_{n-4},y_{n-3},\cdots,y_n)=\eta(y_{n-3},\cdots ,y_n)
$$
where $y_1,\cdots ,y_n$ denote the euclidean coordinates of
$T_xM$ such that $y_1,\cdots, y_{n-4}$ are in $F$.

Finally, by the second equality in (\ref{eq:3.2.24}), we have
that ${\rm supp}(\eta )\subset F$.  Therefore,
$\eta=\Theta(\mu,x)H^{n-4}\lfloor F$.

\demo{{\rm 3.3.}\ Rectifiability}
We have shown that tangent cones exist at $H^{n-4}$-a.e. $x$ in $S$;
moreover, if (\ref{eq:3.2.24}) holds and $\Theta (\mu , \cdot)$ is 
$H^{n-4}$-approximately continuous at $x \in S$, then any
tangent cones at $x$ are $(n-4)$-subspaces in $T_xM$
(Proposition~\ref{prop:3.2.4}). We adopt the notation of the last section.  
In this section, we will prove
that $S$ is rectifiable, i.e., tangent cones are unique at
$H^{n-4}$-a.e. $x$ in $S$. This in fact follows from 
the work of D. Priess \cite{P}, since 
$\Theta (\nu , \cdot)$ exists almost
everywhere and $\nu$ is Borel regular.
However, for the reader's convenience, we give a direct proof
here by using the structure theorem of Federer (cf.~\cite{Fe}, \cite{Li}).

We may write $S=S_u \cup S_r$, where $S_r$ is a rectifiable set
and $S_u$ is a purely unrectifiable set.
We denote by $G(T_xM, n-4)$ the
Grassmannian of all $(n-4)$-dimensional subspaces in $T_xM$. \enddemo

\specialnumber{3.3.1} \proclaim{Lemma}
\label{lem:3.3.1}
For any $x \in M$ and $V$ in $G (T_xM, n-4)${\rm ,}  
$$
H^{n-4}(P_V(\exp^{-1}_x (B_r(x) \cap S_u))) =0
$$
where $r>0$ is sufficiently small and
$P_V$ denotes the orthogonal projection of $T_xM$ onto $V$
with respect to $g_{x,0}${\rm . }
\endproclaim

This lemma can  easily be proved by modifying the
arguments in the proof of [Fe, 3.3.5] or \cite{Si2}. 
We omit it here.

We want to show that $H^{n-4}(S_u)=0$.  Suppose that it is not
true. Then for $H^{n-4}$-a.e. $x$ in $S_u$, $r>0$ small and
any $V \in G(T_xM,n-4)$,
\begin{equation}
  \label{eq:3.3.1} 
  H^{n-4} (P_V (\exp^{-1}_x (S_u \cap B_r (x)))) =0 \spn{3.3.1}
\end{equation}
and
\begin{equation}
  \label{eq:3.3.2}
  \overline{\lim}_{\lambda \to 0+} \frac{H^{n-4}(S_r \cap B_\lambda (x))}
     {\lambda ^{n-4}} = 0 . \spn{3.3.2}
\end{equation}
Since $H^{n-4} (S_u) >0$, we can choose $x$ in $S_u$ such that
(\ref{eq:3.2.24}), (\ref{eq:3.3.1}) and (\ref{eq:3.3.2}) hold,
and $\Theta (\mu , \cdot)$ is $H^{n-4}$-approximately continuous
at $x$.  As before, we define $\mu_{x, \lambda}$ by
\begin{equation}
  \label{eq:3.3.3}
  \mu_{x, \lambda} (E) = \lambda^{n-4} \mu (\exp_x (\lambda E)) \spn{3.3.3}
\end{equation}
where $E \subset T_xM$.  Let $\left\{ \lambda_k \right\}$ be a
sequence of positive numbers such that $\lim_{k \to\infty} \lambda_k =0$ 
and $\mu_{x, \lambda_k}$ converges weakly
to a tangent measure $\eta$ on $T_xM$. By our choice of $x$ and
the proof of Proposition \ref{prop:3.2.4}, we have that $\eta =
\Theta (\mu , x) H^{n-4} \lfloor V$ for some $(n-4)$-subspace $V$ in
$T_xM$.  We claim:
\begin{equation}
  \label{eq:3.3.4}
  \overline{\lim}_{k \to \infty} 
   \frac{H^{n-4}(P_V (\exp^{-1}_x (S \cap B_{\lambda_k}(x))))}
    {\lambda^{n-4}_k} >0 \, . \spn{3.3.4}
\end{equation}
If this is true, then 
\begin{equation}
  \label{eq:3.3.5}
   \overline{\lim}_{k \to \infty} 
   \frac{H^{n-4}(P_V (\exp^{-1}_x (S_u \cap B_{\lambda_k}(x))))}
    {\lambda^{n-4}_k} >0 \, , \spn{3.3.5}
\end{equation} 
because of (\ref{eq:3.3.2}).  However, this contradicts 
(\ref{eq:3.3.1}).

Now we prove the claimed inequality in (\ref{eq:3.3.4}). As in the 
last section, we may find a sequence of Yang-Mills connections
$A_{i, x, \lambda_k}$ (cf.\ (3.2.3)) such that the $\left| F_{A_{i,x,
      \lambda_k}} \right|^2 \, dV_{x, \lambda_k}$ converge to
$\mu_{x, \lambda_k}$ weakly as $i \to \infty$.  Note that for $k$
large enough, the  $A_{i,x, \lambda_k}$ are well defined in $B_4
(0,g_{x, \lambda_k}) \subset T_xM$. Let us identify $T_xM$ with 
$V \times V^{\perp}$, so that
each point $z \in T_xM$ is of the form $(z',z'')$ with 
$z' \in V$ and $z'' \in V^{\perp}$,
where $V^\perp$ is the orthogonal complement 
of $V$ in $T_xM$. Choose orthonormal
coordinates $z_1, \cdots , z_n$ of $T_xM$ with respect to
$g_{x,0}$, such that   $z_1, \ldots , z_{n-4}$ are coordinates of
$V $ and $z_{n-3}, \ldots , z_n$ are coordinates of
$V^{\perp}$. We usually denote $z'$ by
$(z_1, \ldots , z_{n-4})$ and $z''$ by 
$ (z_{n-3} , \ldots , z_n)$.
We put
$$
  B^2_2 (0) = \left\{ z'' \in V^{\perp} | 
    \left| z'' \right| <2\right\} \, .
$$
Clearly, when $k$ is sufficiently large (so that  $g_{x, \lambda_k}$ is
sufficiently close to the flat metric $g_{x,0}$), we have that
$(z',0) + \{0\}\times B^2_2 (0) \subset B_4 (0, g_{x, \lambda_k})$ 
for any $(z',0)\in V \times \{0\}\cap B_2 (0,g_{x, \lambda_k})$. 

Consider
\begin{equation}
  \label{eq:3.3.6}
  m_{i,k} (z') = \int_{B^2_2 (0)} \left| F_{A_{i,x, \lambda_k}}
    \right|^2 (z' ,z'') \phi^2 (z'') \, dV_k (z'') \, , \spn{3.3.6}
\end{equation}
where $dV_k(z'')$ denotes the induced volume form on $B^2_2 (0)$
by the metric $g_{x, \lambda_k}$, and $\phi \in C^{\infty}_0
(B^2_2 (0))$ with $\int _{B^2_2(0)} \phi^2 dV_{g_{x,0}}=1$.  
Then $m_{i,k}$ is a smooth function of $z'$ in $V \cap
B_2 (0,g_{x, \lambda_k})$.  

For simplicity, we will denote by $D$ the covariant derivative
associated to each $A_{i,x, \lambda_k}$ unless further
specification is needed.  For simplicity, we often abbreviate
$\frac{\partial}{\partial z_{\alpha}}$ as $\partial_{\alpha}$.
One computes
  \begin{eqnarray}
&&\spn{3.3.7}\\
     \frac{\partial m_{i,k}(z')}{\partial z_{\alpha}}  
    &=&  \frac{\partial}{\partial z_{\alpha}}
         \left( \int_{B^2_2 (0)} \left| F_{A_{i,x, \lambda_k}}
         \right|^2 (z', z'') \phi^2 (z'') \, dV_k (z'')
       \right)\nonumber \\[2ex]
 &=& 2 \int_{B^2_2 (0)} \frac{\partial}{\partial z_\alpha }
\left ( F_{A_{i,x, \lambda_k}} \, , \, 
F_{A_{i,x, \lambda_k}} \right )
   (z',z'') \phi^2 (z'') \, dV_k (z'') \, .\nonumber
    \label{eq:3.3.7}
  \end{eqnarray}
Since $g_{x, \lambda_k}$ converges to the flat metric $g_{x,0}$ on
$T_xM$ as $k \to \infty$,  
\begin{eqnarray}
    g_{x, \lambda_k} \left(\partial_{\alpha} ,
    \partial_{\beta} \right) 
&=& \delta_{\alpha \beta} + o(1) , ~~\alpha, \beta = 1,2,
  \cdots , n \, , \spn{3.3.8}
  \label{eq:3.3.8}
\\
 \nabla^k_{\partial_{\alpha}}
  \partial_{\beta} &=& o (1) \, , \, 
  \alpha, \beta = 1,2, \cdots , n \, ,  \spn{3.3.9} \label{eq:3.3.9}
\end{eqnarray} 
where $\nabla^k$ denotes the Levi-Civita connection of
$g_{x, \lambda_k}$, and $o(1)$ always denotes a quantity which
converges to zero as $k \to \infty$.  It follows from
(\ref{eq:3.3.8}) and (\ref{eq:3.3.9}) that
\begin{eqnarray*} 
 \frac{\partial}{\partial z_\alpha }
&=& 2 \sum^n_{\beta, \gamma,\beta', \gamma'=1} \left ( \left | F_{A_{i,x, \lambda_k}} \right |^2  
     D_{\partial_{\alpha}}
      F_{A_{i,x, \lambda_k}}  
    \left(\partial_{\beta}  \, , \, 
      \partial_{\gamma} \right) \, ,  
    F_{A_{i,x, \lambda_k}}  
    \left( \partial_{\beta'} \, , \, 
     \partial_{\gamma'} \right) \right ) \\[2ex]
&&\cdot \ 
g_{x, \lambda_k}^{\beta\beta'} g_{x, \lambda_k}^{\gamma\gamma'}
    + o(1) \left |  F_{A_{i,x, \lambda_k}} \right |^2\nonumber ,
  \end{eqnarray*}  
where $\{g_{x, \lambda_k}^{\alpha\beta}\}$ is the inverse matrix
of $\{g_{x, \lambda_k}(\partial _\alpha, \partial_\beta)\}$.
By the second Bianchi identity $DF_{A_{i,x, \lambda_k}} =0$, we
deduce from the above
\begin{eqnarray}
 && \spn{3.3.10} \label{eq:3.3.10}\\
  \frac{\partial}{\partial z_\alpha } \left | F_{A_{i,x,
      \lambda_k}}  \right |^2 
      &=& 4 \sum^n_{\beta, \gamma,\beta', \gamma' =1}  
        \left(D_{\partial_{\beta}} 
        F_{A_{i,x, \lambda_k}}
          \left( \partial_{\alpha} \, , \, 
            \partial_{\gamma}\right) \, , \, 
            F_{A_{i,x, \lambda_k}} 
              \left( \partial_{\beta'} \, , \, 
             \partial_{\gamma'}\right)  \right)\nonumber \\[2ex]
&&\cdot \ 
        g_{x, \lambda_k}^{\beta\beta'} g_{x, \lambda_k}^{\gamma\gamma'}
           + o(1) \left|  F_{A_{i,x, \lambda_k}} \right|^2\nonumber\\[2ex]
      &=& 4 \sum \frac{\partial}{\partial z_\beta } \left (  F_{A_{i,x, \lambda_k}} 
          \left(\partial_{\alpha}  \, , \, 
           \partial_{\gamma}\right ) \, , \, 
           F_{A_{i,x, \lambda_k}} 
           \left(\partial_{\beta'} \, , \, 
            \partial_{\gamma'} \right)  \right)
      g_{x, \lambda_k}^{\beta\beta'} g_{x, \lambda_k}^{\gamma\gamma'}\nonumber\\[2ex]
      &&  -\ 4 \sum^2_{\beta , \gamma =1} \left ( 
             F_{A_{i,x, \lambda_k}} 
            \left( \partial_{\alpha} \, , \, 
           \partial_{\gamma} \right) \, , \,
          D_{\partial_{\beta}}
             F_{A_{i,x, \lambda_k}} 
           \left(\partial_{\beta'}  \, , \, 
             \partial_{\gamma'}\right)  \right) \nonumber\\ [2ex]
&&\cdot \ 
     g_{x, \lambda_k}^{\beta\beta'} g_{x, \lambda_k}^{\gamma\gamma'}
           + o(1) \left| F_{A_{i,x, \lambda_k}}\right|^2\nonumber \, .
\end{eqnarray} 
Since $A_{i,x, \lambda_k}$ is a Yang-Mills connection with
respect to $g_{x, \lambda_k}$,  
\begin{equation}
\label{eq: 3.3.11}
g_{x, \lambda_k}^{\beta\beta'} D_{\partial_{\beta}}  F_{A_{i,x,
 \lambda_k}} \left( \partial_{\beta'} \, , \, 
             {\partial_{\gamma}}\right) = 0. \spn{3.3.11}
\end{equation}
Combining this with (\ref{eq:3.3.10}), we deduce for $\alpha
\leq n-4$,
\begin{eqnarray*}
 \frac{\partial m_{i,k}(z')}{\partial z_{\alpha}} 
 &=&   4 \sum^{n}_{\beta, \gamma =1} \int_{B^2_2 (0)}
       \partial_\beta \left( F_{A_{i,x, \lambda_k}}
        \left( \partial_{\alpha} \, , \, 
             \partial_{\gamma}\right) \, , \, 
           F_{A_{i,x, \lambda_k}} 
           \left( \partial_{\beta'} \, , \, 
             \partial_{\gamma'}\right)  \right) \\ [2ex]
    &&\cdot \  g_{x, \lambda_k}^{\beta\beta'} g_{x, \lambda_k}^{\gamma\gamma'}
           \phi^2 (z'') dV_k(z'')    +  
  o(1) \int_{B^2_2 (0)} \left|  F_{A_{i,x, \lambda_k}}
  \right|^2 \phi^2 (z'') d V_k (z'') \\[2ex]
 &=&  -4 \sum_{\beta=n-3}^{n} \int_{B^2_2 (0)} 
  \left(\partial_{\alpha} 
  \rfloor  F_{A_{i,x, \lambda_k}} \, , \,
  \partial_{\beta'} 
  \rfloor F_{A_{i,x, \lambda_k}}  \right) g_{x, \lambda_k}^{\beta\beta'} 
  \partial_{\beta}\phi^2 (z'') dV_k(z'') \\[2ex]
  &&  +\ 
  4 \sum^{n-4}_{\beta=1} \frac{\partial}{\partial z_{\beta}} 
  \left( \int_{B^2_2 (0)} 
   \left( \partial_{\alpha}\rfloor  F_{A_{i,x, \lambda_k}} \, , \, 
   \partial _\beta \rfloor
    F_{A_{i,x, \lambda_k}}  \right) \phi^2 (z'')
    dV_k (z'') \right) \\[2ex]
  && + \
  o(1) \int_{B^2_2 (0)} \left|  F_{A_{i,x, \lambda_k}}
    \right|^2 \phi^2 (z'') dV_k (z'') \, .
\end{eqnarray*} 
To estimate these derivatives, we need the following: 

\specialnumber{3.3.2} \proclaim{Lemma}
\label{lem: 3.3.2}
Let $\{A_{i,x,\lambda,_k}\}${\rm ,} $x$ etc{\rm .}\ be defined as above{\rm .} Then for any
$\alpha \leq n-4${\rm ,}
\begin{equation}
        \label{eq:3.3.12}
        \lim_{k \to \infty} \lim_{i \to \infty} \int_{B_4
          (0,g_{x, \lambda_k})} 
        \left| \frac{\partial}{\partial z_{\alpha}} \rfloor
          F_{A_{i,x, \lambda_k}} \right|^2
        dV_{x, \lambda_k} =0 \, . \spn{3.3.12}
\end{equation}
\endproclaim

\demo{Proof}
By our assumption, $\left | F_{A_{i,x, \lambda_k}}\right |^2 dV_k$ 
converges to $\mu_{x, \lambda_k}$ weakly
as $i \to \infty$ and $\mu_{x, \lambda_k} \to \eta$ as $k\to \infty$.  
Moreover, $\eta$ is of the form $\Theta (\mu ,
x) H^{n-4} \lfloor V$. Therefore, for any $\delta > 0$, 
\begin{equation}
\label{eq:3.3.13}
\lim_{k \to \infty} \lim_{i \to \infty} \int_{B_4
(0,g_{x, \lambda_k})\backslash T_\delta (V)} 
        \left|F_{A_{i,x, \lambda_k}} \right|^2
        dV_{x, \lambda_k} =0, \spn{3.3.13}
\end{equation}
where $T_\delta (V)$ denotes the $\delta$-tubular neighborhood
of the subspace $V$.

Let $x_1^k, \cdots, x_{n-4}^k \in B_{\lambda_k}(x)\cap S$ be chosen
as in (3.2.25) and (3.2.26). We let $V^k$   be the subspace
in $T_xM$ spanned by $\xi_1^k=\exp _x^{-1}(x_1^k),\ldots,\xi_{n-4}^k=\exp _x^{-1}(x_{n-4}^k)$. Then the  $V^k$
converge to $V$, and these $V^k$ are spanned by $\xi_j = \lim_{k\to \infty } \xi_j^k$. Moreover, we may assume that
$d_{g_{x,0}}(\xi_i, \xi_j)
\ge s$ for $i\not= j$ and $d_{g_{x,0}}(\xi_i, 0) \ge s$, where $s$
is  as   given in Lemma 3.2.3. We have shown in the proof of Lemma 3.2.3,
\begin{equation}
\label{eq:3.3.14}
6^{4-n} \mu_{x,\lambda_k} 
(B_6(\xi _j^k, g_{x,\lambda_k})) \ge \Theta (\mu, x)- \varepsilon(\lambda_k). \spn{3.3.14}
\end{equation} 
Note that $\varepsilon(\cdot )$ is a nondecreasing function with
$\lim _{r\to 0}\varepsilon (r)=0$.
Choose $i(k)$ such that for any $i \ge i(k)$,
\begin{equation}
\label{eq:3.3.15}
6^{4-n} \int _{B_6(\xi _j^k, g_{x,\lambda_k})} |F_{A_{i,x,\lambda_k}}|^2 dV_k
 \ge \Theta (\mu, x)-  2 \varepsilon(\lambda_k). \hskip.5in\spn{3.3.15}
\end{equation} 
Since $\lim\limits_{k\to \infty} \mu_{x,\lambda_k} = \eta$
and $\Theta (\eta, \xi_j) = \Theta (\mu, x)$ (cf.\ (3.2.28)),
by increasing $\varepsilon(r)$ if necessary,
we may assume that 
\begin{equation}
\label{eq:3.3.16}
\left |6^{4-n} \mu_{x,\lambda_k} (B_6(\xi _j^k, g_{x,\lambda_k})) - 
\Theta (\mu, x)\right | \le \varepsilon(\lambda_k). \spn{3.3.16}
\end{equation}
By taking $i(k)$ big enough, we may further have that for $i \ge i(k)$,
\begin{equation}
\label{eq:3.3.17}
6^{4-n} \int _{B_6(\xi_j^k, g_{x,\lambda_k})} |F_{A_{i,x,\lambda_k}}|^2 dV_k
\le \Theta (\mu, x) +  2 \varepsilon(\lambda_k).\hskip.5in\quad \spn{3.3.17}
\end{equation} 
Then we deduce from this and the monotonicity (Theorem 2.1.2) that

\begin{equation}
\label{eq:3.3.18}
\int _{B_6(\xi_j^k, g_{x,\lambda_k})\backslash B_s(\xi_j^k, g_{x,\lambda_k})} 
(\rho _j^k)^{4-n} |\frac{\partial }{\partial \rho _j^k} \rfloor 
F_{A_{i,x,\lambda_k}}|^2 dV_k
\le 2 \varepsilon(\lambda_k), \hskip.5in\spn{3.3.18}
\end{equation}
where $\xi_0 =0$,
$\rho^k_j$ is the distance from $\xi_j^k$ of $g_{x,\lambda_k}$
($j=0,1,\cdots, n-4$).

Then the lemma follows from (3.3.18) and (3.3.13) and
the fact that the  $\xi_j$ span the subspace $V$.
\enddemo

Notice that the integral 
$$
\int_{B_4(0,g_{x,
      \lambda_k})}\left| F_{A_{i,x, \lambda_k}} \right|^2 dV_{x,
    \lambda_k}
$$
is uniformly bounded. Thus we have
      \begin{equation}
        \label{eq:3.3.19}
        \grad m_{i,k} = f_{i,k} + \divv (u_{i,k}) \, , \spn{3.3.19}
      \end{equation}
where $f_{i,k}:V \cap B_2 (0,g_{x, \lambda_k}) \to V$ and
$u_{i,k} : V \cap B_2 (0,g_{x, \lambda_k}) \to V \times V$ are
functions, such that
\begin{equation}
  \label{eq:3.3.20}
   \lim_{k \to \infty} \lim_{i \to \infty} 
   \int_{V \cap B_2 (0,g_{x, \lambda_k})} 
   \left( \left| f_{i,k} \right| + \left| u_{i,k} \right| \right)
     dV_{x,0} = 0 \, . \qquad\hskip.5in\spn{3.3.20}
\end{equation}
Then it follows (cf.\ \cite{AL}) that
there are constants $C_{i,k}$, such that
\begin{equation}
  \label{eq:3.3.21}
  \lim_{k \to \infty} \lim_{i \to \infty} 
  \| m_{i,k} - C_{i,k} \|_{L^1 (V \cap B_\frac{4}{3}(0,g_{x, \lambda_k}))} =0. \spn{3.3.21}
\end{equation}
In fact, since 
$$
\lim_{k \to \infty} \lim_{i \to \infty} \left|
  F_{A_{i,x, \lambda_k}} \right|^2 \, dV_{x, \lambda_k} = \eta
$$
and 
$$
\eta = \Theta (\mu , x) H^{n-4} \lfloor V,
$$
we have 
$$
\lim_{k \to \infty} \lim_{i \to \infty} C_{i,k} = \Theta (\mu ,
x) >0.
$$

For $k$ sufficiently large, the ball $B_{\frac{3}{2}}(0,g_{x,0})$
is contained in every $B_{\frac{4}{3}}(0,g_{x, \lambda_k})$. Then
for any $\xi \in C^{\infty}_0 (V \cap
B_{\frac{3}{2}}(0,g_{x,0}))$,
\begin{eqnarray}
  \label{eq:3.3.22}
&&\spn{3.3.22} \\
  && \Theta (\mu , x) \int_{V \cap B_{\frac{3}{2}}(0,g_{x,0})}
     \xi (z') \, dz'  \nonumber
\\[2ex]  &&\hskip.75in = ~ \lim_{k \to \infty} \lim_{i \to \infty} 
         \int_{V \cap B_{\frac{3}{2}(0,g_{x,0})}}
         \xi (z') m_{i,k} (z') \, dz'\nonumber  \\[2ex]
     &&\hskip.75in = ~ \lim_{k \to \infty} \lim_{i \to \infty} 
         \int_{B_2(0,g_{x, \lambda_k})} 
         \left| F_{A_{i,x, \lambda_k}} \right|^2 (z', z'')
    \xi (z') \phi^2 (z'') \, dV_{x, \lambda_k}
     \nonumber \\[2ex]
      &&\hskip.75in = ~\lim_{k \to \infty}  \int_{B_2(0,g_{x, \lambda_k})}
     \xi (z') \phi^2 (z'') \, d \mu_{x, \lambda_k} (z', z'')\nonumber \, .
\end{eqnarray}
However, as a weak limit of Radon measures $\left|  F_{A_i} \right|^2
\, dV_g$, the measure $\mu$ is of the form $\left|  F_A \right|^2
\, dV_g + \nu$.  After scaling, we have
\begin{equation}
  \label{eq:3.3.23}
  \mu_{x, \lambda_k} = \left| F_{A_{x, \lambda_k}} \right|^2 \, 
  dV_{x, \lambda_k} + \nu_{x, \lambda_k} \, , \spn{3.3.23}
\end{equation}
where $A_{x, \lambda_k}$ is a connection on $T_xM \backslash
\lambda^{-1}_k \exp^{-1}_x (S)$  as defined in (3.2.3),
and $\nu_{x, \lambda_k}$ is a Radon measure on $T_xM$ of the form
\begin{equation}
  \label{eq:3.3.24}
  \Theta (\mu_{x, \lambda_k}, \cdot) H^{n-4} \lfloor \lambda^{-1}_k
  \exp^{-1}_x (S) \, . \spn{3.3.24}
\end{equation}
Using the second equation in (\ref{eq:3.2.24}) which holds at $x$
by our assumption, we see that 
\begin{equation}
  \label{eq:3.3.25}
  \lim_{k \to \infty} \int_{B_2 (0,g_{x, \lambda_k})}
  \xi (z') \phi^2 (z'') \left| F_{A_{x, \lambda_k}} \right|^2
  dV_{x, \lambda_k} =0 \, . \hskip.5in \spn{3.3.25}
\end{equation}
Hence, by (\ref{eq:3.3.25})
\begin{eqnarray*}
 &&\hskip-25pt \Theta (\mu ,x) \int_{V \cap B_{\frac{3}{2}}(0,g_{x,0})}
   \xi (z') \, dz' \\[2ex]
&&\quad=\ \lim_{k \to \infty} \int_{B_{\frac{3}{2}}(0,g_{x,0}) \cap
  \lambda^{-1}_k \exp^{-1}_x (S)} \xi(z') \Theta (\mu_{x,
  \lambda_k}, (z',z'') )\, dH^{n-4} (z',z'')\nonumber \\[2ex]
&&\quad=\  \Theta (\mu ,x) \lim_{k \to \infty}
     \int_{B_{\frac{3}{2}}(0,g_{x,0}) \cap
  \lambda^{-1}_k \exp^{-1}_x (S)}  \xi (z') \, dH^{n-4} (z',z'')\nonumber \, .
\end{eqnarray*} 
Since $\Theta (\mu ,x) >0$, this implies
\begin{eqnarray*}
  &&\hskip-.5in\overline{\lim_{k \to \infty}}
  \frac{H^{n-4}(P_V (\exp^{-1}_x(S \cap B_{\lambda_k} (x))))}
  {\lambda^{n-4}_k} \\ [2ex]
&&\quad=\   \overline{\lim_{k \to \infty}}
  H^{n-4}(P_V (\lambda^{-1}_k \exp^{-1}_x
  (S \cap B_1 (0,g_{x, \lambda_k})) \\[2ex]
&&\quad\geq\  \mbox{~Vol~} (V \cap B_{\frac{1}{2}}(0,g_{x,0}))  > 0.
\end{eqnarray*} 
Thus (\ref{eq:3.3.4}) is proved and we obtain a contradiction to
(\ref{eq:3.3.1}). Hence,\break $H^{n-4}(S_u) =0$ and we have shown
the following:
\specialnumber{3.3.3} \proclaim{Proposition}
  \label{prop:3.3.2}
Let $(S_b, \Theta)$ be the blow\/{\rm -}\/up locus of a weakly convergent
sequence $\left\{ A_i \right\}${\rm .}  Then its support $S_b$ is $H^{n-4}$\/{\rm -}\/rectifiable{\rm .}  In
particular{\rm ,} for $H^{n-4}${\rm -a.e.} $x$ in $S_b${\rm ,} there is a unique
tangent subspace $T_xS_b \subset T_xM${\rm .}
\endproclaim

\section{Structure of blow-up loci}
\label{chap:4}

In this chapter, we study the  geometry of blow-up loci. 

\demo{{\rm 4.1.}  Bubbling Yang-Mills connections}
We assume that $\{A_i\}$ converges to 
an admissible Yang-Mills connection
$A$ with the blow-up locus $(S, \Theta)$\break (cf.\ Lemma 3.1.4).
It is shown in Section 3.3 that $S$ is $H^{n-4}$-rectifiable.
We will adopt the notation of the  last chapter.

If $n=4$, $S$ consists of finitely many points. 
K. Uhlenbeck further showed that when $i$ is sufficiently large, $A_i$ approaches
  a connected sum of $A$ with certain Yang-Mills connections
on the unit sphere $S^4$. These later connections are called
bubbling connections. 

In this section, we analyze the structure of $A_i$ near $S$
when $i$ is sufficiently large. We will construct bubbling
connections on $\RR^n$ as $A_i$ approaches   $A$. 

Recall that $\mu$ is the weak limit of Radon measures
$|F_{A_i}|^2 dV_g$ and is of the form $|F_A|^2 dV_g + 
\Theta (\mu, \cdot) H^{n-4}\lfloor S$. \enddemo

\specialnumber{4.1.1} \proclaim{Proposition}
\label{prop:4.1.1}
Let $x\in S$ satisfy\/{\rm :}\/
\begin{itemize}
\item[{\rm (1)}] 
The tangent plane $V=T_xS\subset T_xM$ exists uniquely\/{\rm ;}\/
\item[{\rm (2)}] 
{\rm (\ref{eq:3.2.24})} holds for $\mu$ and $A${\rm .}
\end{itemize}
  
Then there are linear transformations $\sigma_i: T_xM\mapsto T_xM$
such that a subsequence of $\sigma_i^*\exp_x^* A_i$ converges to
a Yang\/{\rm -}\/Mills connection $B$ on $T_xM$ such that 
$F_B\not= 0$ and $v\rfloor F_B\equiv 0$ for any $v \in V${\rm .}  
Such a connection $B$ is called a bubbling connection at $x\in S${\rm .}
\endproclaim

The rest of this section is devoted to the  proof of Proposition 4.1.1.

Let $A_{i,x,\lambda}$ be the scaled connections on $T_xM$ defined 
in (3.2.4), i.e.,  
\begin{equation}
\label{eq:4.1.1}
A_{i,x,\lambda}= \tau _\lambda ^\ast \exp^\ast_x A_i, \spn{4.1.1}
\end{equation}
where $\tau_\lambda (v) = \lambda v $ for any $v$ in $T_xM$.
Each $A_{i,x,\lambda}$ is a Yang-Mills connection with respect to
the scaled metric $g_{x,\lambda}$. As $i$ tends to  infinity,
$|F_{A_{i,x,\lambda}}|^2dV_{x,\lambda}$ converges to
$\mu_{x,\lambda}$ weakly. On the other hand, as $\lambda$ tends
to zero, $\mu_{x,\lambda}$ converges to
$\Theta(\mu,x)H^{n-4}\lfloor V$ weakly.  
Therefore, there is a sequence $\lambda_{i}$ such that the Radon measure
$|F_{A_{i,x,\lambda_{i}}}|^2dV_{x,\lambda_i}$ converges to
$\Theta(\mu,x)H^{n-4}\lfloor V$ weakly. Moreover,  
modulo gauge transformations,
$A_{i,x,\lambda_i}$ converges to $0$ uniformly
on any compact subsets in $T_xM\backslash V$. This implies particularly that
for $i$ sufficiently large,
\begin{equation}
\label{eq:4.1.2}
|F_{A_{i,x,\lambda_{i}}}|(v)\leq\frac{\varepsilon(r)}{r^2}.\spn{4.1.2}
\end{equation}
We also have (cf.\ Lemma 3.3.2)
\begin{equation}
\label{eq:4.1.3}
\lim_{i\to\infty}\left(\sum^{n-4}_{\alpha=1}\int_{B_2(0,g_{x,0})}|
       \frac{\partial}{\partial
    z_{\alpha}}\rfloor F_{A_{i,x,\lambda_{i}}}|^2
dV_{g_{x,\lambda_{i}}}\right)=0, \spn{4.1.3}
\end{equation}
where $\{z_1,\cdots,z_{n-4}\}$ is an
orthogonal coordinate system of $V$.

As in the last section, we denote by $z=(z^\prime,
z^{\prime\prime})$ a point in $T_xM$ with $z^\prime\in V,
z^{\prime\prime}\in V^{\bot}$. We will identify $V$ and $V^\bot$ with the subspaces
$V\times\{0\}$ and $\{0\}\times V^\bot$ in $T_xM$.
\specialnumber{4.1.2} \proclaim{Lemma}
\label{lem:4.1.2}
There are points $z^{\prime}_{i}$ in $V\cap
B_{\frac{1}{2}}(0,g_{x,0})$ with $\lim_{i\to \infty} z_i^\prime = 0${\rm ,}
such that 
%
\begin{equation}
\label{eq:4.1.4}
{\ninepoint \lim_{i\to\infty}\left( \sup_{0<r\leq\frac{1}{2}}
r^{4-n}\int_{V\cap B_r(z^{\prime}_{i},g_{x,0})}dx^{\prime}
\int_{V^\bot\cap B_{\frac{1}{2}}(0,g_{x,0})}
\sum^{n-4}_{\alpha=1}|\frac{\partial}{\partial
  z_{\alpha}}
\rfloor F_{A_{i,y,\lambda_{i}}}|^2dV_{x,\lambda_{i}}
\right)=0. } \spn{4.1.4}
\end{equation} 
\endproclaim    

\demo{Proof}
We prove this by contradiction. Suppose that the lemma is false,
then we can find $\delta >0$ and $s \in (0,\frac{1}{2})$, 
such that for any $i$ and
$z^\prime\in V\cap B_s(0,g_{x,0})$, there is at least
one $r=r(i,z^\prime)$ such that
\begin{equation}
\label{eq:4.1.5}
r^{4-n}\int_{V\cap
  B_r(z^\prime,g_{x,0})}dx^\prime\int_{V^{\bot}\cap B_
\frac{1}{2}(0,g_{x,0})}\sum^{n-4}_{\alpha=1}|\frac{\partial}{\partial
z_{\alpha}}\rfloor F_{A_{i,x,\lambda_{i}}}|^2 d x''\geq\delta.\hskip.4in \spn{4.1.5}
\end{equation} 
By (\ref{eq:4.1.3}), $\lim_{i\to\infty}r(i,z^\prime)=0$ for any
$z^\prime$.  For each $i$, we cover $V\cap
B_{\frac{1}{2}}(0,g_{x,0})$ by finitely many disjoint balls
$V\cap B_{r(i,z^{\prime}_{i\alpha})}(z^\prime_{i\alpha},g_{x,0})$
($\alpha=1,2,\cdots,m_{i}$), such that 
\begin{equation}
\label{eq:4.1.6}
V\cap B_s(0,g_{x,0})\subset\bigcup^{m_{i}}_{\alpha=1}V\cap 
B_{2r(i,z^\prime_{i\alpha})}(z^\prime_{i\alpha,},g_{x,0}). \spn{4.1.6}
\end{equation}
Then,
\begin{eqnarray}
  \delta \left(\frac{s}{2}\right)^{n-4}
   &\leq& \delta 2^{4-n}\sum^{m_{i}}_{\alpha=1}(2r(i,z^\prime_{i\alpha}))^{n-4}
         =\delta \sum^{m_i}_{\alpha=1}r(i,z^{\prime}_{i\alpha})^{n-4}
         \nonumber\\[2ex]
  &\leq&\sum^{m_{i}}_{\alpha=1}\int_{V\cap
    B_{r(i,z^\prime_{i\alpha})}(z^\prime_{i\alpha},g_{x,0})}
    dx^{\prime}\int_{V^{\bot}\cap
    B_{\frac{1}{2}}(0,g_{x,0})}\sum^{n-4}_{\beta=1}
    \left|\frac{\partial}{\partial z_{\beta}}\rfloor
    F_{A_{i,x,\lambda_{i}}} \right|^2 dx''
    \nonumber\\[2ex]
&\leq& \int_{B_2(0,g_{x,0})}\left(\sum^{n-4}_{\beta=1}
    \left| \frac{\partial}{\partial z_\beta}\rfloor
     F_{A_{i,x,\lambda_{i}}} \right|^2\right) dV_{x,\lambda_{i}}.
   \nonumber
\end{eqnarray}
This is impossible when $i$ is sufficiently large because
of (\ref{eq:4.1.3}). 
\enddemo

Observe that for any $\delta>0$,
\begin{equation}
\label{eq:4.1.7}
\max_{z^{\prime\prime}\in V^{\bot}\cap
  B_{\frac{1}{2}}(0,g_{x,0})}\delta^{4-n}
\int_{B_{\delta}(z_i' + z'',g_{x,0})}
|F_{A_{i,x,\lambda_{i}}}|^2dV_{x,\lambda_{i}}\geq\varepsilon,\hskip.5in \spn{4.1.7}
\end{equation}
where $\varepsilon$ is as in Theorem \ref{th:2.2.1}.
Otherwise, $A_{i,x,\lambda_{i}}$ would converge to a smooth Yang-Mills
connection on $(V\cap B_\delta(z^{\prime}_{i},g_{x,0}))\times 
(V^{\bot}\cap B_{\frac{1}{2}}(0,g_{x,0}))$,
contradicting our assumption  on
$A_{i,x,\lambda}$.

Because of (\ref{eq:4.1.7}), we can find
$\delta_i\in(0,\frac{1}{2})$ and $z^{\prime\prime}_{i}\in
V^{\bot}\cap B_{\frac{1}{4}}(0,g_{x,0})$, such that
\begin{eqnarray}
\label{eq:4.1.8}
&& \quad\delta^{4-n}_{i}\int_{B_{{\delta}_{i}}(z^{\prime}_{i}+z^{\prime\prime}
  _{i},g_{x,0})}|F_{A_{i,x,\lambda_{i}}}|^2dV_{x,\lambda_{i}} \spn{4.1.8}\\[2ex]
&&\qquad \quad\qquad =\ \max_{z^{\prime\prime}\in V^{\bot}\cap
  B_{\frac{1}{2}}(0,g_{x,0})}\delta^{4-n}_{i}\int_{B_{{\delta}_{i}}
  (z_i' + z'',g_{x,0})}|F_{A_{i,x,\lambda_{i}}}|^2
  dV_{x,\lambda_{i}}
 = \frac{\varepsilon}{4}.\nonumber
\end{eqnarray} 
One may even take $z^{\prime\prime}_{i}$ such that
$\lim_{i\to\infty}z^{\prime\prime}_{i}=0$.
Now we define new connections
\begin{equation}
\label{eq:4.1.9}
B_i(y)=A_{i,x,\lambda_{i}}(z^{\prime}_{i}+z^{\prime\prime}_{i}+\delta_{i}y). \spn{4.1.9}
\end{equation}
Each $B_i$ is a Yang-Mills connection with respect to the scaled
metric $g^{\prime}_{i}=\delta^{-2}_{i}g_{x,\lambda_{i}}$ on
$B_{4R_{i}}(0,g_{x,0})$, where $R_i=(4\delta_{i})^{-1}$.  Note
that the based manifolds $(T_xM, g^{\prime}_{i}, z_i'+ z_i'')$ converge to 
$(T_xM, g_{x,0},0)$ as
$i\to\infty$.

Using (\ref{eq:4.1.4}) and (\ref{eq:4.1.8}), we have
\begin{equation}
\label{eq:4.1.10}
\lim_{i\to\infty}\left(\sum^{n-4}_{\alpha=1}\int_{B_{R_{i}}(0,g_{x,0})}|
    \frac{\partial}{\partial z_{\alpha}}\rfloor 
    F_{B_{i}}|^2dV_{g^{\prime}_{i}} \right)=0, \spn{4.1.10}
\end{equation}
\begin{equation}
\int^{\phantom{\int}}_{B_{1}(0,g_{x,0})}|F_{B_{i}}|^2dV_{g^{\prime}_{i}}\spn{4.1.11}
\label{eq:4.1.11}
= 
\ \max_{y\in V^{\bot}\cap
B_{R_{i}-1}(0,g_{x,0})}\int_{B_{1}((0,y),g_{x,0})}
|F_{B_{i}}|^2dV_{g^{\prime}_{i}}= \frac{\varepsilon}{4}.  
\end{equation} 
It follows from the monotonicity formula that
\begin{equation}
\label{eq:4.1.12}
\sup_{i}\left\{\int_{B_R(0,g_{x,0})}|F_{{B}_{i}}|^2dV_{g^{\prime}_{i}}\right\}
\leq C(\Lambda)R^{n-4}, \spn{4.1.12}
\end{equation}
for $0<R<R_i$, where $C(\Lambda)$ denotes a constant depending only on $\Lambda$.

By (\ref{eq:4.1.12}), Proposition 3.1.2 and 
by taking a subsequence if necessary, 
we may assume that $B_i$ converges to an admissible Yang-Mills connection\break $B$. 
It follows from (\ref{eq:4.1.11})
that $B$ is a smooth Yang-Mills connection on\break $(V\cap
B_1(0,g_{x,0}))\times V^{\bot}\subset T_xM$ with respect to $g_{x,0}$.

Moreover, (\ref{eq:4.1.10}) implies that for any $v\in V$,
\begin{equation}
\label{eq:4.1.13}
v\rfloor F_{B}=0, \spn{4.1.13}
\end{equation}
whenever $B$ is well-defined. 

On the band $\left(\left(V\cap B_1(0,g_{x,0})\right)\right)\times
V^{\bot}$, we write

$$
B= \sum _{\alpha =1}^n B^\alpha dy_\alpha,
$$
where $B^\alpha \in {\rm Lie }(G)$ and $y_1, \cdots, y_n$ are
euclidean coordinates such that\break $y_1,\cdots, y_{n-4}$ are tangent to
$V$ along $V$. Let us eliminate $B^\alpha$ for $\alpha \le n-4$
inductively. First, by a gauge transformation, we may assume that
$B^1 =0$; then (\ref{eq:4.1.13}) implies that all $B^\alpha$ are independent
of $y_1$. Again taking a gauge transformation, we can get rid of $B^2$,
and so on. Eventually, by finitely many gauge transformations, we arrive at a
connection, still denoted by $B$, which is a pull-back of some connection
on $V^{\bot}$. This implies that $B$ extends
to   a smooth connection on $T_xM$. Proposition 4.1.1 is proved.

\demo{{\rm 4.2.}\ Blow-up loci of anti-self-dual instantons}
Now we assume that $\{A_i\}$ is a sequence of $\Omega$-anti-self-dual instantons
which converge to an admissible $\Omega$-anti-self-dual instanton $A$,
where $\Omega$ is a form on $M$ of degree $n-4$. The closedness of 
$\Omega $ is not needed in this section.
Let $S\subset M$ be the blow-up locus of $\{A_i\}$. Here, 
we will show that $\Omega$ restricts to the induced volume form on~$S$.
If $\Omega$ is a calibrating form as in [HL], then $S$ is
calibrated by $\Omega$ and is  particularly minimal. 

First we observe that there is more information on
the bubbling connection constructed in
Proposition \ref{prop:4.1.1} in case of anti-self-dual instantons. \enddemo

\specialnumber{4.2.1} \proclaim{Proposition}
\label{prop:4.2.1}
Let $M${\rm ,} $g${\rm ,} $\Omega$, $\{A_i\}${\rm ,} $A$ and $S$ be as above{\rm .} Suppose that
  $x\in S$ satisfies\/{\rm :}\/
\begin{itemize}
\item[{\rm (1)}] The tangent cone $T_xS\subset T_xM$ exists uniquely\/{\rm ;}\/
\item[{\rm (2)}] {\rm (\ref{eq:3.2.24})} holds for $\mu$ and $A${\rm ,} where $\mu$
  is the weak limit of Radon measures $|F_{A_{i}}|^2 dV_g${\rm . }
\end{itemize} 

Then there is an $\Omega_x$\/{\rm -}\/anti\/{\rm -}\/self\/{\rm -}\/dual instanton $B$ on $T_xM${\rm ,}
where $\Omega _x = \Omega |_{T_xM}${\rm ,}
such that $F_B \not= 0${\rm ,} $\tr(F_B) =0$
and $v\rfloor F_B=0$ for any $v\in T_xS${\rm .}
\endproclaim

\demo{Proof} The proof is basically the same as the 
proof of Proposition \ref{prop:4.1.1}. 

First, we observe that $\tr(F_{B_i})$ converges to zero uniformly
as $i$ tends to infinity, where $B_i$ are the scaled connections
defined in (4.1.9). This is because $\tr (F_{A_i})$ are harmonic $2$-forms
with uniformly bounded $L^2$-norm.

Secondly, we observe that $B_i$ are $\Omega_i'$-anti-self-dual with respect to
the metric $g'_i$ and the closed form $\Omega_i'$ of degree $n-4$
on $B_{4R_i}(0, g_{x,0})$ defined by
$$
\Omega_i' = \tau_{(z_i', z_i'')}^{\delta_i*}\exp_x^*\Omega ,
$$
where $\tau_{(z_i', z_i'')}^{\delta_i}: T_xM\mapsto T_xM,~y \mapsto (z_i', z_i'') + \delta_i y$.

Since $(z_i', z_i'')$ goes to $0$ as $i$ tends to $\infty$, 
$\Omega_i'$ converges to $\Omega_x$. Therefore, the limit connection
$B$ is $\Omega_x$-anti-self-dual with respect to $g_{x,0}$ and $\tr(F_B) =0$.

The rest of the  proof follows the same arguments as
those in the proof of Proposition \ref{prop:4.1.1}.
\enddemo

\specialnumber{4.2.2} \proclaim{{C}orollary}
\label{cor:4.2.2}
Let $x\in S$ be as in the last proposition\/{\rm ;} then $\Omega_x$ restricts 
to a volume form on $T_xS \subset T_xM$ which is induced by the
flat metric $g_{x,0}${\rm .}
\endproclaim
\demo{Proof}
We identify $T_xM$ with $\RR^n$, where $n$ is the dimension of $M$, 
such that $g_{x,0}$ is the standard euclidean metric $g_{0}$.
Let $\ast$ be the Hodge operator of $g_0$.
Then the connection $B$ satisfies
\begin{equation}
\label{eq:4.2.1}
F_{B}= - \ast ( \Omega_x \wedge F_B). \spn{4.2.1}
\end{equation}
Define a degree $n-4$, constant form $\Phi_{S,x}$ on $T_xM$ as follows:
let $x_1, \cdots, x_n$ be any euclidean coordinates of $T_xM$ such that
$x_1, \cdots, x_{n-4}$ are tangent to $T_xS$;  then
$$
\Phi_{S,x} = dx_1\wedge \cdots \wedge dx_{n-4}.
$$
Now we decompose $\Omega_x = \alpha \Phi_{S,x} + \Omega_0$, where
$\alpha$ is a constant and $\Omega _0 |_{T_xS}=0$.

Since $v\rfloor F_B =0$ for any $v\in T_xS$, by taking
a gauge transformation if necessary, we may assume that $B= \pi_L^*B_L$
for some nontrivial connection $B_L$ on $L$, where $L$ is the orthogonal
complement of $T_xS$ and $\pi_L$ is the orthogonal
projection from $T_xM $ onto $L$. Then (4.2.1) becomes
\begin{eqnarray}
\label{eq:4.2.2}
F_{B_L}& = &- \alpha \ast_L F_{B_L}, \spn{4.2.2}\\[2ex]
0&=&\ast (\Omega_0 \wedge F_B),  \spn{4.2.3}
\end{eqnarray}
where $\ast_L$ is the Hodge operator of $L$.

Since $F_{B_L }\not=0$, we deduce from (4.2.2) that $\alpha = \pm 1$.
The corollary is proved.
\enddemo 
 
\specialnumber{4.2.3} \proclaim{Theorem}
\label{th:4.2.3}
Let $(M, g)$ be a compact Riemannian manifold{\rm ,} $\Omega$ be
a closed form of degree $n-4$ and $\{A_i\}$ be a sequence of $\Omega$\/{\rm -}\/anti\/{\rm -}\/self\/{\rm -}\/dual
instantons{\rm .} Then by taking a subsequence if necessary{\rm ,} $A_i$
converges to an admissible $\Omega$\/{\rm -}\/anti\/{\rm -}\/self\/{\rm -}\/dual instanton $A$
with the blow\/{\rm -}\/up locus $(S,\Theta)${\rm ,} such that {\rm (1)} $S$ is rectifiable
and $\Omega|_S$ is one of its volume forms induced by $g${\rm .} In particular{\rm ,} $S$
carries a natural orientation\/{\rm ; (2)} $\frac{1}{8\pi^2} \Theta$ is integer\/{\rm -}\/valued\/{\rm ;}
{\rm (3)}  $C_2(S, \Theta)$ is closed in $M$, where $C_2(S,\Theta)$ is an integral current defined by
\begin{equation}
\label{eq:4.2.3}
C_2(S, \Theta) (\varphi ) = \frac{1}{8\pi^2} \int _S (\varphi, \Omega|_S) \Theta d(H^{n-4}\lfloor S),
\spn{4.2.4}
\end{equation}
where $\varphi$ is any smooth form with compact support in $M${\rm .} Moreover, as currents{\rm ,}  we have
\begin{equation}
\label{eq:4.2.4}
\lim_{i\to \infty} C_2(A_i) = C_2(A) + C_2(S, \Theta),  \spn{4.2.5}
\end{equation}
where $C_2(A)$ is as defined in Corollary {\rm 2.3.2.}
\endproclaim

\demo{{R}emark {\rm 4}}
Applying (4.2.4) to the smooth form $4\pi^2 \Omega$, we obtain the conservation
of the action:
$$
  \lim_{i \to \infty} \int _M \mid F_i \mid^2 dV_g = \int _M \mid F_A \mid^2 dV_g
  +\int_S\Theta d(H^{n-4} \lfloor S).
$$
\enddemo

The rest of this section is devoted to the  proof of Theorem 4.2.3.
We will adopt the notations in  the proof of Proposition 4.1.1 and Corollary 4.2.2.

It is clear that (1) follows from 
Proposition 3.3.3, Proposition 4.2.1, Corollary 4.2.2 and results of the  last chapter,
so it suffices to prove (2) and (3).
 
First we show that the density $\frac{1}{8 \pi^2} \Theta (\mu, \cdot)$
is integer-valued. Let $x$ be any point in $S$ such that (\ref{eq:3.2.19})
holds and there is a unique tangent space $T_xS$. Then (\ref{eq:4.1.3})
holds. Now,
\begin{equation}
  \label{eq:4.2.5}
  \Theta (\mu, x) = \lim_{i\to\infty} \int_{B_1(0,g_{x,0})} \mid
    F_{A_{i,x,\lambda_i}} \mid^2 d V_{x,\lambda_i} . \spn{4.2.6}
\end{equation}
Since $A_{i,x,\lambda_i}$ converges to zero uniformly on any
compact subset away from $V=T_xS$, for any $z^{\prime} \in V \cap
B_1 (x,g_{x,0})$, $A_{i,x,\lambda_i} \mid_{\{z^{\prime}\}\times
    V^\bot \cap B_{\sqrt{1-\mid z^{\prime}\mid^2}} (0,g_{x,0})}$
  converges to zero uniformly away from $(z^{\prime},0)$. Then by
  the standard transgression arguments, we can deduce
  \begin{equation}
    \label{eq:4.2.6}
    \lim_{i\to\infty} \frac{1}{8\pi^2}\int_{\{z^{\prime} \} \times V^\bot \cap
      B_{\sqrt{1-\mid z^{\prime}\mid^2}} (0,g_{x,0})} \tr
      (F_{A_{i,x,\lambda_i}}\wedge F_{A_{i,x,\lambda_i}}) \in \mathbb{Z}. \qquad\quad\spn{4.2.7}
  \end{equation}
Clearly, the limit on the right of (\ref{eq:4.2.6}) is a
topological number and does not depend on $z^{\prime}$.

For simplicity, we will denote by $F_{A_{i,x,\lambda_i}}^V$ the
curvature of the restricted connection $A_{i,x,\lambda_i}|_{z'\times V^\bot}$.
Since $A_{i,x,\lambda_i}$ is $\tau_{\lambda_i}^*{\rm exp}^*\Omega$-anti-self-dual 
with respect to $g_{x,\lambda_i}$ and $\lim_{i \to\infty} g_{x,\lambda_i} =
g_{x,0}$, we obtain
\begin{eqnarray}
\label{eq:4.2.7}
&&\spn{4.2.8}\\
&&\hskip-30pt \frac{1}{8\pi^2} |F_{A_{i,x,\lambda_i}}|^2 dV_{x,\lambda_i} =
 -\frac{1}{8\pi^2} \tr (F_{A_{i,x,\lambda_i}}
\wedge F_{A_{i,x,\lambda_i}})\wedge \tau_{\lambda_i}^*{\rm exp}^*\Omega\nonumber \\[2ex]
&=& \frac{1}{8\pi^2} \left(-\tr (F_{A_{i,x,\lambda_i}}^V \wedge
F_{A_{i,x,\lambda_i}}^V ) \phantom{\sum^4_1} \right.\nonumber \\ [2ex]
&&\hskip.45in \left. + (O(1)
\sum^{n-4}_{\alpha =1} \mid \frac{\partial}{\partial
z_{\alpha}} \rfloor F_{A_{i,x,\lambda_i}} \mid + o(1)|F_{A_{i,x,\lambda_i}}|)\mid F_{A_{i,x,\lambda_i}}
 \mid \right)dV_{x,\lambda_i},\nonumber
\end{eqnarray}  
where $o(1)$ denotes a quantity which converges to zero as $i$ tends to infinity.
Together with (\ref{eq:4.2.6}) and (\ref{eq:4.1.3}), this implies
\begin{eqnarray*}
  \frac{1}{8\pi^2} \Theta (\mu,x) ]
& = & \lim_{i\to \infty} \frac{1}{8 \pi^2} \int_{B_{1}(0,g_{x,0})}
       \mid F_{A_{i,x,\lambda_i}} \mid^2 dV_{x,\lambda_i}\\[2ex]
&=& \lim_{i\to \infty} \int_{V \cap B_{1}(0,g_{x,0})} d(H^{n-4}\lfloor V) \\ [2ex]
&&\cdot \
\left (\frac{1}{8\pi^2}\int_{\{z'\} \times V^\bot \cap B_{\sqrt{1-\mid z\mid^2}}
  (0,g_{x,0})} \tr ( F_{A_{i,x,\lambda_i}}\wedge F_{A_{i,x,\lambda_i}}) \right ).
\end{eqnarray*}  
Hence, by (4.2.6), $\frac{1}{8 \pi^2} \Theta (\mu,\cdot)$ is integer-valued.

Next we show that $C_2(S, \Theta)$ is closed, i.e., for any smooth form
$\psi$ of degree $n-5$ and with compact support in $M$,
\begin{equation}
\label{eq:4.2.8}
\partial C_2(S, \Theta) (\psi)= C_2(S, \Theta) (d\psi) = 0. \spn{4.2.9}
\end{equation}
This will follow from (4.2.4) and
Corollary 2.3.2,
since
$$
\int_{M}d\psi\wedge\tr(F_{A_{i}}\wedge F_{A_{i}})=0\quad\hbox{for
  any } i.
$$
We also have
$$\lim_{i\to \infty} \int_M \tr(F_{A_i})\wedge \tr(F_{A_i})
=\int_M \tr(F_{A})\wedge \tr(F_{A}).$$
Therefore, it suffices to prove that
by taking a subsequence if necessary, for any smooth $\varphi$ 
of degree $n-4$,
\begin{eqnarray}
\label{eq:4.2.9}
&& \spn{4.2.10}\\
\noalign{\vskip4pt}
&& \hskip-12pt\frac{1}{8\pi^2}\lim_{i\to\infty}\int_M\varphi\wedge\tr(F_{A_{i}}\wedge
F_{A_{i}})=\frac{1}{8\pi^2}\int_M\varphi\wedge \tr(F_A\wedge F_A) + C_2(S,\Theta)(\varphi). 
\nonumber
\end{eqnarray}  
Define currents $T_i$ by
$$
T_i(\varphi)=\frac{1}{8\pi^2}\int_{M}\varphi\wedge\left (\tr(F_{A_{i}}\wedge
F_{A_{i}})-\tr(F_A\wedge F_A)\right );
$$
then by Proposition~\ref{prop:2.3.1}, $\partial T_i=0$.
Moreover, the total mass of $T_i$ is uniformly bounded; i.e., for any
$\varphi$ with $||\varphi||_{C^{0}}\leq 1$,
\begin{equation}
\label{eq:4.2.10}
|T_i(\varphi)|\leq\frac{1}{8\pi^2} \int_{M}\left(|F_{A_i}|^2-|F_A|^2\right)\,
dV_g\leq \Lambda. \spn{4.2.11}
\end{equation}
This implies, after taking a subsequence if necessary, that $T_i$
converges weakly to a closed current $T$.  Clearly, the mass of $T$
is also bounded by $\Lambda$ and  we have $\partial T =0$. Hence, 
by Theorem 3.2.1 in [Si2], $T$
is rectifiable; more precisely, there is a rectifiable set $S'$
with orientation vector $\eta:S'\to\Lambda^{n-4}T^\ast S'$ and a density
function $\Theta'(x)$, such that
$$
T(\varphi)= \frac{1}{4\pi^2}\int_{S'}(\varphi,\eta)\Theta'\, d(H^{n-4}\lfloor S').
$$
Take $\varphi$ to be $f \Omega$, where $f$ is a smooth function
with compact support; then
\begin{equation}
\label{eq:4.2.11}
T(f\Omega ) = \frac{1}{4\pi^2}\int_{S'}f(\Omega,
\eta)\Theta'(x) d(H^{n-4}\lfloor S'). \spn{4.2.12}
\end{equation}
On the other hand, since $\tr(F_{A_i})$ converges to
$\tr(F_A)$ uniformly on $M$, we have
\begin{eqnarray}
\label{eq;4.2.12}
\quad\quad T(f\Omega ) 
&=&\lim_{i\to \infty} T_i(f\Omega )\spn{4.2.13}\\[2ex]
&= &\frac{1}{8\pi^2} \lim_{i\to \infty}\int_{M} f\Omega \wedge \left 
(\tr(F_{A_i}\wedge F_{A_i}) - \tr (F_A\wedge F_A) \right )\nonumber\\[2ex]
&=&\frac{1}{8\pi^2}\lim_{i\to \infty} \int_{M}f \left(|F_{A_{i}}|^2-|F_A|^2\right) dV_g\nonumber\\[2ex]
&=&\frac{1}{8\pi^2} \int_{S}f (x) \Theta(\mu,x) d(H^{n-4}\lfloor S).  
\end{eqnarray}
Comparing this with (\ref{eq:4.2.11}), we conclude that $S'= S$
and $\Theta(\mu,\cdot)= (\Omega,\eta)\Theta'$. Finally,
since $\Omega_S$ is one of the volume forms of $S$, we
obtain that $(\Omega,\eta)=1$ and consequently, $T=C_2(S,\Theta)$.
This finishes the proof of Theorem 4.2.3.

\demo{{R}emark {\rm 5}}
We need the compactness of $M$ and the closedness of $\Omega$ only to derive
an {\it a priori} bound on $YM(A_i)$ in the above proof. 
\enddemo

\demo{{\rm 4.3.}\ Calibrated geometry and blow-up loci}
Let $(M,g)$ be an $n$-dimensional Riemannian manifold and $\Omega$ 
be a closed form of degree $n-4$. We further assume that for any
$x\in M$ and subspace $F$ of $T_xM$ of codimension $4$,
$\Omega |_F \le dV_F$, where $dV_F$ denotes the induced
volume form on $F$ by $g$. Following [HL], we say that $(F, dV_F)$ is calibrated by 
$\Omega$ if $\Omega |_F = dV_F$. Moreover, if $\Phi=(S,\xi, \Theta)$ is
an integral current with orientation $\xi$ and density $\Theta$,
where $S$ is the support of $\Phi$ and rectifiable,
then we say that $\Phi$ is $\Omega$-calibrated if $(T_xS, \xi(x))$
is calibrated by $\Omega$ for $H^{n-4}$-a.e. $x\in S$.

The following lemma is trivial. \enddemo
\specialnumber{4.3.1} \proclaim{Lemma}
\label{lem:4.3.1}
Any integral current calibrated by $\Omega$ is minimizing in its homology
class{\rm .} In particular{\rm ,} its generalized mean curvature vanishes{\rm .}
\endproclaim
\demo{Proof} Let $\Phi=(S,\xi, \Theta)$ be an integral current calibrated by 
$\Omega$, and $\Psi=(S',\xi', \Theta')$ be another integral current homologous to $\Phi$;
i.e., there is a current $R$ of degree $n-5$ such that for any smooth form
$\varphi$ on $M$,
\begin{equation}
\label{eq:4.3.1}
\int _S (\varphi, \xi) \Theta d H^{n-4} - \int_{S'} (\varphi, \xi') \Theta' dH^{n-4}
= R(d\varphi). \spn{4.3.1}
\end{equation}
By our assumption, $(\Omega, \xi') \le 1$ and $(\Omega, \xi)=1$. Hence,
\begin{equation}
\label{eq:4.3.2}
\int_S \Theta dH^{n-4} \le \int_{S'} \Theta' dH^{n-4} + R(d\Omega) =
\int_{S'} \Theta' dH^{n-4}, \hskip.5in\spn{4.3.2}
\end{equation}
and it follows that $\Phi$ is minimizing.
\enddemo

Clearly, such a $\Phi$ is determined by $S$ with multiplicity $\Theta$.
We will also call $(S,\Theta)$ an $\Omega$-calibrated cycle. 
It is known from the geometry measure theory that for
such a cycle, $S$ is regular in an open and dense subset. In fact, it follows from
[Am] that $S$ can be decomposed as $\bigcup _a S_a$, such that each $S_a$ is closed and
smooth outside a closed subset of Hausdorff codimension at least two
and $\Theta$ restricts to a positive integer on each $S_a$.

\specialnumber{4.3.2} \proclaim{Theorem}
\label{th:4.3.2}
Let $(M, g)${\rm ,} $\Omega$ be as above{\rm ,} and 
$\{A_i\}$ be a sequence of $\Omega$\/{\rm -}\/anti\/{\rm -}\/self\/{\rm -}\/dual
instantons{\rm .} Further assume that either $M$ is compact or the  $YM(A_i)$
are uniformly bounded{\rm .} Then by taking a subsequence if necessary{\rm ,} $A_i$
converges to an admissible $\Omega$\/{\rm -}\/anti\/{\rm -}\/self\/{\rm -}\/dual instanton $A$
with the blow\/{\rm -}\/up locus $(S, \Theta)${\rm ,} such that $(S,\Theta) $ is
is an $\Omega$\/{\rm -}\/calibrated cycle{\rm ,} and 
\begin{equation}
\label{eq:4.3.3}
\lim_{i\to \infty} C_2(A_i) = C_2(A) + C_2(S, \Theta). \spn{4.3.3}
\end{equation}
\endproclaim

This follows from Theorem 4.2.3 and the  discussions above. 

In the case of Hermitian-Yang-Mills connections, we have

\specialnumber{4.3.3} \proclaim{Theorem}
\label{th:4.3.3}
Let $(M, g)$ be a complex $m$\/{\rm -}\/dimensional compact K{\rm \"{\it a}}hler manifold 
with the K{\rm \"{\it a}}hler form $\omega ${\rm ,}
and $\{A_i\}$ be a sequence of Hermitian\/{\rm -}\/Yang\/{\rm -}\/Mills connections on a
given unitary bundle $E${\rm .} Then by taking a subsequence if necessary{\rm ,} $A_i$
converges weakly to an admissible Hermitian\/{\rm -}\/Yang\/{\rm -}\/Mills connection $A$
with the blow-up locus $(S,\Theta)${\rm ,} such that $S=\bigcup _\alpha S_\alpha$
and $\Theta |_{S_\alpha} = 8\pi^2 m_\alpha${\rm ,} where each $S_\alpha$ is a holomorphic 
subvariety in $M$ and $m_\alpha$ is a positive integer{\rm .} Moreover{\rm ,} for any smooth $\varphi${\rm , } 
\begin{equation}
\label{eq:4.3.3}
\lim_{i\to \infty} \int_M \varphi \wedge C_2(A_i) = 
\int_M \varphi \wedge C_2(A) + \sum_\alpha m_\alpha 
\int _{S_\alpha} \varphi.  \spn{4.3.4}
\end{equation}
\endproclaim
\demo{Proof} By Theorem 4.3.2, we may assume that $A_i$ converges
to an admissible Hermitian-Yang-Mills connection $A$ with an 
$\frac{\omega^{m-2}}{(m-2)!}$-calibrated cycle $(S,\Theta)$ as its blow-up
locus. It suffices to show that $(S,\Theta)$ is a holomorphic cycle.

A straightforward computation shows that for any $x \in M$ and subspace 
$F \subset T_xM$ of codimension $4$, $\frac{\omega^{m-2}}{(m-2)!}|_F \le dV_F$
and the equality holds if and only if $F$ is a complex subspace in $T_xM$.
Therefore, $T_xS$ is a complex subspace in $T_xM$ for $H^{2m-4}$-a.e. $x\in S$.
Since $C_2(S,\Theta)$ is a closed integral current, it follows from a
result of J. King [Ki] or Harvey and Shiffman [HS] that there are holomorphic subvarieties
$S_\alpha$ and positive integers $m_\alpha$ such that
$$C_2(S, \Theta) (\varphi) = \sum _\alpha m_\alpha \int_{S_\alpha} \varphi$$
for any $\varphi$. The theorem is proved.
\enddemo

\demo{{R}emark {\rm 6}}
Let $A$ be the Hermitian-Yang-Mills connection in the above theorem. 
It follows from a result of Bando and Siu [BS] that there is a gauge transformation
$\sigma$ on $M\backslash S$ such that $\sigma(A)$ extends to a smooth Hermitian-Yang-Mills connection
outside a holomorphic subvariety in $M$ of codimension at least three. In fact,
the (0,1)-part of $A$ induces a holomorphic structure on the underlying complex vector
bundle. Then the induced holomorphic bundle on $M\backslash S$ extends to
a coherent sheaf which is locally free outside a subvariety of codimension at least
three.
\enddemo

\demo{{\rm 4.4.}\ Cayley cycles and complex anti-self-dual instantons}
In this section, we assume that $(M,g)$ is a Calabi-Yau 
$4$-fold with
the K\"ahler form $\omega$ and a holomorphic $(4,0)$-form $\theta$.
We normalize 
\begin{equation}
\label{eq:4.4.1}
\theta \wedge \overline \theta = \frac{\omega^4}{4!}. \spn{4.4.1}
\end{equation}
As in Section 1.3, we put
\begin{equation}
\label{eq:4.4.2}
\Omega = {2}(\theta + \overline \theta) + \frac{\omega^2}{2}. \spn{4.4.2}
\end{equation} \enddemo

\specialnumber{4.4.1} \proclaim{Lemma}
\label{lem:4.4.1}
For any $4$\/{\rm -}\/dimensional subspace $L\subset TM${\rm ,}  $\Omega\mid_L\le dV_L${\rm .}
\endproclaim

\demo{Proof} 
This should be well-known. For the reader's convenience, we include an
elementary proof here.
Without loss of generality, we may assume that
$M = \CC ^4$ and $L \subset \CC^4$. In any euclidean coordinates
$z_1 , \cdots , z_4$ of $\CC^4$,
$$
\omega = \frac{\sqrt{-1}}{2}
\sum^4_{i=1} \, dz_i \wedge \, d\bar{z}_i,~~~ \theta =
\frac{1}{4} dz_1 \wedge dz_2 \wedge dz_3 \wedge dz_4 . 
$$
Let $J$ be the standard complex structure on $\CC ^4$.
Then $\dim_{\RR} JL\cap L =0$ or $2$ or $4$.  
Since $\pi_{L} \cdot (J|_{L})$ is skewsymmetric, where $\pi_{L}$ denotes the
orthogonal projection onto $L$,
one can choose an orthonormal basis $\left\{ u_1, u_2, u_3, u_4 \right\}$
of $L$, such that
\begin{eqnarray}
\label{eq:4.4.3}
  J u_1 &=& u^{\perp}_1 + \lambda u_2 , \quad
  J u_2 = u^{\perp}_2 - \lambda u_1 , \spn{4.4.3}\\[2ex]
  J u_3& = &u^{\perp}_3 + \lambda'  u_4 , \quad
  J u_4 = u^{\perp}_4 - \lambda' u_3 ,\nonumber
\end{eqnarray}
where $u^{\perp}_1, u^{\perp}_2, u^{\perp}_3, u^{\perp}_4$ are in
the orthogonal complement $L^{\perp}$.  

First we assume that $\dim_{\RR} JL\cap L =0$, i.e., $L$ is totally real.
Then $| \lambda | <1$, $| \lambda' | <1$.  Define 
\begin{eqnarray}
\label{eq:4.4.4}
v_1&=&u_1 \, , \, v_2 = \sqrt{1- \lambda^2} u_2 - \frac{\lambda u^{\perp}_1}{\sqrt{1-
    \lambda^2}} \, , \spn{4.4.4}\\[2ex]
v_3&=&u_3 \, , \, v_4 = \sqrt{1- {\lambda'}^2} u_4 -
\frac{\lambda u^{\perp}_3}{\sqrt{1- {\lambda'}^2}} \,\nonumber .
\end{eqnarray}
Then $\left\{ v_i,J_0 v_i \right\}_{1 \leq i \leq 4}$
is an orthonormal basis of $\CC^4$, such that
\begin{equation}
  \label{eq:4.4.5}
  J v_2 = \frac{u^{\perp}_2}{\sqrt{1- \lambda}^2} \in
  L^{\perp} \, , \, 
  J v_4 = \frac{u^{\perp}_4}{\sqrt{1- {\lambda'}^2}} \in
  L^{\perp} \, . \spn{4.4.5}
\end{equation}
Let $\left\{ v^\ast_i , (J v_i )^\ast \right\}$ be its the dual
basis. Put $\varphi^\ast_i = v^\ast_i - \sqrt{-1} (J v_i)^\ast $. Then
\begin{eqnarray}
\label{eq:4.4.6}
\omega &=& \frac{\sqrt{-1}}{2} \sum^4_{r=1} \varphi^\ast_i
 \wedge \overline{\varphi^\ast_i}\spn{4.4.6} \\ [2ex]
\theta &=& \frac{1}{4} e^{\sqrt{-1} \gamma} \varphi^\ast_1 \wedge \varphi^\ast_2
\wedge \varphi^\ast_3 \wedge \varphi^\ast_4 \, , \, \gamma \in \RR .\spn{4.4.7}  
\end{eqnarray}
Using (4.4.5), (4.4.6) and (4.4.7), one shows
\begin{eqnarray}
\label{eq:4.4.8}
\theta |_L&=& \frac{e^{\sqrt{-1} \gamma}}{4} \sqrt{(1-\lambda^2) (1-{\lambda'}^2)}
\,dV_L,\spn{4.4.8}\\[2ex]
\omega^2|_L&=& 2 \lambda \lambda' \,dV_L. \spn{4.4.9}
\end{eqnarray}
It follows that 
$$\Omega |_L = \left ( \cos \gamma \sqrt{(1-\lambda^2)(1-{\lambda'}^2)} + \lambda\lambda'
\right ) dV_L,$$
and so $\Omega |_L \le dV_L$.
 
Two other   cases can be easily reduced to this case by perturbing $L$ slightly. \phantom{needtopush}
\enddemo

Following [HL], we see that $(S,\Theta)$ is a Cayley cycle if it
is calibrated by the $\Omega$. In this case, 
for $H^4$-a.e. $x\in S$, the tangent space $T_xS$ is a
Cayley plane in $T_xM$.

The following observation is of considerable interest, though simple. 
\specialnumber{4.4.2} \proclaim{Proposition}
\label{prop: 4.4.2}
Let $(S,\Theta)$ be a Cayley cycle{\rm .} Then 
$ C_2(S,\Theta)\cdot [\theta]$
is a nonnegative real number{\rm .} Moreover{\rm ,} $C_2(S,\Theta)\cdot [\theta]=0$
if and only if $(S, \Theta)$ is a holomorphic cycle{\rm .}
\endproclaim

\demo{Proof} We adopt the notation  in the proof of Lemma 4.4.1. 
If $T_xS$ exists and is a Cayley subspace, then $e^{\sqrt{-1} \gamma}=1$
and $\lambda=\lambda'$; this implies 
$$(\theta |_L , \xi_S) = (1-\lambda^2) \ge 0.$$
The first statement is proved. If $C_2(S,\Theta)\cdot [\theta]=0$,
then $\lambda ^2 =1$; i.e., $T_xS$ is a complex subspace. Then the proposition
follows from the main result in [HS] or [Ki].
\enddemo

\demo{{R}emark {\rm 7}}
Since $S$ is of codimension greater than $3$, we may simply define $C_1(S, \Theta)=0$.
Then the above result can be rephrased as 
$$\left (2C_2(S, \Theta) - {r-1\over r} C_1(S, \Theta)^2 \right ) \cdot [\theta] \ge 0,$$
and the equality holds if and only if $(S, \Theta)$ is a holomorphic cycle.
This is analogous to (1.3.8).
\enddemo
\demo{{R}emark {\rm 8}}
Similarly, one can easily show that for any Cayley cycle $(S,\Theta)$,
$C_2(S,\Theta) \cdot [\omega^2] \ge 0$. Moreover, the equality
holds if and only if $S$ is special Langrangian.
\enddemo

The next result follows from Theorem 4.2.3 and the above discussions.
\specialnumber{4.4.3} \proclaim{Theorem}
\label{th:4.4.3}
Let $(M,g)$ be a compact Calabi\/{\rm -}\/Yau $4$-fold with K{\rm \"{\it a}}hler form $\omega$
and a holomorphic $(4,0)$\/{\rm -}\/form $\theta${\rm .} Let $\{ A_i\}$ be a sequence of
complex anti\/{\rm -}\/self\/{\rm -}\/dual instantons{\rm .} Then by taking a subsequence
if necessary{\rm ,} the $A_i$ converge to an admissible
complex anti\/{\rm -}\/self\/{\rm -}\/dual instanton $A$ with the blow\/{\rm -}\/up
locus $(S,\Theta)${\rm ,} such that $(S,\Theta)$ is a Cayley cycle
and
$$
\lim_{i \to \infty} C_2(A_i) = C_2(A) + C_2(S,\Theta) .
$$
\endproclaim

\demo{{R}emark {\rm 9}}
The above theorem also holds for general ${\rm Spin}(7)$-manifolds, which
contain Calabi-Yau $4$-folds as special examples of ${\rm Spin}(7)$-manifolds.
\enddemo

Next we assume that $M$ is an $8$-dimensional compact manifold which admits Calabi-Yau
structures. A Calabi-Yau structure on $M$ is given by a complex structure $J$,
a K\"ahler metric $g$ compatible with $J$ and a holomorphic $(4,0)$-form $\theta$
satisfying (4.4.1). We denote by ${\cal M}$ the moduli space of all Calabi-Yau
structures modulo obvious equivalence relations. 

Let $E$ be a fixed $U(r)$-bundle over $M$. We define 
\begin{eqnarray}
\noalign{\vskip-3pt}
&&\label{eq:4.4.11} \spn{4.4.10}\\
\noalign{\vskip-3pt}
&&\hskip-3pt {\cal M}(E) \nonumber
\\
= &&\hskip-3pt  \{ (J,g,\theta) \in {\cal M} \,|\, \left (2C_2(E)-\frac{r-1}{r}C_1(E)^2 
\right ) \cdot \varphi \ge 0, \hbox{ {for} $\varphi = [\theta]$ {\rm or}  $[\omega^2_g] \}$,}
\nonumber
\end{eqnarray}
where $\omega_g$ denotes the K\"ahler form of $g$,
$C_1(E)$ and $C_2(E)$ denote the first and second Chern character of $E$.
It is easy to show that ${\cal M}(E)$ is a connected, analytic variety.

For a fixed Calabi-Yau structure $(J,g,\theta)$, we denote by ${\cal Y}_{J,g,\theta}(E)$
the moduli space of all complex anti-self-dual connections of $E$ on
the Calabi-Yau $4$-fold $(M, J,g,\theta)$ modulo gauge transformations.
By Theorem 4.4.3, modulo gauge transformations,
its compactification $\overline {\cal Y}_{J,g,\theta}(E)$
consists of all triples $(A, S, \Theta)$ satisfying: (1) $A$ is an
admissible complex anti-self-dual instanton on $M$; (2) $(S,\Theta)$ is a
Cayley cycle; (3) $C_i(E) = [C_i(A)] + [C_i(S,\Theta)]$ in $H^{2*}(M, \RR )$
for $i=1,2$. The topology of $\overline {\cal Y}_{J,g,\theta}(E)$ is the one
determined by the convergence property given at the beginning of Section 3.1.

We define
$$\overline {\cal Y}(E) = \bigcup _{( J,g, \theta) \in {\cal M}} 
\overline {\cal Y}_{J,g,\theta}(E).$$
By Proposition 4.4.2, there is an obvious map $f_c: \overline {\cal Y}(E)
\mapsto {\cal M}(E)$. Denote by ${\cal M}_c(E)$ its image.
Then the next result follows from the same arguments as those in the proof of Theorem 4.4.3.
\specialnumber{4.4.4} \proclaim{Theorem}
\label{th:4.4.4}
The set ${\cal M}_c(E)$ is closed in ${\cal M}(E)${\rm .}
\endproclaim
It is easy to show that ${\cal M}(E)$ is an analytic variety.
We conjecture that {\it ${\cal M}_c(E)$ is an analytic subvariety
in ${\cal M}(E)$}.

\demo{{R}emark {\rm 10}}
All the   discussions above work as well for general anti-self-dual instantons.
We single out the complex anti-self-dual case because of its plausible
connection to the Hodge conjecture on holomorphic cycles on Calabi-Yau $4$-folds.
\enddemo

\demo{{\rm 4.5.}\ General blow-up loci}
Let $\{A_i\}$ be a sequence of smooth Yang-Mills connections which converge weakly to
an admissible Yang-Mills connection $A$ with blow-up locus $(S,\Theta)$.

\specialnumber{4.5.1} \proclaim{Theorem}
\label{th:4.5.1}
For any vector field $X$ with compact support in $M${\rm ,} 
\begin{equation}
\label{eq:4.5.1}
-\int _S \divv _S X \Theta \,dH^{n-4} = \int _M \left ( |F_A|^2
\divv X - 4 (F_A(\nabla X, \cdot ), F_A)\right ) dV_g, \quad \spn{4.5.1}
\end{equation}
where $(F_A(\nabla X, \cdot ), F_A)$ is defined in 
any local orthonormal basis $\{e_i\}$ of $M$ as
$$
\sum _{i,j=1}^n (F_A(\nabla_{e_i} X, e_j), F_A(e_i, e_j))
$$
and $\divv _SX$ denotes the
divergence of $X$ along $S${\rm .} That is{\rm ,}  if $T_pS$ exists and $\{v_i \}$ is any orthonormal
basis of $T_pS${\rm ,} $\divv _S X (p)= \sum_{i=1}^{n-4} ( \nabla_{v_i}X, v_i ) (p)${\rm .}
\endproclaim
\demo{Proof} As above, $c$ always denotes a uniform constant.
Since $S$ is rectifiable, we can find a countable set of submanifolds
$\{N_\alpha\}$ such that $S= S_0 \bigcup _\alpha S_\alpha$,
where $S_\alpha = N_\alpha \cup S$ and $H^{n-4}(S_0)=0$ (cf.\ [Si2]).
Moreover, we may assume that $T_xS = T_xN_\alpha$ for $H^{n-4}$-a.e.
$x\in S_\alpha$. 

Fixing any $\delta >0$, we can arrange $N_\alpha$ such that 
for some $\alpha _\delta> 0$,
\begin{eqnarray}
\label{eq:4.5.2}
S_\alpha \cap S_{\alpha'} &=& \emptyset,~~{\rm for}~\alpha, 
\alpha' \le \alpha_\delta;\spn{4.5.2}\\[2ex]
H^{n-4}(\bigcup _{\alpha > \alpha _\delta}S_\alpha ) &\le& \delta. \nonumber
\end{eqnarray}
It follows that by taking a subsequence if necessary, we have
\begin{equation}
\label{eq:4.5.3}
\lim_{\varepsilon \to 0}\lim_{i\to \infty} 
\int _{B_{\varepsilon} (\bigcup _{\alpha > \alpha _\delta}S_\alpha)}
|F_{A_i}|^2 dV_g \le 2 \delta. \spn{4.5.3}
\end{equation} 
Since $\delta$ can be arbitrarily small,   it
 suffices to prove that for each $\alpha \le \alpha_\delta$,
\begin{eqnarray}
\label{eq:4.5.4}
&&\spn{4.5.4}\\
&&\lim_{\varepsilon\to 0}\lim_{i\to \infty} \int _{B_{\varepsilon}
(S_\alpha)} \left ( |F_{A_i}|^2 \divv X - 4
\sum_{k,l}(F_{A_i}(\nabla_{e_k} X, e_l), F_{A_i}(e_k, e_l))\right ) dV_g\nonumber \\[2ex]
&&\quad =\int _{S_\alpha} \divv _S X \Theta \,dH^{n-4}. \nonumber
\end{eqnarray}

Without loss of generality, we may assume that $e_1,\cdots, e_{n-4}$
are tangent to $N_\alpha$, while $e_{n-3}, \cdots, e_n$ are normal
to $N_\alpha$. Then it follows from Lemma 3.3.2 that
(4.5.4) is the same as 
\begin{eqnarray}
\label{eq:4.5.5} 
&&  \spn{4.5.5}\\
&& \hskip-12pt \lim_{\varepsilon\to 0}\lim_{i\to \infty}\int _{B_{\varepsilon}
(S_\alpha)}\! \left ( |F_{A_i}|^2 \divv ^\perp X - 4
\sum_{k,l=n-3}^n (F_{A_i}(\nabla_{e_k} X, e_l), F_{A_i}(e_k, e_l))
\!\right ) dV_g \nonumber\\[2ex]
= && \hskip-12pt
0,\nonumber
\end{eqnarray}
where $\divv^\perp X = \sum_{k=n-3}^n g(\nabla_{e_k}X , e_k)$ is 
the divergence of $X$ in normal directions of $N_\alpha$. 

Write $\nabla_{e_k}X = X_{i,k} e_i$, then 
$\divv^\perp X = \sum_{l=n-3}^n X_{l,l}$ and (4.5.5) becomes
\begin{eqnarray}
\label{eq:4.5.6}
&& \spn{4.5.6}\\
&&\hskip-.5in \lim_{\varepsilon\to 0}\lim_{i\to \infty} \int _{B_{\varepsilon}
(S_\alpha)} \sum_{k,l=n-3}^n X_{k,l}\nonumber\\[2ex]
&&\quad \cdot \  \left ( |F_{A_i}|^2 \delta_{kl}
- 4 \sum_{j=n-3}^n (F_{A_i}(e_k, e_j), F_{A_i}(e_l, e_j))
\right ) dV_g  =  0.\nonumber
\end{eqnarray} 
By taking subsequences if necessary, we may assume that there
are measures $\mu_{kl}$ ($k,l=n-3,\cdots, n$), defined by \smallbreak
\vglue-8pt
\begin{equation}
\label{eq:4.5.7}
\mu_{kl}(h) = \lim_{i\to \infty} \int _{B_{\varepsilon} (N_\alpha)} h
( |F_{A_i}|^2 \delta_{kl}
- 4 \sum_{j=n-3}^n (F_{A_i}(e_k, e_j), F_{A_i}(e_l, e_j))) dV_g , \spn{4.5.7}
\end{equation}
where $h$ is any function with compact support in $B_{\varepsilon} (N_\alpha)$.
It follows from the monotonicity (Theorem 2.1.2) that for any 
$x\in S$ and $r$ sufficiently small,
$$
\mu_{kl} (B_r(x)) \le c e^{ar^2} r ^{n-4}.
$$
Hence, in order to prove (4.5.6), it suffices to show
that the upper-density $\overline 
\Theta (\mu_{kl}, x)$ ($k,l=n-3,\cdots, n$)
vanishes for $H^{n-4}$-a.e. $x\in S_\alpha$, where
\begin{equation}
\label{eq:4.5.8}
\overline \Theta (\mu_{kl}, x) = \lim\sup_{r\to 0}
r^{4-n} |\mu_{kl} (B_r(x))|. \spn{4.5.8}
\end{equation}

We will prove (4.5.8) by contradiction. If (4.5.8) is false, 
there is an $S'_\alpha \subset S_\alpha$ such that 
$H^{n-4}(S'_\alpha) > 0$ and for some $k, l$,
$\overline \Theta (\mu_{kl}, x) > 0$
for any $x \in S'_\alpha$. 
By orthogonal transformations, we may assume that
$k=l=n$. We can also have
that for $x\in S'_\alpha$,
the tangent space $T_xS= T_xS'_\alpha$ exists and
\begin{equation}
\label{eq:4.5.9}
\lim_{r\to 0} r^{4-n} \int_{B_r(x)} |F_A|^2 dV_g = 0. \spn{4.5.9}
\end{equation}
Then, by using the arguments in the proof of Lemma 4.1.2
and taking a subsequence if necessary,
we can find $\varepsilon_i, r_i>0$ with $\lim \varepsilon_i =0$
and $\lim \frac{r_i }{\varepsilon_i} =0$, $x_i \in S'_\alpha$,
such that \smallbreak
\vglue-6pt
\begin{equation}
\label{eq:4.5.10}
r_i^{4-n}\left |\int_{B_{r_i}(x_i)}\left (|F_{A_i}|^2 
- 4 \sum_{j=n-3}^n (F_{A_i}(e_n, e_j), F_{A_i}(e_n, e_j))
\right ) dV_g \right |  \ge \eta_0, \quad\spn{4.5.10}
\end{equation}
\begin{equation}
\lim_{i\to \infty} \varepsilon_i^{4-n} \int_{B_{\varepsilon_i}(x_i)} 
\sum_{j=1}^{n-4} |e_j\rfloor F_{A_i}|^2 dV_g  =  0. \spn{4.5.11}
\end{equation} 
For simplicity, we assume that $M\subset \RR^n$ and $g$ is flat. The general
case can be treated with slight modifications.
Put $B_i(y) = r_i A_i(x_i + r_i y)$. Then $B_i$ converges to zero outside
a subspace $\RR^{n-4}\times \{0\} = \lim_{i\to \infty} T_{x_i}N_\alpha$.

Let $X$ be a vector field with compact support in
$B_2(0)\subset \RR^n$. Since $B_i$ is Yang-Mills, we have that for 
any $j \le n-4$,
\begin{eqnarray*}
\int_{B_2(0)} |F_{B_i}|^2 X_{j,j} dV_g 
&=& -2\int_{B_2(0)} \sum _{k,l=1}^n (F_{B_i}(e_k, e_l), \nabla_{e_j} 
F_{B_i}(e_k, e_l)) X_j dV_g\\[2ex]
{\rm (Bianchi~ identity)}
&=&-4 \int_{B_2(0)} \sum _{k,l=1}^n (F_{B_i}(e_k, e_l), \nabla_{e_l} 
F_{B_i}(e_k, e_j)) X_jdV_g\\[2ex]
&=& 4\int_{B_2(0)} \sum _{k,l=1}^n (F_{B_i}(e_k, e_l), 
F_{B_i}(e_k, e_j)) X_{j,l}dV_g\\[2ex]
&\mapsto&0,~~~{\rm as}~i\to \infty.
\end{eqnarray*}   
Then we have
\begin{eqnarray*}
0&=& \int_{B_2(0)} \left(|F_{B_i}|^2 \divv X-4\sum^{n}_{k,l=1} 
(F_{B_i}(\nabla_{e_{k}}X, e_l), F_{B_i}(e_k,e_l))\right)\, dV_g\\[2ex]
&=&\int_{B_2(0)}\sum^{n}_{k,l=n-3} 
 X_{k,l}\left(|F_{B_i}|^2 \delta_{kl}- 4 \sum_j (F_{B_i}(e_k, e_j), F_{B_i}(e_l,e_j))\right)\, dV_g.
\end{eqnarray*}
Let $\eta$ be a nonnegative function on $\RR^1$ satisfying: $\eta(t) =1$
for $t\le 1$ and $\eta(t) =0$ for $t> \frac{4}{3}$.
Choose 
$$
X = \eta (|y'|) \eta (|y''|) y_n e_n , 
$$
where $y'=(y_1, \cdots, y_{n-4})$,
$y''=(y_{n-3}, \cdots ,y_n)$. Then the above implies
$$
\lim_{i\to \infty}\int_{B_1(0)}\left(|F_{B_i}|^2 - 
4 \sum_{j}(F_{B_i}(e_n, e_j), F_{B_i}(e_n,e_j))\right)\, dV_g =0.
$$
This contradicts   (4.5.10) and the theorem is proved.
\enddemo
 
We say that $A$ is stationary if the following holds 
for any vector field $X$ with compact support in $M$:
\begin{equation}
\label{eq:4.5.12}
\int_M\left(|F_A|^2 \divv X-4\sum^{n}_{i,j=1} (
  F_A(\nabla_{e_{i}}X, e_j), F_A(e_i,e_j))\right)\, dV_g =0\hskip.25in \spn{4.5.12}
\end{equation}
where $\{e_i\}$ is any orthonormal basis of $M$. If $A$
is a smooth Yang-Mills connection, this follows
from the first variation formula for Yang-Mills action.

\demo{{R}emark {\rm 11}} 
More generally, inspired by R. Schoen's notion of stationary harmonic maps,
we may define a stationary Yang-Mills connection as a weak solution 
of the Yang-Mills equation which satisfies (\ref{eq:4.5.12}). 
It is interesting to develop a regularity
theory for such weak solutions. But we will confine ourselves to 
admissible connections.
\enddemo

If $A$ is stationary, then the right side of (4.5.1) vanishes for any $X$.
This implies:

\specialnumber{4.5.2} \proclaim{{C}orollary}
\label{cor:4.5.2}
If $A$ is stationary{\rm ,} then $S$ is stationary{\rm ,} i{\rm .}e{\rm ., }
$S$ has no boundary in $M$ and its generalized mean curvature vanishes{\rm .}
\endproclaim

This also provides another proof of Theorem 4.3.2 with slightly weaker conclusion.

\section{Removable singularities of Yang-Mills equations}
\label{chap:5}

In this chapter, we investigate the extension problem of admissible
Yang-Mills connections. Since the extension
problem is local in nature, we may assume that $M$ is an open
subset in $\RR^n$ with a metric $g$, which may be nonflat.
\demo{{\rm 5.1.}\ Stationary properties of Yang-Mills connections}
Let $A$ be an admissible Yang-Mills connection as in Section 2.3,
and $r_p, c(p), a$ be as in Theorem 2.1.2. 
Then by the   arguments of  Section 2.1,
we have: \enddemo
\specialnumber{5.1.1} \proclaim{Proposition}
\label{prop:5.1.1}
Let $A$ be any admissible Yang\/{\rm -}\/Mills connection satisfying {\rm (4.5.2),} i{\rm .}e{\rm .,}
\begin{equation}
\label{eq:5.1.1}
\int_M\left(|F_A|^2 \divv X-4\sum^{n}_{i,j=1}\langle
  F_A(\nabla_{e_{i}}X, e_j), F_A(e_i,e_j)\rangle\right)\, dV_g =0 \hskip.5in \spn{5.1.1}
\end{equation}
where $\{e_i\}$ is any orthonormal basis of $M${\rm . }
Then for any $0 < \sigma < \rho < r_p${\rm ,}  
\begin{eqnarray}
\label{eq:5.1.2}
&& \hskip-36pt\rho^{4-n}e^{a\rho^{2}}\int_{B_{\rho}(p)}|F_A|^{2}dV_g-
\sigma^{4-n}e^{a\sigma^{2}}\int_{B_{\sigma}(p)}|F_A|^{2}dV_g \spn{5.1.2}\\[2ex]
&&\qquad\ge\  4\int_{B_\rho(p)\backslash B_{\sigma}(p)}r^{4-n}|\frac{\partial}{\partial r}\rfloor
F_A|^{2}dV_g. \nonumber
\end{eqnarray}
Moreover{\rm ,} if $M=\RR^n$ and $g$ is flat{\rm ,} then 
the equality holds in {\rm (5.1.2)} for $\rho\in(0,\infty)$ and $a=0${\rm .}
\endproclaim 

Next we prove that any admissible $\Omega$-anti-self-dual instantons are
stationary; i.e., they satisfy (5.1.1).

Let $A$ be an admissible $\Omega$-anti-self-dual instanton
with singular set $S=S(A)$.

Given any vector field $X$ with compact support in $M$, let $\phi_t:M\to M$ be 
its integral curve. As in Section
2.1, we define $A^t$ to be the
connection $\phi_t^*(A)$. Then by the same arguments as those
in Section 2.3, one can show that ${\rm Ch}_2(A^t)$
defines a closed $4$-form on $M$ in the sense of distribution.

First we claim that ${\rm Ch}_2(A^t)$ is independent of $t$, i.e., for
any closed $(n-4)$-form $\varphi$,
\begin{equation}
\label{eq:5.1.3}
\int_M\varphi\wedge\left({\rm Ch}_2(A^t)-{\rm Ch}_2(A)\right)=0. \spn{5.1.3}
\end{equation}
Since $\phi_t$ is an identity
near the boundary $\partial M$ of $M$, 
$$
{\rm Ch}_2(A^t)-{\rm Ch}_2(A)=0~~~{\rm near}~~\partial M.
$$

Without loss of generality, we may assume that the bundle $E$
under consideration is trivial over $M$. We constructed in
Section 2.3 a Chern-Simon $3$-form $\Psi$ such that
\begin{equation}
\label{eq:5.1.4}
d\Psi={\rm Ch}_2(A)\quad\hbox{on } M\backslash S, \spn{5.1.4}
\end{equation}
and for some uniform constant $c$,
\begin{equation}
\label{eq:5.1.5}
|\Psi(x)|\leqslant\frac{c}{d(x,S)^3},~~~~ x\in M\backslash S. \spn{5.1.5}
\end{equation}
Noticing that ${\rm Ch}_2(A^t)= \phi_t^* {\rm Ch}_2(A)$, we have
that for $\Psi_t = \phi_t^* \Psi$,
\begin{equation}
\label{eq:5.1.6}
d(\Psi_t -\Psi )={\rm Ch}_2(A^t)-{\rm Ch}_2(A) ~~~~\hbox{ in } 
M\backslash(S\cup\phi_t(S))  \hskip.5in \spn{5.1.6}
\end{equation}
and
\begin{equation}
\label{eq:5.1.7}
|\Psi_t-\Psi|(x)\leq \frac{2c}{d(x, S\cup\phi_t(S))^3}, 
~~~x\in M\backslash (S\cup\phi_t(S)).  \hskip.5in \spn{5.1.7}
\end{equation}
Furthermore, $\Psi_t-\Psi=0$ near $\partial M$ and for
$H^{n-4}$-a.e. $x \in S\cup\phi_t(S)$,
\begin{equation}
\label{eq:5.1.8}
\lim_{x\to x_0}d(x,S\cup\phi_t(S))^3\left(\Psi_t -\Psi\right)(x)=0 .\spn{5.1.8}
\end{equation}
Now (5.1.3) follows easily from (5.1.6)--(5.1.8)
and the same arguments as in the proof of Proposition 2.3.1.

\specialnumber{5.1.2} \proclaim{Proposition}
\label{prop:5.1.2}
Assume that $\Omega$ is a closed form of degree $n-4${\rm .}
Then any admissible $\Omega$\/{\rm -}\/anti\/{\rm -}\/self\/{\rm -}\/dual instanton $A$ on
$M$ is stationary{\rm . }
\endproclaim
\demo{Proof} For simplicity, we assume that $\tr(F_A)=0$.
The general case follows from identical arguments because $\tr(F_A)$ 
is smooth on $M$. 

Now we have
\begin{equation}
\label{eq:5.1.9}
\int_M \tr(F_{A^t}\wedge F_{A^t})\wedge \Omega = 
\int_M \tr(F_{A}\wedge F_{A})\wedge \Omega. \spn{5.1.9}
\end{equation}
This is the same as
\begin{equation}
\label{eq:5.1.10}
\int_M \left (F_{A^t}, T(F_{A^t}) \right ) dV_g
=\int_M\left (F_{A}, T(F_{A}) \right ) dV_g, \spn{5.1.10}
\end{equation}
where $T$ is the operator $-\ast\cdot \Omega\wedge$ acting on $2$-forms.
Then the $\Omega$-anti-self-duality of $A$ states $T(F_A) = F_A$.
It follows that 
\begin{eqnarray*}
 YM(A^t) 
&=&\frac{1}{4\pi^4}\int_M |F_{A^t}|^2 dV_g \\[2ex]
&=&\frac{1}{4\pi^2}\int_M \left
(F_{A^t},(Id - T)(F_{A^t}) \right ) dV_g
+ \frac{1}{4\pi^2} \int_M \left (F_{A^t}, T(F_{A^t}) \right ) dV_g  \\[2ex]
&=&\frac{1}{4\pi^2}\int_M \left
(F_{A^t},(Id - T)(F_{A^t}) \right ) dV_g - {\rm Ch}_2(A).
\end{eqnarray*}
Since $(Id - T) (F_A)=0$ and $T$ is symmetric,
the last integral above is at the order $t^2$.
Therefore, 
$$
\frac{d}{dt} YM(A^t)\big |_{t=0} = 0.
$$
This implies that $A$ is stationary.
\enddemo

In fact, we believe 
that {\it any admissible Yang-Mills connection {\rm (}\/possibly under certain mild conditions\/{\rm )} is stationary.}
If this is true, we can conclude from Corollary 4.5.2
that the blow-up locus of any Yang-Mills connections are stationary, 
in other words, it is a generalized minimal variety.

\demo{{\rm 5.2.}\ A removable singularity theorem}
In this section, we always assume that $A$ is an admissible
Yang-Mills connection on $M$ and stationary. Fix $p \in S=S(A)$,
where $S(A)$ denotes the singular set of $A$.
Let $r_p, c(p), a$ be as in Theorem 2.1.2. Our goal of this section 
is to prove a removable singularity theorem under appropriate assumptions.
We assume that $S\cap B_{\frac{r_p}{2}}(p)$ satisfies the following uniform covering 
(UC) property: for any $y \in S\cap B_{\frac{r_p}{2}}(p)$ and $\delta \le r <
\frac{ r_p}{2}$,
there are always balls $B_\delta(x_i)$ ($i=1,\cdots, l$) such that
$x_i \in S$, $S\cap B_{r}(y) \subset \bigcup _i B_\delta (x_i)$
and $l \delta ^{n-4} \le c r^{n-4}$ for some uniform constant $c>0$.
One can easily show that this (UC) property holds, if there is a measure $\mu$ with support
$S$ such that the  total measure $\mu (S\cap B_{r_p}(p)) < \infty$, and
for every $x\in S\cap  B_{r_p}(p)$, $r^{4-n}\mu(S\cap B_r(x))$ is decreasing with $r$, and the density 
$\Theta(x) = \lim_{r\to 0}r^{4-n}\mu(S\cap B_r(x))>0$.
In particular, if $A$ is the limit of smooth Yang-Mills connections $A_i$
outside $S$, then $S$ has the  (UC) property, 
since $\mu = \lim_{i\to\infty} |F_{A_i}|^2 dV_g$ satisfies the above conditions.
\specialnumber{5.2.1} \proclaim{Theorem}
\label{th:5.2.1}
Let $A${\rm ,} $S$ be as above{\rm .} Then there is an $\varepsilon >0${\rm , }
depending only on $n=\dim M${\rm ,} such
that for any $p\in S$ and $0<r< r_p${\rm ,} if
\begin{equation}
\label{eq:5.2.1}
r^{4-n}\int_{B_{r}(p)}\mid F_A\mid^2\,dV_g<\varepsilon, \spn{5.2.1}
\end{equation}
then there is a gauge transformation $\sigma$ near $p$ such that
$\sigma(A)$ extends to be a smooth connection near $p${\rm .}
\endproclaim

A direct corollary of this is the next result:

\specialnumber{5.2.2} \proclaim{Theorem}
\label{th:5.2.2}
Let $A${\rm ,} $S$ be as in Theorem {\rm 5.2.1.}
Then there is a gauge transformation $\sigma$ such that 
$\sigma(A)$ is smooth outside 
a closed subset $S'$ of $H^{n-4}$\/{\rm -}\/measure zero{\rm .}
\endproclaim
\demo{Proof}
Let $\varepsilon$ be given as in Theorem~\ref{th:5.2.1}.  
Then for any $x\in M$,
the limit 
\smallbreak
\centerline{${\displaystyle \lim_{r\to 0}r^{4-n}e^{ar^{2}}\int_{{B_{r}}(x)}\mid F_A\mid^2\, dV_g}$}
\noindent 
exists. Define
\begin{equation}
\label{eq:5.2.2.}
S'=\{x\in M\mid\lim_{r\to 0}r^{4-n}e^{ar^{2}}\int_{B_r(x)}\mid
F_A\mid^2\, dV_g \ge \varepsilon \}. \spn{5.2.2}
\end{equation}
Then by Theorem 5.1.1, $S'$ is closed.
Moreover, by using standard arguments as those in the proof of 
Lemma~\ref{lem:3.1.4} $(c)$, we can show that $H^{n-4}(S')=0$.

By Theorem~\ref{th:5.2.1}, there is a countable covering
$\{U_\alpha\}$ of $M\backslash S'$, so that for each $\alpha$,
there is a gauge transformation $\sigma_\alpha$ on $(M\backslash
S')\cap U_\alpha$  such that $E$ is trivial over $U_\alpha$ and
$D_{\sigma_{\alpha}(A)}=d+A_\alpha$ for some smooth $A_\alpha$.  
It follows that for any $\alpha, \beta$, we have the
transition function
$g_{\alpha\beta}=\sigma_\alpha\cdot\sigma_{\beta}^{-1}:
U_{\alpha}\cup U_\beta \backslash S'\to G$, where $G$ is the structure group
of $E$, such that
\begin{equation}
\label{eq:5.2.3}
A_\alpha=g_{\alpha\beta}^{-1}dg_{\alpha\beta}+g^{-1}_{\alpha\beta}A_\beta
g_{\alpha\beta}. \spn{5.2.3}
\end{equation}
Therefore, $g_{\alpha\beta}$ extends to   a smooth map on
$U_\alpha\cap U_\beta$, since $g_{\alpha\beta}$ takes values in a
compact group $G$. Furthermore, $\{g_{\alpha\beta}\}$ satisfies the
cocycle condition
$$
g_{\alpha\beta}\cdot g_{\beta\gamma}=g_{\alpha\gamma}\hbox{ on }
U_\alpha\cap U_\beta\cap U_\gamma.
$$
Therefore, $\{g_{\alpha\beta}\}$ defines a $G$-bundle $E^\prime$
over $M\backslash S$ extending $E\mid_{M\backslash S(A)}$, and
$\{A_\alpha\}$ defines a Yang-Mills connection for $E^\prime$.
The theorem is proved.
\enddemo

The rest of this section is devoted to the proof of
Theorem~\ref{th:5.2.1}. By scaling, we may assume that
$r=5$, $M=B_5(p)$ and $E$ is trivial over $M$.
For simplicity, we may further assume that the metric $g$
is flat. The general case can be proved by identical arguments.

We will always denote by $c$ a uniform constant. As before, we write
$S=S(A)$ as the singular set of $A$.
\specialnumber{5.2.3} \proclaim{Lemma}
\label{lem:5.2.3}
There is a gauge transformation $\sigma$ on $M\backslash S$ such
that for any $x\in B_3(p)\backslash S${\rm ,}
\begin{eqnarray}
\label{eq:5.2.4}
\rho(x)^{n-2}\mid A^\sigma\mid^2(x) &\le & 
c\int_{B_{\frac{1}{2}\rho(x)}(x)}\mid F_A\mid^2 dV_g,\spn{5.2.4}\\ 
\qquad \int_{B_{\frac{2}{5}\rho(x)}(x)} \left (
\frac{\mid A^\sigma\mid^2}{\rho(x)^2}
+ \mid\nabla A^\sigma\mid^2\right )\,dV_g&\le &
c\int_{B_{\frac{1}{2}\rho(x)}(x)}\mid F_A\mid^2 dV_g, \spn{5.2.5}
\end{eqnarray}
where $\rho(x)= d(x,S)$ and $D_{\sigma (A)}=d+A^\sigma$ with
$ A^\sigma\in\Omega^1\left(M\backslash S, {\rm Lie}(G)\right)${\rm .}
\endproclaim
\demo{Proof}
We may assume that $2^{n-4} e^a \varepsilon \le \varepsilon(n)$, where
$\varepsilon(n)$ is as  given in Theorem 2.2.1.
Then by the monotonicity (5.1.2), for any $x\in B_3(p)\backslash S$,
$$
\rho(x)^{4-n} \int_{B_ {\rho(x)}(x)} |F_A|^2 dV_g < \varepsilon (n).
$$
It follows from Uhlenbeck's curvature estimate (Theorem 2.2.1)
that 
$$
|F_A|(y) \le \frac{c\sqrt{\varepsilon(n)}}{\rho(x)^2},~~{\rm for}~
y\in B_{\frac{1}{2}\rho(x)}(x).
$$
Note that $c$ always denotes a uniform constant in this proof.

Next, using Theorem 1.2.7 in [Uh1, p.~18], we can have
a gauge transformation
$\sigma_x$ over $ B_{\frac{1}{2}\rho(x)}(x)$, such that
$D_{\sigma_x(A)} = d + A^{\sigma_x}$ and for any 
$y \in B_{\frac{1}{40}\rho(x)}(x)$,
$$
\rho(x) |A^{\sigma_x}|(y) + \rho(x)^2 |\nabla A^{\sigma_x}|(y) \le
c \left (\left (\frac{\rho(x)}{20}\right )^{4-n} \int_{B_ {\frac{\rho(x)}{20}}(x)} 
|F_A|^2 dV_g \right )^{\frac{1}{2}}.
$$
Now we outline the construction of $\sigma$ from those $\sigma_x$.
We cover $M\backslash S$ by balls $B_{r_i}(x_i)$ satisfying: (1) $x_i\in M\backslash S$
and $r_i = \frac{1}{40} \rho(x_i)$; (2) For any $x\in M\backslash S$, the number of 
those $B_{r_i}(x_i)$ containing $x$ is uniformly finite. For each $i$, denote by $\sigma_i$ the above
$\sigma_{x_i}$. If $B_{r_i}(x_i)\cap B_{r_j}(x_j)$ is nonempty, then
$r_i \le 2 r_j$ and $r_j \le 2r_i$; so by the above estimate for $A^{\sigma_i}$ and 
$A^{\sigma_j}$, we can obtain
$$
r_i |d \sigma_i\cdot \sigma_j^{-1}| + r_i^2 |\nabla d \sigma_i\cdot \sigma_j^{-1}| \le
c \left (r_i^{4-n}\int_{B_{r_i}(x_i)\cup B_{r_j}(x_j)} |F_A|^2 dV_g \right )^{\frac{1}{2}},
$$
on the overlap $B_{r_i}(x_i)\cap B_{r_j}(x_j)$.
Notice that $B_{r_i}(x_i)\subset B_{\rho(x)/10}(x)$ whenever $x\in B_{r_i}(x_{i})$.
Thus we can glue these $\sigma_i$ to get a gauge transformation
$\sigma$ such that 
$$
\rho(x) |A^{\sigma}|(x) + \rho(x)^2 |\nabla A^{\sigma}|(x) \le
c \left (\left (\frac{\rho(x)}{10}\right )^{4-n} \int_{B_ {\frac{\rho(x)}{2}}(x)} 
|F_A|^2 dV_g \right )^{\frac{1}{2}}.
$$
\vglue4pt\noindent This implies that for any $y\in B_{\frac{2}{5}\rho(x)}(x)$, we have
$$\rho(x) |A^{\sigma}|(y) + \rho(x)^2 |\nabla A^{\sigma}|(y) \le
c \left (\left (\frac{\rho(x)}{2}\right )^{4-n} \int_{B_ {\frac{\rho(x)}{2}}(x)} 
|F_A|^2 dV_g \right )^{\frac{1}{2}}.
$$
\vglue4pt\noindent Then the lemma follows easily.
\enddemo 

For simplicity, we assume that $\sigma$ can be taken to be ${\rm Id}$
in the above lemma.
\specialnumber{5.2.4} \proclaim{Lemma}
\label{lem:5.2.4}
Let $A$ be as above{\rm .} Then
\begin{equation}
\label{eq:5.2.6}
\int_{B_{1}(x)}\left (\frac{\mid A\mid^2}{\rho(y)^2}
 + \mid \nabla A\mid^2\right ) \,
dV_g\leq c \int_{B_{3}(x)}\mid F_A\mid^2\, dV_g, \spn{5.2.6}
\end{equation}
where $\rho(y)=d(y,S)${\rm .}
\endproclaim
\demo{Proof} Since $\rho$ is Lipschitz and $\mid\nabla\rho\mid\equiv
1$, by the co-area formula (cf.\ [Si2]), we have
\begin{equation}
\label{eq:5.2.7}
\int_{B_1(x)}\frac{\mid A\mid^2}{\rho(y)^2}\, dV_g=
\int^{1}_{0}\frac{dr}{r^2}\int_{\rho^{-1}(r)\cap B_{1}(x)}
\mid A\mid^2 d \mathcal H^{n-1} \spn{5.2.7}
\end{equation} 
where $d\mathcal H^{n-1}$ denotes the induced measure on the level
surface $\rho^{-1}(r)$.

For any $r\le 1$, there is a covering
$\{B_{\frac{2r}{5}}(x_{ir})\}_{1\leq i\le N_{r}}$ of 
$\rho^{-1}([\frac{2}{3}r, \frac{4}{3}r])\break\cap B_1(x)$, such that
$\rho(x_{ir})=r$ and for any $y\in\rho^{-1}([\frac{2}{3}r,
\frac{4}{3}r])$, the number of balls $B_{\frac{r}{2}}(x_{ir})$
containing $y$ is uniformly bounded.  Hence,
\begin{eqnarray*}
 &&\hskip-1.5in \int^{1}_{0}\frac{dr}{r^2}\int_{B_{1}(x)\cap\rho^{-1}(r)}
\mid A \mid^2 d\mathcal H^{n-1} \\[1ex]
&=& \frac{1}{ln 2}\int^{1}_{0}\frac{dr}{r^2}
      \int^{\frac{3}{2}r}_{\frac{3}{4}r}\frac{ds}{s}\int_{B_{1}(x)
      \cap\rho^{-1}(r)}\mid A\mid^2 d\mathcal H^{n-1} \\[1ex]
&=&\frac{1}{ln 2}\int^{\frac{3}{2}}_{0}\frac{ds}{s}
        \int_{B_{1}(x)\cap\rho^{-1}([\frac{2}{3}s,\frac{4}{3}s])}
\frac{1}{\rho^2}\mid A\mid^2\, dV_g\\[1ex]
&\le& \frac{1}{ln 2}\int^{\frac{3}{2}}_{0}\frac{ds}{s}
       \left( \sum_i\int_{B_{\frac{2s}{5}}(x_{is})}
\frac{\mid A\mid^2 }{\rho^2(y)}
       \, dV_g \right)\\[1ex]
\hbox{By (5.2.5)}&\le& c\int^{\frac{3}{2}}_{0}\frac{ds}{s}\left(
          \sum_i\int_{B_{\frac{s}{2}}(x_{is})}\mid F_A \mid ^2\, dV_g \right) \\
&\le& c\int^{\frac{3}{2}}_{0}\frac{ds}{s}\int_{B_3(x)
          \cap\rho^{-1}([\frac{s}{2},\frac{3s}{2}])}\mid
        F_A\mid^2\, dV_g\\[1ex]
&\le& c\int^{3}_{0}dr\int^{2r}_{\frac{2}{3}r}\frac{ds}{s}
          \int_{B_{3}(x)\cap\rho^{-1}(r)}\mid F_A\mid^2\, dV_g\\[1ex]
&\le& c\int_{B_3(x)}\mid F_A\mid^2\, dV_g.
\end{eqnarray*}
Similarly, we can derive
$$
\int_{B_{1}(x)}\mid\nabla A\mid^2\, dV_g\leqslant
c\int_{B_{3}(x)} \mid F_A\mid^2\, dV_g.
$$
The lemma is proved.
\enddemo

\specialnumber{5.2.5} \proclaim{Lemma}
\label{lem:5.2.5}
Let $A$ be as above{\rm .}  Then there are a function $\alpha$ and a
$2$\/{\rm -}\/form $\beta$ such that
\begin{equation}
\label{eq:5.2.9}
 A~=~d\alpha +d^\ast \beta,~~ d\beta =0 \quad\hbox{ on }\ B_{1}(x),\spn{5.2.8}
\end{equation}
\begin{equation}
\label{eq:5.2.10}
 \mid\mid \alpha \mid\mid_{H^{1,2}(B_{1}(x))}+
\mid\mid \beta \mid\mid_{H^{1,2}(B_{1}(x))}  \le  
c\mid\mid A\mid\mid_{L^{2}(B_{2}(x))} .\hskip.75in \spn{5.2.9}
\end{equation}
 \endproclaim

\demo{Proof}
Let $\eta:B_3(x)\to\RR^1$ be a cut-off function: $\eta(y)=1$ 
for $d(x,y)\le 1$, $\eta(y)=0$ for $d(x,y)\ge
2$ and $\mid\nabla \eta\mid\le 1$. 
By the extension of the classical Hodge-de Rham decomposition due
to Iwaniec and Martin, we have unique $\alpha$ and $\beta$ on $\RR^n$ such
that $\eta A=d\alpha +d^\ast \beta $ on $\RR^n$, $d\beta=0$,
and
$$
\mid\mid \alpha \mid\mid_{H^{1,2}(\RR^n)} + 
\mid\mid \beta \mid\mid_{H^{1,2}(\RR^n)}\ \le\
c \mid\mid \eta A \mid\mid_{H^{1,2}(R^{n})}.
$$
Then the lemma follows easily.
\enddemo 

Put $\tilde{A}= A-d\alpha$; then $d^\ast \tilde{A}=0$. 
Since $A$ is a Yang-Mills connection in the weak sense, 
\begin{eqnarray}
\label{eq:5.2.13}
0 &=& D^{\ast}_{A}F_A\spn{5.2.10}\\ [2ex]
  &=& d^\ast F_A +[F_A, A]\nonumber\\ [2ex]
  &=& d^\ast d{A} + d^\ast (A\wedge A) +[F_A,A]\nonumber\\ [2ex]
  &=& d^\ast d\tilde{A} + d^\ast (A\wedge A) +[F_A,A]\nonumber\\ [2ex]
  &=& (d^\ast d +d d^\ast)\tilde{A} + d^\ast(A\wedge A) + [F_A,A]. \nonumber
\end{eqnarray}
We decompose
\begin{equation}
\label{eq:5.2.14}
\tilde{A}=\tilde{A}_{0}+ \tilde{A}_{1}, \spn{5.2.11}
\end{equation}
such that
\begin{equation}
\label{eq:5.2.15}
(d^\ast d + dd^\ast)\tilde{A}_{0}=0,~~~ \hbox{ in }B_{1}(x), \spn{5.2.12}
\end{equation}
and
\begin{eqnarray}
\label{eq:5.2.16}
(d^\ast d+dd^\ast)\tilde{A}_{1}&=&
-[F_A,A]- d^\ast(A\wedge A),~~\hbox{ in }B_{1}(x)\spn{5.2.13}\\[2ex]
\tilde{A}_{1}&=&0 ~~~~~~\hbox{ on }\partial B_{1}(x).\nonumber
\end{eqnarray}
\specialnumber{5.2.6} \proclaim{Lemma}
\label{lem:5.2.6} There exists
\begin{equation}
\label{eq:5.2.17}
\mid\mid
\tilde{A}_{1}\mid\mid_{H^{1,2}(B_{1}(x))}\ \le\
c\sqrt{\varepsilon}\mid\mid F_A\mid\mid_{L^{2}(B_{3}(x))}, \spn{5.2.14}
\end{equation}
where $\varepsilon$ is as given in {\rm (5.2.1).}
\endproclaim
\demo{Proof}
First we have from (5.2.4),
\begin{equation}
\label{eq:5.2.18}
\mid A \mid (y)\le \frac{c\sqrt{\varepsilon}}{\rho(y)},
~~~\forall y\in B_{3}(x)\backslash S. \spn{5.2.15}
\end{equation}
Multiplying (\ref{eq:5.2.16}) by $\tilde{A}_{1}$ and integrating
by parts, we obtain
\begin{eqnarray*}
&&\hskip-.5in \int_{B_{1}(x)}\mid \nabla \tilde{A}_{1} \mid^2 \, dV_g\nonumber\\[2ex]
&&\quad =\ -\int_{B_{1}(x)}\left((\tilde{A}_{1},[F_A,A])
     + (\tilde{A}_{1}, d^\ast (A\wedge A))\right)\, dV_g \nonumber\\[2ex]
&&\quad =\  -\int_{B_{1}(x)}\left((\tilde{A}_{1}, [F_A,A])
     +( d\tilde{A}_{1}, A\wedge A)\right)\, dV_g \nonumber\\[2ex]
&&\quad\le\ c\sqrt{\varepsilon}\left(\int_{B_{1}(x)}
\frac{\mid\tilde{A}_{1}\mid\mid F_A\mid}{\rho(y)}\, dV_g +
\int_{B_{1}(x)}\frac{\mid
d\tilde{A}_{1}\mid\mid A\mid}{\rho(y)}\, dV_g \right)\nonumber\\[2ex]
\noalign{\hbox{by (5.2.6)}}\\
&&\quad\le\  c\sqrt{\varepsilon}\left(\int_{B_{3}(x)}\mid
F_A\mid^2 \, dV_g \right)^{\frac{1}{2}}
\left(\int_{B_{1}(x)}\left(\frac{\mid\tilde{A}_{1}\mid^2}
{\rho(y)^2} + \mid \nabla \tilde{A}_{1}\mid^2 \right) \, dV_g \right)^{\frac{1}{2}}   .  
\end{eqnarray*}
Then (\ref{eq:5.2.17}) follows from the next lemma.
\enddemo 
\specialnumber{5.2.7} \proclaim{Lemma}
\label{lem:5.2.7}
For any function $f$ vanishing on $\partial B_{1}(x)${\rm ,} 
\begin{equation}
\label{eq:5.2.19}
\int_{B_{1}(x)}\frac{\mid f\mid^2}{\rho(y)^2}\, dV_g
\le c \int_{B_{1}(x)}\mid \nabla f\mid^2\, dV_g.\spn{5.2.16}
\end{equation}
\endproclaim
\demo{Proof}
This lemma follows directly from a result of C. Fefferman and D. Phong [FP]
(also see [CW, Th.~1.4], [CWW], [Fef]),
once we verify the following: for any $y \in B_1(x)$ and $r\le 1$,
\begin{equation}
\label{eq:5.2.20}
\int_{B_{r}(y)} \frac{1}{\rho^3} dV_g \le c r^{n-3}  \spn{5.2.17}
\end{equation}
where $c$ is a uniform constant. 

Let us check (\ref{eq:5.2.20}). If $\rho (y) \ge 2r$,
then $r\le \rho(z) \le 3r $ for any $z\in B_{r}(y)$.  Now
$$
\int_{B_{r}(y)} \frac{1}{\rho^3} dV_g \le \frac{c}{r^3} \int_{B_r(y)} dV_g
\le c r^{n-3}.
$$

Next, we assume that $\rho (y) \le 2r$. By our assumption on $S$,
for any $\delta < 4r$, there are $L_\delta$ balls $B_\delta(x_i)$
such that $S\cap B_{r}(y)\subset \bigcup_i B_\delta(x_i)$
and $L_\delta \le c \left (\frac{r}{\delta}\right )^{n-4}$. Then
by the co-area formula,
\begin{eqnarray*}
 \int_{B_{r}(y)} \frac{1}{\rho^3} dV_g  
&=&\int_0^{4r}
\frac{1}{s^3} ds \int_{\rho^{-1}(s)\cap B_{r}(y)} dH^{n-1}\\[2ex]
&=&(4r)^{-3}\int_{B_{r}(y)} dV_g + 
3\int_0^{4r}
\frac{1}{s^4} ds \int_{\rho^{-1}([0,s])\cap B_{r}(y)} dV_g\\[2ex]
&\le&c r^{n-3} + 3 \int_0^{4r}\frac{1}{s^4} ds \left (
\sum_{i=1}^{L_s} \int_{B_{2s}(x_i)}dV_g\right ) \\[2ex]
&\le& c r^{n-3} + c r^{n-4} \int_0^{5r} ds\\[2ex]
&\le& cr^{n-3}.
\end{eqnarray*}
Thus (\ref{eq:5.2.20}) follows.
\enddemo

Let $\theta\in(0, 1)$ be fixed. 
Since $\tilde{A}_{0}$ is harmonic, we have, from
standard elliptic estimates, that
%
\begin{eqnarray}
\label{eq:5.2.22}
\frac{1}{\theta^{n-4}}\int_{B_{\theta}(x)}\mid
d\tilde{A}_{0}\mid^2
& \le& \theta^4\int_{B_{1}(x)}\mid d\tilde{A}_{0}\mid^2
dV_g \spn{5.2.18}\\
&\le&  \theta^4\int_{B_{1}(x)}\mid d\tilde{A}\mid^2 dV_g.\quad
\nonumber\end{eqnarray}
Then 
\begin{eqnarray}
\noalign{\vskip-3pt}
&& \spn{5.2.19}\label{eq:5.2.23}
\\
\noalign{\vskip-3pt}
  && \theta^{4-n}\int_{B_{\theta}(x)}\mid F_A\mid^2 \, dV_g\nonumber \\ 
 &=&  \theta^{4-n}\int_{B_{\theta}(x)}\left(\mid
  dA\mid^2 + 2(F_A, A\wedge A) - \mid A\wedge 
  A\mid^2  \right) \, dV_g  \nonumber \\[2ex]
& \le& \theta^{4-n}\int_{B_{\theta}(x)}\left( \mid
  d\tilde{A}\mid^2 + 2(F_A, A\wedge A) \right) \, dV_g \nonumber\\[2ex]
\hbox{by (5.2.4), (5.2.1)}
&\le& \theta^{4-n}\int_{B_{\theta}(x)}\left(\mid d 
  \tilde{A}\mid^2 + \frac {c \sqrt{\varepsilon} |A||F_A|}{\rho (y)}
 \right )\,dV_g 
\nonumber\\[2ex]
\hbox{by (5.2.6)}&\le& 
\theta^{4-n} \int_{B_{\theta }(x)}\mid d\tilde{A}\mid^2 \, dV_g
   + c \sqrt{\varepsilon}\theta^{4-n}\int_{B_{3}(x)}\mid
   F_A\mid^2 \, dV_g. \nonumber
\end{eqnarray} 
Similarly, we have
\begin{equation}
\label{eq:5.2.24}
 \int_{B_{1}(x)}\mid d\tilde{A}_{0}\mid^2
dV_g \le \int_{B_{1}(x)}\mid d\tilde{A}\mid^2 dV_g
+c \sqrt{\varepsilon} \int_{B_{3}(x)}\mid
F_A\mid^2 \, dV_g. \qquad\spn{5.2.20}
\end{equation} 
\pagebreak
On the other hand, using Lemma 5.2.6 and Lemma 5.2.4, we deduce
\begin{eqnarray}
\label{eq:5.2.25}
&&\spn{5.2.21}\\
 \int_{B_{\theta}(x)}\mid d\tilde{A}\mid^2 dV_g  
&=&\int_{B_{\theta}(x)}\left (\mid d\tilde{A}_0\mid^2 
+ |d\tilde A_1|^2 + 2 ( d\tilde A_0, d\tilde A_1)
\right ) dV_g\nonumber \\[2ex]
&\le&\int_{B_{\theta}(x)}\mid d\tilde{A}_0\mid^2 dV_g
+ 2||\tilde A_1||_{H^{1,2}(B_1(x))}
\nonumber\\[2ex]
&&\cdot \ \left ( \int_{B_1(x)}
|d\tilde A_0|^2 dV_g \right )^{\frac{1}{2}} +  
||\tilde A_1||^2_{H^{1,2}(B_1(x))}\nonumber\\[2ex]
&\le&\int_{B_{\theta}(x)}\mid d\tilde{A}_0\mid^2 dV_g
+ c \sqrt{\varepsilon} \int_{B_{3}(x)}\mid
F_A\mid^2 \, dV_g.\nonumber
\end{eqnarray}
It follows from the above four inequalities that 
$$
\theta^{4-n}\int_{B_{\theta}(x)}\mid F_A\mid^2\, dV_g
\le \theta^4 \int_{B_{1}(x)}\mid F_A\mid^2 \, dV_g
+ c \sqrt{\varepsilon} \theta^{4-n} \int_{B_{3}(x)}\mid
F_A\mid^2 \, dV_g.
$$
By scaling, we obtain
that for $r\le 1$ and $y\in B_{1}(p)$,
\begin{eqnarray}
\label{eq:5.2.26}
&&\spn{5.2.22}\\
  (\theta r)^{4-n}\int_{B_{\theta r}(x)}\mid F_A\mid^2\, dV_g
&\le& \theta^4 r^{4-n} \int_{B_{r}(x)}\mid F_A\mid^2 \, dV_g\nonumber \\[2ex]
& &+\ c \sqrt{\varepsilon} \theta^{4-n} r^{4-n} \int_{B_{3r}(x)}\mid
F_A\mid^2 \, dV_g. \nonumber
\end{eqnarray} 
Then,  from the monotonicity for $A$,  we have
$$
(\lambda r)^{4-n}\int_{B_{\lambda r}(y)}\mid F_A\mid^2 \, dV_g
\le  \left (3^4 + c \sqrt{\varepsilon(r)} \lambda^{-n}
\right ) \lambda ^4  
r^{4-n} \int_{B_{r}(x)}\mid F_A\mid^2 \, dV_g,
$$
where $\lambda = \frac{\theta}{3} < \frac{1}{3}$ and 
$$
\varepsilon(r) = r^{4-n} \int_{B_{r}(x)}\mid F_A\mid^2 \, dV_g\le 
8 \varepsilon.
$$
A simple iteration yields
\begin{eqnarray}
\label{eq:5.2.27}
&&   (\lambda^k r)^{4-n}\int_{B_{\lambda ^k r}(y)}\mid F_A\mid^2 \, dV_g
 \spn{5.2.23}\\[2ex]
&&\quad \quad\le\ \prod_{i=0}^{k-1} \left (1 + 
c \sqrt{\varepsilon(\lambda ^{i} r)} \lambda^{-n}
\right ) (3 \lambda) ^{4k}  
r^{4-n} \int_{B_{r}(x)}\mid F_A\mid^2 \, dV_g,
\nonumber \end{eqnarray}
where $k \ge 1$.  

Choose $\lambda$ and $\varepsilon $ such that 
$6^4 \lambda < 1$ and $8 c \sqrt{\varepsilon} \lambda^{-n}< 1$. 
This implies that for any $i \le k-1$,
$$
(1 + c \sqrt{\varepsilon(\lambda ^{i} r)} \lambda^{-n})3^4\lambda < 1.
$$

For any $r \le 1$, we define $k$, $r_0 \in (\frac{1}{3}, 1]$ by 
$\lambda ^k r_0 = r$. Then
\begin{eqnarray*}
r^{4-n}\int_{B_{r}(y)}\mid F_A\mid^2 \, dV_g &\le&  \lambda^{3k}r_0^{4-n}
\int_{B_{r_0}(y)}\mid F_A\mid^2 \, dV_g\\[2ex]
&\le&  r^3 r_0^{-3} \int_{B_{1}(y)}\mid F_A\mid^2 \, dV_g\\[2ex]
&\le& c r^3.
\end{eqnarray*}
Now we replace $\varepsilon( r)$ in (\ref{eq:5.2.26}) by
$c r^3$ and obtain
$$
(\theta r)^{4-n}\int_{B_{\theta r}(x)}\mid F_A\mid^2\, dV_g
\le \theta^4 r^{4-n} \int_{B_{r}(x)}\mid F_A\mid^2 \, dV_g
+ c \theta^{4-n} r ^{\frac{9}{2}}.
$$
Choose $\theta = \frac{1}{2}$ and $c'$ such that 
$c (\frac{1}{2})^{4-n}
+c' (\frac{1}{2})^{\frac{9}{2}} \le c'( \frac{1}{2})^4$. Then
$$
(\frac{r}{2})^{4-n}\int_{B_{\frac{r}{2 }}(x)}\mid F_A\mid^2\, dV_g +
 c' \left (\frac{r}{2}\right)^{\frac{9}{2}}
\le (\frac{1}{2})^4 \left( r^{4-n} \int_{B_{r}(x)}\mid F_A\mid^2 \, dV_g
+ c' r ^{\frac{9}{2}}\right ).
$$
It follows from this and a simple iteration that
$$
 r^{4-n} \int_{B_{r}(x)}\mid F_A\mid^2 \, dV_g
\le c''  r ^ 4,
$$
where $c''$ is some uniform constant.

Therefore, the curvature $F_A$ is bounded in $B_1(p)$.
Using results in\break [Uh2], we can construct a gauge transformation
$\sigma$ such that $d^* A_\sigma =0$ and\break $||A_\sigma||_{C^1(B_{1}(p))}$ 
is bounded. Since $D_{\sigma (A)}^* F_{\sigma(A)} = 0$,
$A_\sigma$ is smooth, and consequently,
$\sigma(A)$ extends to a smooth connection near $p$.
Theorem 5.2.1 is proved.

\demo{{\rm 5.3.}\ Cone-like Yang-Mills connections}
In this section, we study the infinitesimal structure of 
stationary Yang-Mills connections at their singular points.
Let $A$ be a stationary Yang-Mills connection on $M$ with
$L^2$-bounded curvature $F_A$. 
It follows from Theorem 5.1.1 that for any $x \in S$,
the limit
$$
\lim_{r\to 0} r^{4-n}\int_{B_r(x)} |F_A|^2 dV_g
$$
exists. Therefore, we can define
$$
S([A]) = \{ x\in M ~|~\lim_{r\to 0} r^{4-n} \int_{B_r(x)} |F_A|^2 dV \ge \varepsilon\},
$$
where $\varepsilon$ is as  given in Theorem 5.2.1. Then $S(A)$ contains
$S([A])$. Denote by $S$ the set $S([A])$.
By Theorem 5.2.2, we have $H^{n-4}(S)=0$; moreover, there is a gauge transformation
$\sigma$ on $M\backslash S(A)$ such that $\sigma (A)$
extends to a smooth connection on $M\backslash S$. Without loss of generality,
we may assume that $S=S(A)$.
 
Now we explain why $S$ is expected to be of Hausdorff codimension of 
at least $5$.

To analyze $A$ near $x$, we scale the metric $g$ and $A$
as follows: for any $\lambda \in (0,1)$, define
$$
g_\lambda =\lambda^{-2} g,~~~ A_\lambda = \tau_\lambda^*{\rm exp}_x^* A,
$$
where $\tau_\lambda : T_xM \mapsto T_xM$ maps $v$ to $\lambda v$.
Clearly, $A_\lambda$ is a stationary Yang-Mills connection
with respect to $g_\lambda$. Moreover, for any $R>0$, 
\begin{equation}
\label{eq:5.3.1}
R^{4-n} \int_{B_R(x, g_\lambda)} |F_{A_\lambda}|^2 dV_{g_\lambda}
= (\lambda R)^{4-n}\int_{B_{\lambda R}(x)} |F_{A}|^2 dV_{g} \le c,\hskip.5in \spn{5.3.1}
\end{equation}
whenever $\lambda$ is sufficiently small. Here
and in the following, $c$ always denotes a uniform constant.

Then we can deduce the following from results in Section 3.1: for
any sequence $\{\lambda_i\}$ with $\lim_{i\to \infty} \lambda (i) =0$,
taking a subsequence and gauge transformations
if necessary, we may assume that $A_{\lambda(i)}$ converges to
a connection $A^c$ outside $S_c\subset T_xM$. Here, $A^c$
is Yang-Mills with respect to the flat metric $g_0$ on $T_xM= \RR ^n$
and $H^{n-4} (S^c \cap B_R(0, g_0)) < \infty$. Moreover, we may assume that
$|F_{A_{\lambda (i)}}|^2 dV_g $ converges weakly to 
$|F_{A^c}|^2dV_{g_0} + \Theta_c H^{n-4} \lfloor S_c$,
where $\Theta_c$ is a function with its support in $S_c$. \enddemo

\specialnumber{5.3.1} \proclaim{Lemma}
\label{lem:5.3.1}
With the above notation{\rm ,}  {\rm (1) }
$\frac{\partial }{\partial r} \Theta_c =0$\/{\rm ; } 
{\rm (2)} $a \cdot S_c = S_c${\rm ,} where $a\cdot S_c$
denotes the set of points $a z$ with $z\in S_c$\/{\rm ; (3)}
$\frac{\partial }{\partial r} \rfloor F_{A^c} = 0${\rm .}
\endproclaim
\demo{Proof} By Theorem 4.5.1, for any vector field
$X$ with compact support,  
\begin{eqnarray}
\label{eq:5.3.1}
&&\spn{5.3.2}\\
&&\hskip-24pt -\int _{S_c} \divv _{S_c} X \Theta_c \,dH^{n-4}\nonumber  \\[2ex]
&=& \int _{T_xM} 
\left ( |F_{A^c}|^2
\divv X - 4\sum_{i,j=1}^n 
(F_{A^c}( \frac{\partial X}{\partial x_i}, 
\frac{\partial }{\partial x_j}), F_A(\frac{\partial }{\partial x_i}, 
\frac{\partial }{\partial x_j}))\right ) dV_{g_0},\nonumber
\end{eqnarray} 
where $x_1, \cdots, x_n$ are euclidean coordinates of $T_xM=\RR ^n$. 
Choosing $X(x) = \xi(r) r \frac{\partial }{\partial r}$, where
$r = \sqrt{\sum_i x_i^2}$ and $\xi$ has compact support, we obtain
\begin{eqnarray}
\label{eq:5.3.2}
&&  \int _{S_c} \left (\xi' r + (n-4)\xi\right ) \Theta_c \,dH^{n-4} + 
\int _{T_xM} 
\left (\xi' r + (n-4) \xi \right ) |F_{A^c}|^2 dV_{g_0} \spn{5.3.3}\\[2ex]
&&\hskip.5in = \int _{S_c} \xi ' r |\nabla^{\perp } r|^2 \Theta_c \,dH^{n-4}
 + 4\int_{T_xM} \xi' r \left |\frac{\partial }{\partial r}\rfloor 
F_{A^c} \right |^2 dV_{g_0}. \nonumber
\end{eqnarray} 
Following the arguments in deriving (2.1.20) from (2.1.15), we can deduce
from (5.3.3) that for any $\sigma < \rho$,\pagebreak
\begin{eqnarray}
\label{eq:5.3.3}
&&\int _{S_c\cap (B_\rho(0, g_0)\backslash B_\sigma(0, g_0))} 
r^{4-n} |\nabla^{\perp } r|^2 \Theta_c \,dH^{n-4} \spn{5.3.4}\\
&&\qquad\quad~+ ~
4 \int _{B_\rho(0, g_0)\backslash B_\sigma(0, g_0)} 
r^{4-n}\left |\frac{\partial }{\partial r}\rfloor 
F_{A^c} \right |^2 dV_{g_0}\nonumber \\[2ex]
&&\qquad =~\rho^{4-n}\left (\int _{S_c\cap B_\rho(0, g_0)} \Theta_c \,dH^{n-4}
 + \int_{B_\rho(0, g_0)} |F_{A^c}|^2 dV_{g_0} \right )\nonumber\\[2ex]
&& \qquad\quad~-~ \sigma^{4-n}\left (\int _{S_c\cap B_\sigma(0, g_0)} \Theta_c \,dH^{n-4}
 + \int_{B_\sigma (0, g_0)} |F_{A^c}|^2 dV_{g_0} \right ). \nonumber
\end{eqnarray} 
On the other hand, for any $s> 0$,  
\begin{eqnarray*}
&&\hskip-.5in s^{4-n}\left (\int_{B_s(0, g_0)}
|F_{A^c}|^2 dV_{g_0} + 
\int _{S_c\cap B_s(0, g_0)} \Theta_c \,dH^{n-4} \right )\\[2ex]
&&\qquad =~\lim_{i\to \infty} 
(\lambda(i) s)^{4-n}\int_{B_{\lambda(i) s}(x)} |F_A|^2 dV_g\\[2ex]
&&\qquad =~\lim_{s'\to 0}{s'}^{4-n} \int_{B_{s'}(x)}|F_A|^2 dV_g > 0.
\end{eqnarray*}
Therefore, 
\begin{eqnarray}
\label{eq:5.3.4}
&&\int _{S_c\cap (B_\rho(0, g_0)\backslash B_\sigma(0, g_0))} 
r^{4-n} |\nabla^{\perp } r|^2 \Theta_c \,dH^{n-4}\spn{5.3.5}\\
[2ex]
&&\qquad\qquad + 
4 \int _{B_\rho(0, g_0)\backslash B_\sigma(0, g_0)} 
r^{4-n}\left |\frac{\partial }{\partial r}\rfloor 
F_{A^c} \right |^2 dV_{g_0}= 0. \nonumber
\end{eqnarray}
This implies that $\nabla^{\perp } r = 0$ on $S_c$
and $ \frac{\partial }{\partial r}\rfloor 
F_{A^c} =0$; i.e., both (2) and (3) hold.  

Furthermore, arguing as we did in the proof of Lemma 3.2.1, we 
can deduce that $|F_{A^c}|^2dV_{g_0} + \Theta_c H^{n-4}
\lfloor S_c$ is a cone measure. Now, (1) holds. 
\enddemo
 
We will call such an $A^c$ a tangent Yang-Mills connection of $A$ at $x$.
In general, $A$ may have  a different tangent Yang-Mills connection at $x$,
which depends on choices of sequences $\{\lambda(i)\}$.
By Corollary 2.1.3, $A^c$ is gauge equivalent to $d+ B$ for
some $B: S^{n-1} \mapsto T^*S^{n-1}\otimes {\rm Lie}(G)$.
Thus $S(A^c)$ is invariant under radial
scaling and  so is $S([A^c])$. If $A^c$ is also stationary, 
then $H^{n-4}(S([A^c]))=0$.
Together with Uhlenbeck's
removable singularity theorem in [Uh1] (also see Theorem 5.2.1),
this implies that $S([A^c])=\{0\}$ whenever $n=5$.
If the blow-up set $S_c$ is empty, we further deduce 
that $A$ has an isolated singularity at $x$. This leads us to 

\nonumproclaim{{C}onjecture 1} If $A$ is stationary{\rm ,} then 
the Hausdorff codimension of $S([A])$ is at least $5${\rm .}
\endproclaim

We can expect stronger conclusion for $\Omega$-anti-self-dual
instantons. Now let $A$ be an $\Omega$-anti-self-dual instanton.
Then its tangent Yang-Mills connection $A^c$ is\break $\Omega_x$-anti-self-dual.
Moreover, $\Omega_x$ is a nonvanishing constant form.
\specialnumber{5.3.2} \proclaim{Lemma}
\label{lem:5.3.3}
If $v \rfloor F_{A^c} = 0$ for any $v$ in
a subspace $L \subset T_xM$ of dimension $n-5${\rm ,} then modulo
gauge transformations{\rm ,} $A^c$ extends smoothly to a connection on $T_xM${\rm .}
\endproclaim

\demo{Proof} Write $\Omega_x = dV_L \wedge d \ell\ + \ \Omega'_x$,
such that $\Omega'_x$ is perpendicular to any form $dV_L \wedge \varphi$,
where $dV_L$ is a volume form on $L$ and $\varphi$ is a $1$-form.
Then $\ell \not= 0$; otherwise, $F_{A^c}=0$ by our assumption and
$\Omega_x$-anti-self-duality. Furthermore, we have 
$$-\ast (F_{A^c}\wedge dV_L \wedge d\ell) = F_{A^c}.$$
Hence, $\frac{\partial }{\partial \ell} \rfloor F_{A^c} =0$
and modulo a gauge transformation, $A^c$ is the pull-back of 
an anti-self-dual connection on the $4$-subspace perpendicular to $L$ and 
$\frac{\partial }{\partial \ell}$. By the removable singularity theorem of Uhlenbeck
in dimension 4 (also see Theorem 5.2.1), there is a gauge transformation
$\sigma$ such that $\sigma(A^c)$ extends to a smooth connection on $T_xM$.
Hence, the lemma is proved.
\enddemo

\nonumproclaim{{C}onjecture 2}  If $A$ is $\Omega$\/{\rm -}\/anti\/{\rm -}\/self\/{\rm -}\/dual{\rm ,}
then its singular set $S([A])$ has Hausdorff codimension at least $6${\rm .}
\endproclaim

Both conjectures can be affirmed if one can show that
$\lim_{i\to \infty} S(A_i) \subset S(A)$ for any sequence of Yang-Mills 
connections $A_i$ converging to $A$.

Finally, let us discuss briefly   the classification of
tangent $\Omega$-anti-self-dual instantons 
on $\RR ^6$ with the only singularity at $0$. Since $\Omega$ is
a linear $2$-form on $\RR ^6$, we may choose coordinates
$x_1, \cdots, x_6$, such that 
$$
\Omega = a_1 dx_1 \wedge dx_2 + a_2 dx_3\wedge dx_4 + a_3 dx_5 \wedge dx_6.
$$
Let $A^c$ be a nonflat tangent $\Omega$-anti-self-dual instanton. Then 
$$-\ast_5 (\alpha \wedge F_{A^c}) = F_{A^c}$$ 
on $S^5 \subset \RR ^6$, where $\alpha = \frac{\partial }{\partial r} 
\rfloor \Omega$
and $\ast_5$ is the Hodge operator on $S^5$. Since $F_A\not= 0$, 
$\alpha \not=0$. A simple computation shows that
$|\alpha|(x)$ has to be $1$ for any $x\in S^5$. Hence, we may assume
that $a_1=a_2=a_3=1$, and consequently, $\Omega$ is the standard
symplectic form on $\RR^6$. If $J_0$ denotes the
complex structure on $\RR ^6$ such that the  $dx_{2i-1} + \sqrt{-1}dx_{2i}$
($i=1,2,3$) span the induced holomorphic tangent bundle, then
$A$ is Hermitian-Yang-Mills with respect to this complex structure.
Moreover, we have $v \rfloor F_{A^c}=0$ when $v$ is either 
$\frac{\partial}{\partial r}$ or $J_0 (\frac{\partial}{\partial r})$.
This implies that modulo a gauge transformation, $A^c$ is the pull-back
of a Hermitian-Yang-Mills connection on $\CC P^2$.
Conversely, any Hermitian-Yang-Mills connections on $\CC P^2$
give rise to a tangent $\Omega$-anti-self-dual instanton on~$\RR ^6$.

If $\Omega$ is the 3-form of Section 1.4, defining the $G_2$-structure
on $\RR^7$, then\break tangent $\Omega$-anti-self-dual instantons 
are in one-to-one correspondence with\break
Hermitian-Yang-Mills connections on $S^6$ with respect to the almost 
complex structure induced by $\Omega$. These are all the possible
tangent $\Omega$-asd (anti-self-dual) instantons on $\RR ^7$. 

It is also possible to classify all tangent $\Omega$-anti-self-dual 
instantons on~$\RR ^8$. Then one problem is how to show that 
singularities of any $\Omega$-anti-self-dual instantons are modeled on
these tangent connections on a manifold of dimension no more than $8$.

\section{Compactification of moduli spaces}
\label{chap:6}

In this chapter, we first construct a compactification of 
the moduli space of anti-self-dual instantons.
Then we discuss briefly possible extensions of results proved in the  last few chapters.

\demo{{\rm 6.1.}\ Compactifying moduli spaces}
Let $(M,g)$ be a compact Riemannian $n$-manifold and
$\Omega$ be a closed differential form of degree $n-4$. Let
$E$ be a unitary vector bundle over $M$.
Recall that ${\cal M}_{\Omega, E}$ consists of all equivalence classes 
of $\Omega$-anti-self-dual, often abbreviated as $\Omega$-asd, 
instantons on $M$, i.e., solutions of
(1.2.2). Here, two solutions $A_1$ and $A_2$ are equivalent if and only if
there is a gauge transformation $\sigma$ of $E$ such that 
$\sigma(A_1) = A_2$. In general, ${\cal M}_{\Omega, E}$ may not be
compact.

We now describe in detail the compactification outlined in the introduction.
A generalized $\Omega$-asd instanton
is made of (1) an admissible $\Omega$-asd instanton 
$A$ of $E$, which extends to become a smooth connection over
$M\backslash S(A)$ for a closed subset $S(A)$ with $(n-4)$-dimensional
Hausdorff measure $H^{n-4}(S(A))\break = 0$; (2) a closed integral current
$C=(S, \Theta)$ calibrated by $\Omega$ satisfying the
energy identity
$$\frac{1}{4\pi^2} 
\int_M |F_A|^2 dV_g + \int _S\Theta dH^{n-4} = \int_M
{\rm Ch}_2(E) \wedge \Omega,
$$
where ${\rm Ch}_2(E)$ denotes the second Chern character of $E$. 

If the co-norm $|\Omega| \le 1$, $C$ is an area-minimizing integral
current, so that it follows from [Am] that $C$ can be represented 
by $\sum _a m_a C_a$ satisfying: $m_a = \Theta |_{C_a}$ and
each $C_a$ is closed and of the form $C_a^0 \bigcup {\rm Sing}(C_a)$
such that $C_a^0$ is a smooth submanifold calibrated by $\Omega$
and ${\rm Sing}(C_a)$ is a closed subset of Hausdorff codimension at least two.
\enddemo

\demo{{R}emark {\rm 12}}
We believe that the singularity of each 
$\Omega$-calibrated cycle is also of a certain geometric structure. 
If $\Omega$-calibrated cycles $C_a$ are holomorphic, then each singular set
${\rm Sing}(C_a)$ is a holomorphic subvariety.
\enddemo

Two generalized $\Omega$-asd instantons
$(A, C)$, $(A', C')$ are equivalent if and only if
$C=C'$ and there is a gauge transformation
$\sigma$ on $M\backslash S(A) \cup S(A')$,
such that $\sigma(A)= A'$ on $M\backslash S(A) \cup S(A')$. 
We denote by $[A, C]$ the equivalence class represented by
$(A,C)$. We identify $[A, 0]$ with $[A]$ in $ {\cal M}_{\Omega, E}$
if $A$ extends to a smooth connection of $E$ over $M$ modulo
a gauge transformation. 

We define $\overline {\cal M}_{\Omega, E}$
to be set of all equivalence classes of generalized $\Omega$-anti-self-dual
instantons of $E$.

\demo{{R}emark {\rm 13}}
A natural problem occurs when an $\Omega$-calibrated cycle,
or simply a submanifold, is actually the limit of a sequence of 
$\Omega$-asd instantons. More generally, one may ask if a minimal submanifold
of dimension $n-4$ can be the limit of a sequence of Yang-Mills connections.
It is a delicate and interesting problem involving use of the 
implicit function theorem.
\enddemo

One can define the first two Chern forms of $(A,C)$ as follows:
${\rm Ch}_1(A,C)= {\rm Ch}_1(A)$ is given by $\frac{\sqrt{-1}}{2\pi}
\tr(F_A)$ and ${\rm Ch}_2(A,C)= {\rm Ch}_2(A) + {\rm PD}(C)$, where ${\rm Ch}_2(A)$ is
given by $-\frac{1}{4\pi^2} \tr(F_A\wedge F_A)$
and ${\rm PD}(C)$ denotes the Poincar\'e dual of the integral current $C$. 
Since $H^{n-4}(S(A))=0$ for a generalized $\Omega$-asd $(A,C)$,
both ${\rm Ch}_1(A)$ and ${\rm Ch}_2(A)$ are closed currents on $M$.
So they give rise to cohomology classes of $M$.
In fact, ${\rm Ch}_2(A, C)$ always represents ${\rm Ch}_2(E)$ in $H^*(M, \ZZ )$ for
any $[A,C]$ in $\overline {\cal M}_{\Omega, E}$.

The topology of $\overline {\cal M}_{\Omega, E}$ can be defined as 
follows: a sequence $[A_i, C_i]$ converges to $[A,C]$ in 
$\overline {\cal M}_{\Omega, E}$ if and only if (1) 
$C$ can be decomposed into two closed integral currents
$C'+C''$ such that $C_i$ converges to $C'$ in $M$ with respect to 
the standard topology for currents; (2) 
There are gauge transformations $\sigma_i$ such that $\sigma_i(A_i)$
converges to $A$ outside $S(A)$ and the support of $C$, and
the generalized Chern forms
${\rm Ch}_2(\sigma_i(A_i), C_i)$ converge to ${\rm Ch}_2(A,C)$ as currents.
One can show that $\overline {\cal M}_{\Omega, E}$
is then a Hausdorff topological space which follows from results in Chapter 4 and 5.

\specialnumber{6.1.1} \proclaim{Theorem}
\label{th:6.1.1}
For any $M${\rm ,} $g${\rm ,} $\Omega$ and $E$ as above{\rm ,}
$\overline {\cal M}_{\Omega, E}$ 
is compact with respect to this topology{\rm .}
\endproclaim 

Let $T$ be a compact family of metrics and closed $(n-4)$-forms
$g_t, \Omega_t$. Then by the above arguments, we can show:
\specialnumber{6.1.2} \proclaim{{C}orollary}
\label{cor:6.1.2}
For any $M${\rm ,} $T=\{g_t,\Omega_t\}$ and $E$ as above{\rm ,}
$\bigcup_{t\in T}\overline {\cal M}_{\Omega_t, E}$ 
is compact with respect to the topology defined above{\rm .}
\endproclaim

We end this section with a generalization of Theorem 6.1.1.
in the case where $\Omega$ is not necessarily closed.

We will still call $A$ an $\Omega$-asd instanton whenever
$A$ satisfies the equation in Lemma 1.2.1, even if $\Omega$ is not closed.
Suppose now that $\Omega$ has the decomposition $\Omega_1 + \Omega_2$,
such that $\Omega_1$ is closed and for any $2$-form $\varphi $,
\begin{equation}
\label{eq:6.1.1}
-\varphi\wedge \varphi \wedge \Omega_2 < |\varphi|^2 dV_g. \spn{6.1.1}
\end{equation}
Then for any $\Omega$-asd instanton $A$, we still have an {\it a priori}
bound on $YM(A)$ as we did in (1.2.3). Following the arguments
in the proof of Theorem 6.1.1, we can obtain the next result:

\specialnumber{6.1.3} \proclaim{Theorem}
\label{th:6.1.3}
Let $\Omega=\Omega_1+\Omega_2$ be as above{\rm .} Then 
$\overline {\cal M}_{\Omega, E}$ is compact{\rm .}
\endproclaim 

Note that an $\Omega$-asd instanton may not be Yang-Mills if $\Omega$
is not closed. 

\demo{{\rm 6.2.} Final remarks}
We expect that one can define certain deformation invariants by using
$\overline {\cal M}_{\Omega, E}$ as one did in the case of Donaldson, 
Gromov-Witten and Seiberg-Witten invariants, etc. 

More precisely, let $(M, g)$ be a compact Riemannian manifold, and
$\Omega$ be a degree $n-4$ form satisfying (\ref{eq:6.1.1}) and
the ellipticity condition: for any $x$ in $M$, the symmetric
operator $T=-\ast \Omega \wedge$ on $2$-forms has 1 as its eigenvalue,
of multiplicity exactly equal to $\frac{(n-1)(n-2)}{2}$.

We hope that $\overline {\cal M}_{\Omega, E}$ is a smooth manifold
of expected dimension if $g$ and $\Omega$ are in general position.

Let ${\rm ad}(E)$ be the adjoint bundle of $E$, i.e., the associated bundle
$P(E)\times _\rho {\rm Lie}(G)$ with $\rho$ being the adjoint representation
of $G$ in ${\rm Lie}(G)$, where $P(E)$ denotes
the principal bundle of the $G$-bundle $E$. 
For any connection $A$, define a linear
operator 
\begin{eqnarray}
\label{eq:6.2.1}
L_A: \Omega^1(M, {\rm ad}(E)) &\mapsto& \Omega^0(M,{\rm ad}(E))\oplus
\Omega_+^2(M,{\rm ad}(E)),\spn{6.2.1}\\[2ex]
L_A(\varphi) &=& 
(D_A^*\varphi, D_A\varphi + \ast (\Omega \wedge D_A\varphi)). \nonumber
\end{eqnarray}
By our assumption on $\Omega$, each $L_A$ is
elliptic. Its index is the expected dimension of 
$\overline {\cal M}_{\Omega, E}$.

If $M$ is a Calabi-Yau $4$-fold with $\theta$ and $\omega$
as in (1.3.1), then by simple computations, one can show
that the index of $L_A$ is the same as half of the index of the 
$\overline \partial$-operator $D^{0,1}_A$ on $\Omega^{0,*}(M, {\rm End}(E))$. 
The index of $D^{0,1}_A$ can be computed easily by the Atiyah-Singer 
index theorem.

Therefore, to carry out this program, 
we need to prove only transversality for
$\Omega$-asd instantons. This will be studied in a future paper.

Let us end  this section with a simple example of the above program. Let
$E$ be an ${\rm SU}(2)$-bundle over a Calabi-Yau 3-fold $V$. Let $\theta_0$
and $\omega_0$ be, respectively, a holomorphic 3-form and a K\"ahler form on $V$,
satisfying:
$$
\theta_0\wedge \overline \theta_0 = 2 \sqrt{-1} \frac{\omega_0^3}{3!}.
$$
Now let $M=V\times T$, where $T$ is a torus of complex dimension one.
We denote by $dz$ the standard flat (1,0)-form on $T$.
Put 
$$
\Omega = 4 {\rm Re}(\theta_0 \wedge dz) + \frac{1}{2} \left ( \omega_0
+ \frac{\sqrt{-1}}{2} dz\wedge d\bar z\right )^2.
$$
Then $T$-invariant solutions of the $\Omega$-asd equation on $M$
reduce to the solutions of the following equation on $V$:
\begin{equation}
\label{eq:6.2.2}
F_A^{0,2} = \overline \partial ^* f,  ~~~F_A^{1,1}\wedge \omega^2_0 =
[f, \bar f], \spn{6.2.2}
\end{equation}
where $A$ is a connection of $E$ and $f$ is a section of $\wedge^{0,3}({\rm End}(E))$. Note that
(6.2.2) is an elliptic system.
Presumably, counting solutions of (6.2.2) leads to
the so-called holomorphic
Casson invariants as studied in [DT] and by R. Thomas in his thesis.

\AuthorRefNames [CWW]

\bye